\documentclass[twoside,10pt]{article}
\usepackage{amssymb,amsmath,amsfonts,enumerate,bm}

\usepackage[dvips]{graphicx}
\usepackage{psfrag,graphicx}

\newtheorem{theorem}{Theorem}[section]
\newtheorem{lemma}[theorem]{Lemma}

\newtheorem{proposition}[theorem]{Proposition}
\newtheorem{remark}[theorem]{Remark}
 
\numberwithin{equation}{section}

\def\bproof{\noindent{\it Proof. }}
\def\eproof{\null\hfill {$\Box$}\bigskip}

\newcommand{\Er}{\mathbb{R}}

\newcommand{\Ee}{\mathbb{E}}
\newcommand{\Le}{\mathbb{L}}
\newcommand{\Me}{\mathbb{M}}
\newcommand{\en}{\mathbb{N}}

\newcommand{\Ers}{{\mathbb{R}}_{\text{sym}}^{3\times 3}}
\newcommand{\Erdev}{{\mathbb{R}}_{\text{dev}}^{3\times 3}}
\newcommand{\Erd}{{\mathbb{R}}^{3\times 3}}

\newcommand{\tr}{{\textrm{\upshape tr}}}
\newcommand{\cK}{C^{\text{Korn}}}

\newcommand{\norm}[2][{}]{\lVert#2\rVert_{{#1}}}

\newcommand{\abs}[2][{}]{\lvert#2\rvert_{#1}}

\renewcommand{\tilde}{\widetilde}

\newcommand{\dive}[1]{\operatorname{div}#1}

\def\bfI{{\rm{\bm I}}}

\def\og{\leavevmode\raise.3ex\hbox{$\scriptscriptstyle\langle\!\langle$~}}
\def\fg{\leavevmode\raise.3ex\hbox{~$\!\scriptscriptstyle\,\rangle\!\rangle$}}

\def\C{{\mathrm{C}}}
\def\D{{\mathrm{D}}}
\def\H{{\mathrm{H}}}
\def\L{{\mathrm{L}}}
\def\V{{\mathrm{V}}}
\def\X{{\mathrm{X}}}
\def\Y{{\mathrm{Y}}}
\def\W{{\mathrm{W}}}
\def\ee{{\mathrm{e}}}
\def\tra{\mathsf{T}}
\def\eqldef{\overset{\text{\tiny def}}{=}}

\def\dd{\;\!\mathrm{d}}   

\def\calA{\mathcal A}

\def\calD{\mathcal D}

\def\calG{\mathcal G}

\def\calJ{\mathcal J}
\def\calL{\mathcal L}
\def\calQ{\mathcal Q}

\def\calV{\mathcal V}
\def\calW{\mathcal W}

\begin{document}

\title{Global existence result for phase transformations
 with heat transfer in shape memory alloys{\footnote{Dedicated 
to 75th birthday of K.~Gr\"oger}}}

\author{Laetitia Paoli\thanks{Universit\'e de Lyon, LaMUSE,
23 rue Paul Michelon, 42023 Saint-Etienne Cedex 02, France 
({\tt laetitia.paoli@univ-st-etienne.fr}).}
\and
Adrien Petrov\thanks{Weierstra{\ss}-Institut f\"ur Angewandte Analysis 
und Stochastik,
Mohrenstra{\ss}e 39, 10117 Berlin, Germany ({\tt petrov@wias-berlin.de}).}}
   
\pagestyle{myheadings}

\thispagestyle{plain}
\markboth{L. Paoli, A. Petrov}
{\small{Global existence result for phase transformations
 with heat transfer in shape memory alloys}}

\date{Berlin--St Etienne, 28 April 2011} 


\maketitle

\hspace*{-0.6cm}{\bf Abstract.} 
We consider  three-dimensional models for rate-independent
processes describing materials undergoing phase transformations  
with heat transfer. The problem is formulated 
within the framework of generalized standard solids by the coupling of 
the momentum equilibrium equation and the flow rule with the heat transfer 
equation. 
 Under appropriate regularity assumptions on the initial
data, we prove the existence a global solution for this thermodynamically 
consistent system, 
by using a fixed-point argument combined with global energy estimates.
\vspace{0.3cm}

\hspace*{-0.6cm}{\bf {Key words.}}
Existence result, generalized standard materials,  heat equation,
enthalpy transformation, maximal monotone operators, 
doubly nonlinear equations,
shape-memory alloys.
\vspace{0.3cm}

\hspace*{-0.6cm}{\bf {AMS 2000 Subject Classification:}}
35K55, 74C05, 74C10, 74F05, 74N30, 80A17.

\section{Introduction}
\label{s:introduction}

Motivated by the study of shape memory alloys, we consider rate-independent 
processes describing materials undergoing phase transformations. In the 
framework of generalized standard solids due to  Halphen and Nguyen 
(see \cite{HalNgu75SMSG}), the unknowns are the displacement field $u$ and 
an internal variable $z$ and the problem is described by the momentum 
equilibrium equation combined with a flow rule for the evolution of 
the internal variable. A very powerful tool to study such problems 
is the so called energetic formulation, introduced in
\cite{MieThe04RIHM,Miel05ERIS}
and later on developed and intensively applied in 
\cite{FraMie06ERCR,MieRou06RIDP,MieRos07EURC,Miel07MTIP,MiePet07TDPT,
MiRoSt08GLRR}.
Note that coupling rate-independent processes with rate-dependent processes
makes, in general,
the problems much more difficult; see for instance 
\cite{EfeMie06RILS,MiPeMa08CSKV, BarRou08TPSS, Roub09RISS,Roub09TVEL}. 

In this paper, we are interested in coupling the rate-independent process
with the thermal process, which is not rate-independent, and viscous damping. 
The model is based on the Helmholtz free energy \(W(\ee(u),z,\nabla
z,\theta)\), depending on the 
\emph{infinitesimal strain tensor}  \(\ee(u)\eqldef \tfrac12(\nabla u{+}\nabla
u^{\tra})\) for the \emph{displacement} \(u\),  where $(\cdot)^{\tra}$ 
denotes the transpose of the tensor,
the  \emph{internal variable}  \(z\)
and the \emph{temperature} \( \theta\).
For simplicity, we will omit any dependence on the material point
$x\in \Omega$ and \(t\in [0,T]\) with \(T>0\).
We assume that $W$ can be decomposed as follows 
\begin{equation}
\label{eq:free_energy}
W(\ee(u),z,\nabla z,\theta)\eqldef
W_1(\ee(u),z,\nabla z)-W_0(\theta)+\theta W_2(\ee(u),z).
\end{equation}
This partially linearized decomposition ensures that entropy separates 
the thermal and mechanical variables (see \eqref{eq:26}). Let us emphasize 
that the last term $\theta W_2(\ee (u),z)$ allows for coupling effects between 
the temperature and the internal variable, which is motivated by the
phenomenological models for shape memory alloys presented in Section 
\ref{sec:mech-model}. Since coupling terms will appear both in the momentum 
equilibrium equation and in the inclusion describing the evolution of the 
internal variable, this setting is more general than the one presented in 
\cite{Roub10TRIV}.
We make the assumptions of small deformations. The problem is thus described 
by the following system
\begin{subequations}
\label{eq:ent_eqED}
\begin{align}
&-\dive( \sigma_{\textrm{el}}
{+}\Le\mathrm{e}(\dot u))
=\ell, \quad \sigma_{\textrm{el}}\eqldef\D_{\ee (u)} W(\ee(u),z,\nabla z, \theta),
\label{eq:ent_eq_1ED}\\&
\partial\Psi(\dot z)+\Me\dot z+ \sigma_{\textrm{in}}
\ni 0, \quad \sigma_{\textrm{in}}\eqldef\D_{z} W(\ee(u),z, \nabla z, \theta) - 
\dive \D_{\nabla z} W(\ee(u),z, \nabla z, \theta),
\label{eq:ent_eq_2ED}\\&
c(\theta)\dot\theta-\dive(\kappa(\ee(u),z,\theta)\nabla\theta)
=\Le\ee(\dot u){:}\ee(\dot u)
+
\theta\partial_{t}W_2(\ee(u),z)
+\Psi(\dot z)+
\Me \dot z{:}\dot z.\label{eq:ent_eq_3ED}
\end{align}
\end{subequations}
Here  \(\Psi\) denotes the dissipation potential. As it is
common in modeling hysteresis effect in mechanics, 
we assume that \(\Psi\) is positively homogeneous of
degree 1, i.e., \(\Psi(\gamma z)=\gamma\Psi(z)\) for all \(\gamma\geq 0\).
The viscosity tensors are denoted \(\Le\) and \(\Me\), 
\(c(\theta)\) is the \emph{heat capacity} and 
\(\kappa(\ee(u),z,\theta)\) is the  \emph{conductivity}.
As usual, \((\,\dot{}\,)\), \(\D_z^i\)
and \(\partial\) denote the time derivative
\(\tfrac{\partial}{\partial t}\), the \(i\)-th derivative  with respect 
to \(z\)
and the subdifferential in the
sense of convex analysis (for more details see \cite{Brez73OMMS}),
respectively.
Observe that \eqref{eq:ent_eq_1ED}, \eqref{eq:ent_eq_2ED} and 
\eqref{eq:ent_eq_3ED} are usually called 
the momentum equilibrium equation, the flow rule 
and the heat-transfer equation, respectively.

The paper is organized as follows. In Section \ref{sec:mech-model}, 
the thermodynamic consistency is justified and  some illustrative
examples are presented.
Then the mathematical formulation of the problem in terms of displacement,
internal variables and temperature is presented in Section \ref{s:math_form}. 
Our problem is reformulated in terms of enthalpy, which 
is a crucial ingredient to prove the existence result. 
Sections \ref{sec:const-equat}, \ref{sec:enthalpy-equation} and 
\ref{sec:existence-result} are devoted to the proof of a 
local existence result by a fixed point argument. 
More precisely, in Section \ref{sec:const-equat}, 
we consider first the system composed by 
the momentum equilibrium equation and flow rule 
for a given temperature \(\theta\) and we prove existence and 
regularity results.
Next in Section \ref{sec:enthalpy-equation}, we recall 
existence and regularity results for the enthalpy equation for any given 
right hand side.
Then a local existence result follows in Section \ref{sec:existence-result} 
by using a fixed-point argument. 
Finally a global energy estimate is established in Section 
\ref{sec:global-existence-result} leading to a global existence result
for the system \eqref{eq:ent_eqED}.

\section{Mechanical model}
\label{sec:mech-model}

We justify here 
the thermodynamic consistency of the model \eqref{eq:ent_eqED}.
Starting from the Helmholtz free energy $W$, we introduce 
the specific \emph{entropy} \(s\) via the Gibb's relation 
\begin{equation}
\label{eq:17bis}
s\eqldef -\D_{\theta} W(\ee(u),z,\nabla z,\theta),
\end{equation}
and the \emph{internal energy} 
\begin{equation}
\label{eq:17}
W_{\textrm{in}}(\ee(u),z,\nabla z,\theta)
\eqldef W(\ee(u),z,\nabla z,\theta)+\theta s.
\end{equation}
Then the \emph{entropy equation} is given by 
\begin{equation}
\label{eq:18}
\theta\dot s + \dive (j) = \xi,
\end{equation}
where $j$ is the \emph{heat flux} and $\xi$ is the \emph{dissipation rate}.
We get
\begin{equation*}
\xi = \Le\ee(\dot u){:}\ee(\dot u) +\Me \dot z{:}\dot z+\Psi(\dot z) \ge 0,
\end{equation*} 
and, assuming Fourier's law for the temperature, we have
\begin{equation*}
j= - \kappa(\ee(u),z,\theta)\nabla\theta .
\end{equation*}
We can check now that the second law of thermodynamics  is 
satisfied if $\theta>0$. Indeed, assuming that the system is thermally
isolated, 
we may  divide
\eqref{eq:18} by \(\theta\) and Green's formula yields
\begin{equation*}
\begin{aligned}
&\int_{\Omega}\dot s\dd x=
\int_{\Omega}\tfrac{\dive(\kappa(\ee(u),z,\theta)\nabla\theta)}{\theta}\dd x+
\int_{\Omega}\tfrac{
\Le\ee(\dot u){:}\ee(\dot u)
{+}\Me \dot z{:}\dot z{+}\Psi(\dot z)}{\theta}\dd x\\&=
\int_{\Omega}\tfrac{\kappa(\ee(u),z,\theta)\nabla\theta \cdot \nabla 
\theta}{\theta^2}\dd x+
\int_{\Omega}\tfrac{
\Le\ee(\dot u){:}\ee(\dot u)
{+}\Me \dot z{:}\dot z{+}\Psi(\dot z)}{\theta}\dd x \ge 0.
\end{aligned}
\end{equation*}
We differentiate now \(W_{\textrm{in}}(\ee(u),z,\nabla
z,\theta)\) with respect to time, we obtain by using the chain rule and 
\eqref{eq:17} that
\begin{equation}
\label{eq:22}
\begin{aligned}
&\dot W_{\textrm{in}}(\ee(u),z,\nabla
z,\theta)
=\D_{\ee(u)}  W(\ee(u),z,\nabla z,\theta){:}\ee(\dot u)\\&+
\D_z W(\ee(u),z,\nabla z,\theta){:}\dot z+
\D_{\nabla z} W(\ee(u),z,\nabla z,\theta){\cdot}\nabla\dot z+
\theta\dot{s}.
\end{aligned}
\end{equation}
We integrate \eqref{eq:22} over \(\Omega\), thus we use the Green's formula 
and \eqref{eq:18}, we find
\begin{equation}
\label{eq:23}
\begin{aligned}
&\int_{\Omega} \dot W_{\textrm{in}}(\ee(u),z,\nabla
z,\theta)\dd x=\int_{\Omega}
\D_{\ee(u)} W(\ee(u),z,\nabla z,\theta){:}\ee(\dot u)\dd x+
\int_{\Omega} \D_z W(\ee(u),z,\nabla z,\theta){:}\dot z\dd x\\&+
\int_{\Omega} \D_{\nabla z} W(\ee(u),z,\nabla z,\theta){\cdot}\nabla\dot z\dd x
+\int_{\Omega}(\dive(\kappa(\ee(u),z,\theta)\nabla\theta){+}
\Le\ee(\dot u){:}\ee(\dot u)
{+}\Me \dot z{:}\dot z{+}\Psi(\dot z))\dd x.
\end{aligned}
\end{equation}
On the one hand, we multiply \eqref{eq:ent_eq_1ED} by \(\dot u\),
and we integrate this expression over \(\Omega\) to get
\begin{equation}
\label{eq:24}
\int_{\Omega}\D_{\ee(u)} W(\ee(u),z,\nabla z,\theta){:}\ee(\dot u)\dd x
+\int_{\Omega}\Le\ee(\dot
u){:}\ee(\dot u)\dd x=\int_{\Omega}\ell{\cdot}\dot u\dd x.
\end{equation}
On the other hand, the definition of the subdifferential 
\(\partial\Psi(\dot z)\) leads to the variational equality
associated to \eqref{eq:ent_eq_2ED} 
\begin{equation}
\label{eq:25}
\int_{\Omega}\D_z W(\ee(u),z,\nabla z,\theta){:}\dot z\dd x
+ \int_{\Omega} \D_{\nabla z} W(\ee(u),z,\nabla z,\theta){\cdot}\nabla\dot z\dd x
+\int_{\Omega}\Me\dot z{:}\dot z\dd x
+ \int_{\Omega} \Psi(\dot z) \dd x =0.
\end{equation}
We use \eqref{eq:24} and \eqref{eq:25} into \eqref{eq:23}, we obtain
\begin{equation*}
\int_{\Omega} \dot W_{\textrm{in}}(\ee(u),z,\nabla
z,\theta)\dd x=
\int_{\Omega}\ell{\cdot}\dot u\dd x
+ \int_{\partial\Omega} \kappa(\ee(u),z,\theta)\nabla\theta{\cdot}\eta \dd x.
\end{equation*}
This means that the total energy balance can be expressed in terms of the
internal energy, which is the sum of power of external load and heat.

Note that from \eqref{eq:17bis}, we get
\begin{equation}
\label{eq:26}
s\eqldef
\D_{\theta} W_0(\theta)-W_2(\ee(u),z),
\end{equation}
and
\begin{equation}
\label{eq:27}
W_{\textrm{in}}(\ee(u),z,\nabla z,\theta)
\eqldef W_1(\ee(u),z,\nabla z,\theta)+\theta \D_{\theta}W_0(\theta)-W_0(\theta).
\end{equation}
We use \eqref{eq:27} into \eqref{eq:18}, we may deduce 
the heat-transfer equation \eqref{eq:ent_eq_3ED} with the heat capacity
given by \(c(\theta)=\theta\D_{\theta}^2 W_0(\theta)\).

\bigskip
Motivated by the study of shape memory alloys, we will focus in the 
rest of the paper on the special case where  $W_1$ is given by
\begin{equation*}
W_1 (\ee(u),z,\nabla z)\eqldef\tfrac12\Ee(\mathrm{e}(u){-}z){:}
(\mathrm{e}(u){-}z)+
\tfrac{\nu}2\abs{\nabla z}^2+H_1(z),
\end{equation*} 
where the internal variable \(z\) is a deviatoric $d \times d$ tensor, 
$H_1$ is a 
hardening functional and the term $\tfrac{\nu}2\abs{\nabla z}^2$, $\nu>0$, 
takes into account 
some nonlocal interaction effect for the internal variable.
Moreover we will assume 
\begin{equation*}
W_2 (\ee(u) , z)\eqldef\alpha \tr(\ee(u))+H_2(z).
\end{equation*}
Here $\alpha {\rm I}$, with $\alpha \ge 0$ and $ {\rm I}$ the identity 
matrix, is the isothermal expansion tensor.
Furthermore, the sum $H(z, \theta)\eqldef H_1(z)+\theta H_2(z)$ may be 
interpreted 
as a temperature-dependent hardening functional. 
Note that the system \eqref{eq:ent_eqED} is then rewritten as
\begin{equation*}
\begin{aligned}
&-\dive(\Ee(\mathrm{e}(u){-}z){+}\alpha\theta\bfI{+}\Le\mathrm{e}(\dot
u))
=\ell,\\&
\partial\Psi(\dot z)+\Me\dot z-\Ee(\ee(u){-}z)+\D_z H_1(z)+
\theta\D_zH_2(z)-\nu\Delta z\ni 0,\\&
c(\theta)\dot\theta-\dive(\kappa(\ee(u),z,\theta)\nabla\theta)
=\Le\ee(\dot u){:}\ee(\dot u)
+\theta(\alpha\tr(\ee(\dot u)){+}\D_zH_2(z){:}\dot z)+\Psi(\dot z)+
\Me \dot z{:}\dot z.
\end{aligned}
\end{equation*}
We may illustrate the presented model by a  nontrivial example, namely,
a three-dimensional macroscopic  phenomenological  model
for shape-memory polycrystalline materials undergoing phase
transformations. This model was introduced by Souza et al
(\cite{SoMaZo98TDMS}) and 
by  Auricchio et al (see \cite{AurPet02IACR}).  
The internal variable describes the inelastic part of the deformation 
due to the martensitic phase transformation and the
hardening functional $H$ takes the form
\begin{equation*}
H(z, \theta)\eqldef c_1 (\theta) \abs{z}+
c_2 (\theta) \abs{z}^2+\chi(z).
\end{equation*}
Here \(\chi:\Er^{d\times d}_{\text{dev}}\rightarrow [0,+\infty]\) 
denotes the indicator function of the
ball \(\{z\in\Er^{d\times d}_{\text{dev}}:\ \abs{z}\leq c_{3} (\theta)\}\)
and the coefficients $c_i(\theta) $ are positive real numbers.
Let us observe that  \(c_1 (\theta) >0\) is an activation threshold for
initiation of martensitic phase transformations,  \(c_2 (\theta)\) measures the
occurrence of hardening with respect to the internal variable
\(z\) and  \(c_3 (\theta) \) represents the maximum modulus of transformation
strain that can be obtained by alignment of martensitic variants.

The dependence of the coefficients $c_i$, $i=1,2,3$, with respect to 
$\theta$ is due to a strong thermo-mechanical constitutive coupling coming from
the latent heat absorption or release, which is one of the main features of  
the behavior of shape-memory alloys (see 
\cite{AurPet04STPT} for more details).
 
For mathematical purposes, we should regularize the hardening functional as in 
 \cite{MiePet07TDPT}. Namely, we replace \(H(z, \theta)\) by 
\begin{equation*} 
H^{\delta}(z,\theta)\eqldef
c_1 (\theta) \sqrt{\delta^2{+}\abs{z}^2}
+ c_2 (\theta )\abs{z}^2+
\tfrac{((\abs{z}{-}c_3 (\theta))_+)^4}{\delta (1{+}\abs{z}^2)},
\end{equation*}
where \(\delta>0\) is a small parameter.
Then, assuming that the mappings $c_i$, $i=1,2,3$, are of class $\C^1$ we 
can consider an affine approximation of  $H^{\delta}(z, \theta)$ as 
$H^{\delta}_1 (z) + \theta H^{\delta}_2 (z)$  with
\begin{equation*}
H_1^{\delta}(z)=
\bar  c_1  \sqrt{\delta^2{+}\abs{z}^2}
+ \bar c_2 \abs{z}^2+
\tfrac{((\abs{z}{-} \bar c_3 )_+)^4}{\delta (1{+}\abs{z}^2)},
\end{equation*}
where \(\bar c_i >0\), \(i=1,2,3\).

Let us emphasize that the existence proof presented in the next sections 
can be easily extended to more general models of shape memory alloys for 
which $W_1$ is given by
\begin{equation*}
W_1(\ee(u),z,\nabla z)\eqldef 
\tfrac{1}{2} \Ee( \ee(u){-}E(z)){:}(\ee(u){-} E(z)) + 
\tfrac{\nu}2\abs{\nabla z}^2+H_1(z),
\end{equation*} 
where the internal variable \(z\) is a vector of 
$\Er^{N-1}$, with $N \ge 2$, and $E$ is an affine mapping from 
$\Er^{N-1}$ to the set of $d \times d$ deviatoric tensors. In such 
models, $N$ is the total number of phases, i.e., the austenite and all the 
variants of martensite, and the components $z_1, \dots, z_{N-1}$ of $z$ 
and $z_N \eqldef 1-\sum_{k=1}^{N-1} z_k$ are interpreted as phase fractions. 
Then $E(z)$ is the effective transformation strain of the mixture, given by
\begin{equation*}
E(z)\eqldef \sum_{k=1}^{N-1} z_k E_k + \Bigl(1{-}\sum_{k=1}^{N-1} z_k\Bigr) E_N,
\end{equation*}
where $E_k$ is the transformation strain of the phase $k$ and 
the temperature dependent hardening functional $H(z, \theta)$ is 
the sum of a smooth part $w(z, \theta)$ and the indicator function of 
the set $[0,1]^{N-1}$ 
(see \cite{MieThe99MMRI,Miel00EMFM,HalGov02ARFE,GoMiHa02FEMV,MiThLe02VFRI,
GoHaHe07UBFE}).
The hardening functional can be regularized in a similar 
way as in the previous example by replacing $H(z, \theta)$ by 
\begin{equation*}
H^{\delta}(z,\theta) = 
w(z, \theta) + \sum_{k=1}^{N-1} 
\tfrac{((-z_k)_+)^4{+} 
((z_k{-}1)_+)^4}{\delta (1{+}\abs{z_k}^2)},
\end{equation*}
and we may consider an affine approximation of $w(z, \theta)$ as 
$w_1(z) + \theta w_2(z)$. For more details on this example the 
reader is referred to \cite{PaoPet11TMPV}.

\section{Mathematical formulation}
\label{s:math_form}

We consider a reference configuration \(\Omega\subset\Er^3\).
We assume that \(\Omega\) is an  bounded domain such that $\partial 
\Omega$ is of class $\C^{2+\rho}$.
We will denote by \(\Ers\) the space of symmetric 
\(3{\times}3\) tensors endowed with the natural scalar product
\(v{:}w\eqldef \tr(v^{\tra}w)\) 
and the corresponding norm \(\abs{v}^2\eqldef v{:}v\) for all
\(v,w\in\Ers\). In particular, we assume that
\begin{equation*}
W_1(\ee(u),z,\nabla z)\eqldef
\tfrac12\Ee(\mathrm{e}(u){-}z){:}(\mathrm{e}(u){-}z)+
\tfrac{\nu}2\abs{\nabla z}^2+H_1(z)
\quad{\text{ and }}\quad
W_2(\ee(u),z)\eqldef \alpha\tr(\ee(u))+H_2(z),
\end{equation*}
where \( \nu>0\), 
\(\alpha \ge 0\), 
is the isotropic thermal expansion coefficient,
\(\Ee\) denotes the \emph{elastic tensor} and
\(H_i\), \(i=1,2\), two \emph{hardening functionals}.
Given a function \(\ell:\Omega{\times} (0,T)\rightarrow
\Er^3\), we look for 
a \emph{displacement} \(u:\Omega{\times} (0,T)\rightarrow \Er^3\), 
a matrix of  \emph{internal variables}  
\(z:\Omega{\times} (0,T)\rightarrow \Erdev\)
and a
\emph{temperature} \(\theta:\Omega{\times} (0,T)\rightarrow \Er\)
satisfying the following system:
\begin{subequations}
\label{eq:ent_eq}
\begin{align}
&
-\dive(\Ee(\mathrm{e}(u){-}z){+}\alpha\theta\bfI{+}\Le\mathrm{e}(\dot
u))
=\ell,\label{eq:ent_eq_1}\\&
\partial\Psi(\dot z)+\Me\dot z-\Ee(\ee(u){-}z)+\D_z H_1(z)+
\theta\D_zH_2(z)-\nu\Delta z\ni 0,\label{eq:ent_eq_2}\\&
c(\theta)\dot\theta-\dive(\kappa(\ee(u),z,\theta)\nabla\theta)
=\Le\ee(\dot u){:}\ee(\dot u)
+\theta(\alpha\tr(\ee(\dot u)){+}\D_zH_2(z){:}\dot z)+\Psi(\dot z)+
\Me \dot z{:}\dot z.\label{eq:ent_eq_3}
\end{align}
\end{subequations}
We have naturally to prescribe initial conditions for
the displacement, the internal variables, and the temperature, namely
\begin{equation}
\label{eq:init_cond}
u(\cdot,0)=u^0,\quad
z(\cdot,0)=z^0,\quad
\theta(\cdot,0)=\theta^0.
\end{equation}
The problem is to be completed with boundary conditions. More precisely,
we suppose here that
\begin{equation}
\label{eq:boun_cond}
u_{|_{\partial\Omega}}=0,\quad
\nabla z{\cdot}\eta_{|_{\partial\Omega}}=0,\quad
\kappa\nabla\theta{\cdot}\eta_{|_{\partial\Omega}}=0,
\end{equation}
where \(\eta\) denotes the outward normal to the boundary \(\partial\Omega\)
of \(\Omega\). 
The original problem \eqref{eq:ent_eq} 
can be rewritten in terms of enthalpy instead of temperature by employing
the so-called enthalpy transformation
\begin{equation}
\label{eq:2}
g(\theta)=\vartheta\eqldef \int_0^{\theta}c(s)\dd s.
\end{equation}
Clearly, 
$g$ is the unique primitive of 
the function \(c\), which is supposed to be continuous,  such that
\(g(0)=0\). Furthermore, we will assume that for all \(s\geq 0\), 
\(c(s) \ge c^c >0\) where \(c^c\) is a constant. Hence we deduce that
\(g\) is a bijection from \([0,\infty)\) into \([0,\infty)\).
We define 
\begin{subequations}
\label{eq:34}
\begin{align}
\label{eq:1}
&\zeta(\vartheta)\eqldef 
\begin{cases}
g^{-1}(\vartheta) &\text{ if } \vartheta\geq 0,\\
0 &\text{otherwise},
\end{cases}\\&
\label{eq:37}
\kappa^c(\ee(u),z,\vartheta)
\eqldef\tfrac{\kappa(\ee(u),z,\zeta(\vartheta))}{c(\zeta(\vartheta))},
\end{align}
\end{subequations}
where \(g^{-1}\) is the inverse of \(g\). For more details on the 
enthalpy transformation, the reader is referred to \cite{Roub09TVEL} and the
references therein. Therefore 
the system \eqref{eq:ent_eq} is transformed into the following form 
\begin{subequations}
\label{eq:ent_eq1}
\begin{align}
\label{eq:ent_eq1_11}
&-\dive(\Ee(\mathrm{e}(u){-}z){+}\alpha\zeta(\vartheta)\bfI{+}
\Le\mathrm{e}(\dot
u))
=\ell,\\&
\label{eq:ent_eq1_12}
\partial\Psi(\dot z)+\Me\dot z-\Ee(\ee(u){-}z)+\D_z H_1(z)+
\zeta(\vartheta)\D_zH_2(z)-\nu\Delta z\ni 0,\\&
\label{eq:ent_eq1_13}
\dot\vartheta-\dive(\kappa^c(\ee(u),z,\vartheta)\nabla\vartheta)
=\Le\ee(\dot u){:}\ee(\dot u)
+\zeta(\vartheta)(\alpha\tr(\ee(\dot u)){+}\D_zH_2(z){:}\dot z) 
+\Psi(\dot z)+\Me \dot z{:}\dot z,
\end{align}
\end{subequations}
with boundary conditions
\begin{equation}
\label{eq:boun_cond1}
u_{|_{\partial\Omega}}=0,\quad
\nabla z{\cdot}\eta_{|_{\partial\Omega}}=0,\quad
\kappa^c\nabla\vartheta{\cdot}\eta_{|_{\partial\Omega}}=0,
\end{equation}
and initial conditions 
\begin{equation}
\label{eq:init_cond1}
u(\cdot,0)=u^0,\quad
z(\cdot,0)=z^0,\quad
\vartheta(\cdot,0)=\vartheta^0=  g(\theta^0).
\end{equation}
The identity \eqref{eq:ent_eq1_13} is called the enthalpy equation.
As usual  Korn's inequality will play a role in the
mathematical analysis developed in the next sections. 
We have assumed that $\partial \Omega$ is of class $\C^{2+\rho}$, 
so that the Korn's
inequality holds, i.e. 
\begin{equation}
\label{eq:58} 
\exists \cK>0\ \forall u\in \H_0^1(\Omega)  :\
\norm[\L^2(\Omega)]{\ee(u)}^2
\geq
\cK
\norm[\H^1(\Omega)]{u}^2,
\end{equation}
(see \cite{KonOle88BVPS,DuvLio76IMP}). 

Let us introduce now  the assumptions  on the dissipation potential
\(\Psi\), on the hardening functions \(H_i\), \(i=1,2\),
and on the data
\(\Ee\), \(\Le\), \(\Me\), \(\ell\), \(c\eqldef c(\theta)\) and 
\(\kappa\eqldef \kappa(\ee(u),z,\theta)\), 
which will allow us to obtain some regularity properties that are needed 
to prove the existence result. 

We assume that the dissipation potential $\Psi$ is positively
homogeneous of degree $1$
and satisfies the
triangle inequality, namely, we have
\begin{subequations}
\label{eq:Psi} 
\begin{align}
\label{hom}
&\forall \gamma \ge 0\ \forall z\in \Ers: \ 
\Psi(\gamma z)=\gamma\Psi(z),
\\
\label{eq:Psi.bdd}
&\exists C^\Psi>0 \ \forall z\in \Ers:\ 
0 \le \Psi(z)\leq C^\Psi \abs{z},
\\
\label{eq:4}
&\forall z_1,z_2\in \Ers:\ \Psi(z_1{+}
z_2)\leq \Psi(z_1)+\Psi(z_2).
\end{align}
\end{subequations} 
It is clear that \eqref{hom}, \eqref{eq:Psi.bdd} and \eqref{eq:4} 
imply that $\Psi$ is convex and continuous.
We impose that the hardening functionals \(H_i\), 
 \(i=1, 2\),
 belong to  $\C^2 (\Erdev; \Er)$ and satisfy the following inequalities
\begin{subequations}
\label{eq:H}
\begin{align}
\label{eq:H1}
&\exists c^{H_1}, \tilde c^{H_1} >0  \ \forall z\in\Erdev :\ 
H_1(z)\geq c^{H_1}{\abs{z}}^2 - \tilde c^{H_1},
\\
\label{eq:H3}
&
\exists C^{H_i}_{zz}>0\ \forall z\in\Erdev:\ 
{\abs{\textrm{D}_{z}^2 H_i(z)}}\leq C^{H_i}_{zz}.
\end{align}
\end{subequations}
Note that \eqref{eq:H3} leads to
\begin{equation}
\label{eq:H2}
\exists C^{H_i}_{z}>0\ \forall z\in\Erdev :\ 
{\abs{\textrm{D}_{z} H_i(z)}}\leq C^{H_i}_{z}(1{+}{\abs{z}}), \quad 
{\abs{H_i(z)}} \le  C^{H_i}_{z}(1{+}{\abs{z}}^2).
\end{equation}  
The \emph{elastic tensor} \(\Ee:\Omega\rightarrow \calL(\Ers,\Ers)\)
is a symmetric  positive definite operator such that
\begin{subequations}
\label{eq:47}
\begin{align}
\label{eq:3}
&
\exists c^{\Ee}>0\ \forall z\in\L^2(\Omega; \Ers):\ 
c^{\Ee}\norm[\L^2(\Omega)]{z}^2
\leq
\int_{\Omega}\Ee z{:}z\dd x,\\&
\label{eq:46}
\forall i,j,k=1,2,3:\
\Ee(\cdot),
\tfrac{\partial\Ee_{i,j}(\cdot)}{\partial x_k}\in
\L^{\infty}(\Omega).
\end{align}
\end{subequations}
We assume that \(\Le\) and \(\Me\) 
are symmetric  positive definite tensors. 
This implies that
\begin{subequations}
\label{eq:LM} 
\begin{align}
\label{eq:6}
&
\exists c^{\Le}, C^{\Le}>0\ \forall z\in\Er^{3\times 3}:\ 
c^{\Le}\abs{z}^2
\leq
\Le z{:} z\leq
C^{\Le}\abs{z}^2,
\\
\label{eq:8}
&
\exists c^{\Me}, C^{\Me}>0\ \forall z\in\Er^{3\times 3}:\ 
c^{\Me}\abs{z}^2\leq
\Me z{:}z\leq
C^{\Me}\abs{z}^2.
\end{align}
\end{subequations}
We assume furthermore that
\begin{equation} \label{eq:8bis}
\forall z \in \Erdev: \ \Me z \in \Erdev.
\end{equation}
We consider that \(\ell\) is an external loading satisfying
\begin{equation}
\label{eq:ell}
\ell\in \H^1(0,T;\L^2(\Omega)).
\end{equation}
Finally, for  the heat capacity
\(c\) and the conductivity \(\kappa^c\) we assume that
\begin{subequations}
\label{eq:36}
\begin{align}
\label{eq:38}
&c: [0,\infty)\rightarrow [0,\infty) \text{ is continuous},\\&
\label{eq:40}
\exists 
\beta_1
 \ge 2 \ \exists  c^c>0\ \forall \theta\geq 0:\
c^c(1{+}\theta)^{\beta_1{-}1}\leq c(\theta) ,
\\&
\label{eq:39}
\kappa^c: \Ers{\times}\Erd_{\textrm{dev}}{\times}\Er\rightarrow \Ers
\text{ is continuous, bounded and uniformly positive definite, i.e}\\&
\label{eq:41}
\exists c^{\kappa^c}>0\ \forall (e, z, \vartheta) \in \Ers \times \Erdev
\times \Er
\ \forall v\in\Er^3:\ \kappa^c (e,z, \vartheta) v{\cdot}v\geq 
c^{\kappa^c}  \abs{v}^2, \\&
\exists  C^{\kappa^c}>0\ \forall (e, z, \vartheta) \in \Ers \times \Erdev
\times \Er :\ \abs{\kappa^c (e,z, \vartheta)} \leq 
C^{\kappa^c} .
\end{align}
\end{subequations}
Let us emphasize that the three dimensional model for shape memory alloys
presented in Section \ref{sec:mech-model} fulfills the previous assumptions 
so that we may apply  
the abstract result obtained in the next sections.

Let us end this section by some comments about the proof strategy.
In order to obtain a local existence result for the coupled problem 
\eqref{eq:ent_eq1}--\eqref{eq:init_cond1}, 
a fixed point argument will be used. 
More precisely, for any given $\tilde \vartheta$, we
define 
$\theta\eqldef\zeta (\tilde \vartheta)$ and we solve first the system
composed by the momentum equilibrium equation and the flow rule
\eqref{eq:ent_eq_1}--\eqref{eq:ent_eq_2}, then we solve the enthalpy equation 
\eqref{eq:ent_eq1_13} with 
$\kappa^c\eqldef \kappa^c(\ee(u), z, \zeta (\tilde \vartheta))$. 
This allows us to define a mapping
\begin{equation*}
\phi: \tilde \vartheta \mapsto \vartheta, 
\end{equation*}
and our aim is to prove that this mapping satisfies the assumptions of
Schauder's fixed point theorem. Therefore, let us
consider a 
given $\tilde \vartheta \in \L^{\bar q} (0,T; \L^{\bar p} (\Omega))$ with 
$\bar p \ge 1$ and $\bar q\ge 1$. We define $\theta\eqldef\zeta (\tilde
\vartheta)$. Since 
$\zeta$ is a Lipschitz continuous mapping from $\Er$ to $\Er$, we infer that
the mapping 
\begin{equation*}
\begin{aligned}
\phi_1:\ 
\L^{\bar q} (0,T; \L^{\bar p} (\Omega)) &\rightarrow \L^{\bar q}(0,T; 
\L^{\bar p} (\Omega)) \\
 \tilde \vartheta &\mapsto \theta 
\end{aligned}
\end{equation*}
is also Lipschitz continuous. Furthermore \eqref{eq:40} implies that 
\begin{equation*}
\forall \theta \in [0,\infty):\
\tfrac{c^c}{\beta_1} ((1{+}\theta)^{\beta_1}{-}1)\eqldef g_1(\theta) \leq
g(\theta) .
\end{equation*}
Thus we have
\begin{equation*}
\forall \vartheta \in [0,\infty):\ 0
\leq \zeta (\vartheta) \leq \zeta_1 (\vartheta)\eqldef g_1^{-1}(\vartheta),
\end{equation*}
and 
\begin{equation*} 
\forall \vartheta \in \Er:\
\abs{\zeta(\vartheta)}\leq \bigl(\tfrac{\beta_1}{c^c} 
\vartheta^+ {+}1 \bigr)^{\frac1{\beta_1}}-1,
\end{equation*}
with $\vartheta^+\eqldef\max (\vartheta, 0)$ 
for all $\vartheta \in \Er$. It follows that
\begin{equation} \label{eqlp:1}
\forall \beta \in [1, \beta_1] \ \forall \vartheta \in \Er:\
\abs{\zeta(\vartheta)}\leq \bigl(\tfrac{\beta_1}{c^c} 
\vartheta^+ {+}1 \bigr)^{\frac1{\beta}} -1 
\leq \bigl(\tfrac{\beta_1}{c^c} 
\vartheta^+ \bigr)^{\frac1{\beta}}.
\end{equation}
Hence, for all $\beta \in [1, \beta_1]$ and for all
$\tilde \vartheta \in \L^{\bar q} (0,T; \L^{\bar p} (\Omega))$, 
we have $\theta = \zeta(\tilde \vartheta) \in \L^{\beta \bar q} 
(0,T; \L^{\beta \bar p} (\Omega))$ with
\begin{equation*}
\norm[\L^{\beta \bar q} (0,T; \L^{\beta \bar p} (\Omega))]{\theta} 
\leq \bigl(\tfrac{\beta_1}{c^c}\bigr)^{\frac1{\beta}} 
\norm[\L^{ \bar q}(0,T; \L^{ \bar p} (\Omega))]{\tilde \vartheta}^{\frac1{\beta}}.
\end{equation*}

In the rest of the paper, we will assume that $\bar q >4$
and $\bar p =2$. 

When there is not any confusion, 
we will use simply the notation \(\X(\Omega)\) instead of  \(\X(\Omega;\Y)\) 
where \(\X\) is a functional space and \(\Y\) is a vectorial space. 

\section{Existence and regularity results for
the system composed by the momentum equilibrium 
equation and the flow rule}
\label{sec:const-equat}

This section is devoted to the proof of
existence and uniqueness results 
for the system composed by the momentum equilibrium 
equation and the flow rule \eqref{eq:ent_eq_1}--\eqref{eq:ent_eq_2}
under the consideration that \(\theta=\zeta(\tilde \vartheta)\) is given  
in a bounded subset of  
\(\L^q(0,T;\L^p(\Omega))\) with $q = \beta_1 \bar q$ and 
$p\in \bigl[4,  \min( \beta_1 \bar p, 6) \bigr] $. More precisely, 
we look for a solution of the problem $(\textrm{P}_{uz})$:
\begin{subequations}
\label{eq:ent_eq_1_1}
\begin{align}
&-\dive(\Ee(\mathrm{e}(u){-}z){+}\alpha\theta\bfI{+}\Le\mathrm{e}(\dot
u))
=\ell,\label{eq:ent_eq_1_2}\\&
\partial\Psi(\dot z)+\Me\dot z-\Ee(\ee(u){-}z)+\D_z H_1(z)+
\theta\D_zH_2(z)-\nu\Delta z\ni 0,\label{eq:ent_eq_2_2}
\end{align}
\end{subequations}
with initial conditions
\begin{equation}
\label{eq:init_cond2}
u(\cdot,0)=u^0,\quad
z(\cdot,0)=z^0,
\end{equation}
and boundary conditions
\begin{equation}
\label{eq:boun_cond2}
u_{|_{\partial\Omega}}=0,\quad
\nabla z{\cdot}\eta_{|_{\partial\Omega}}=0.
\end{equation}
Furthermore we will establish some a priori estimates and some 
regularity results for the solution of $(\textrm{P}_{uz})$.

As a first step, we use classical results for Partial Differential Equations
(PDE) and Ordinary Differential Equations (ODE) to 
obtain an existence result. 
Let $\calA: \H^1(\Omega) \to (\H^1(\Omega))'$ be the linear continuous 
mapping defined by
\begin{equation*}
\forall (u,v)\in (\H^1(\Omega))^2:\
\langle \calA u,  v \rangle_{(\H^1(\Omega))', \H^1(\Omega)}\eqldef\int_{\Omega} 
\nu \Me^{-1} \nabla u{:}\nabla v \dd x.
\end{equation*}
Classical results about elliptic operators implies that $\calA$ generates an
analytic semigroup on $\L^2(\Omega)$, which extends to a $\C^0$-semigroup of
contractions on $\L^p(\Omega)$. We denote by $\calA_p$ 
(resp. $\calA_{\frac{p}{2}}$) 
the realization of its
generator in $\L^p(\Omega)$ (resp. $\L^{p/2}(\Omega)$) and by 
$\X_{q,p} (\Omega)$ the intersection of
interpolation spaces 
\begin{equation*}
\X_{q,p}(\Omega) \eqldef(\L^p(\Omega), \calD(\calA_p))_{1-
  \frac{2}{q}, \frac{q}{2}} \cap  (\L^{p/2} (\Omega), \calD(\calA_{\frac{p}2}))_{1-
  \frac{1}{q}, q},
\end{equation*}
(see for instance \cite{HieReh08QPMB, PruSch01SMRP} 
and the references therein). Here 
\(\calD(\calA_p)\) (resp. \(\calD(\calA_{\frac{p}2})\)) 
denotes the domain of \(\calA_p\) (resp. \(\calA_{\frac{p}2}\)).

In the sequel, the notations for the constants introduced in the
proofs are valid only in the proof and we also use the set 
\(\calQ_{\tau}\eqldef \Omega\times (0,\tau)\) with \(\tau\in[0,T]\). 

\begin{theorem}[Existence  for $(\textrm{P}_{uz})$]
\label{sec:existence}
Let $\theta$ 
be given in $\L^q( 0,T; \L^p 
(\Omega))$. 
Assume that $u^0 \in \H^1(\Omega)$, $z^0 \in \H^1(\Omega)$ and that
\eqref{eq:Psi}, \eqref{eq:H}, \eqref{eq:LM} 
and \eqref{eq:ell} 
hold.
Then the problem 
\eqref{eq:ent_eq_1_1}--\eqref{eq:boun_cond2}
admits a solution 
\((u,z)\in \H^1(0,T;\H^1_0(\Omega)\times\L^2(\Omega))
\cap\L^{\infty}(0,T;\H^1_0(\Omega){\times}\H^1(\Omega))\).
Furthermore if $z^0 \in \X_{q,p} ( \Omega)$ 
then  $z \in \L^q (0,T; \H^2 (\Omega)) \cap \C^0([0,T], 
\H^1(\Omega))$ and $\dot z \in 
\L^q(0,T; \L^2(\Omega))$.
\end{theorem}

\bproof
Observe first that for all \( f \in\L^2(0,T;
(\H_0^1(\Omega);\Er^3)')\) and for all \(u^* \in \H_0^1(\Omega;\Er^3) \)
the following problem 
\begin{equation*}
- \dive(\Ee\ee(u){+}\Le\ee(\dot u))= f ,
\end{equation*}
with initial conditions
\begin{equation*}
u(\cdot,0)=u^* \in \H_0^1(\Omega;\Er^3),
\end{equation*}
and boundary conditions
\begin{equation*}
u_{|_{\partial\Omega}}=0 ,
\end{equation*}
possesses a unique solution 
\(u\eqldef\calL(u^*,f) \in\H^1(0,T;\H_0^1(\Omega))\cap \C^0([0,T];
\H_0^1(\Omega))\).
Note that 
\begin{equation*}
\forall f_1, f_2  \in\L^2(0,T;(\H_0^1(\Omega);\Er^3)')\
\forall u^* \in \H_0^1(\Omega;\Er^3):\
\calL (u^*, f_1{+}f_2) = \calL (u^*,f_1) + \calL(0,f_2),
\end{equation*}
and the mapping 
\begin{equation*}
\begin{aligned}
\calL_0 (\cdot)\eqldef\calL(0, \cdot):\ \L^2(0,T;
(\H_0^1(\Omega)'))&\rightarrow  
\H^1(0,T;\H_0^1(\Omega))\cap \C^0([0,T];\H_0^1(\Omega)),\\
 f &\mapsto u,
\end{aligned}
\end{equation*}
is linear and continuous. Moreover, the classical energy estimate gives
\begin{equation} 
\label{eqlp:3}
\forall t \in [0,T]:\
c^{\Ee}\norm[\L^2(\Omega)]{\ee (u(\cdot,t))}^2  + c^{\Le}\cK \int_0^t 
\norm[\H^1(\Omega)]{\dot u(\cdot,s)}^2 \dd s
\le \tfrac{1}{c^{\Le}\cK} \int_0^t \norm[(\H^1_0(\Omega))']{f(\cdot,s)}^2  
\dd s,
\end{equation}
for all $f \in \L^2(0,T;(\H_0^1(\Omega)'))$ and $u = \calL_0 (f)$.
It follows that \eqref{eq:ent_eq_1_2} and \eqref{eq:ent_eq_2_2}
can be rewritten as 
\begin{equation}
\label{eq:5}
\partial\Psi(\dot z)+\Me\dot z+\Ee z+\D_z H_1(z)+
\theta\D_zH_2(z)-\nu\Delta z+g_1(\theta)+g_2(z)\ni 0,
\end{equation}
with initial conditions
\begin{equation}
\label{eq:7}
z(\cdot,0)=z^0 \in \H^1(\Omega ;\Erdev),
\end{equation}
and boundary conditions
\begin{equation}
\label{eq:9}
\nabla z{\cdot}\eta_{|_{\partial\Omega}}=0,
\end{equation}
where 
\(g_1(\theta)\eqldef -\Ee\ee(\calL(u^0, \ell))  
-\Ee\ee(\calL_0( \dive(\alpha\theta\bfI)))\in
\H^1(0,T;\L^2(\Omega;\Ers))\)
and 
\begin{equation*}
\begin{aligned}
g_2:\ \L^2(0,T;\L^2(\Omega;\Erdev)) &\rightarrow \H^1(0,T;\L^2(\Omega;\Ers)),\\
z&\mapsto \Ee\ee(\calL_0(\dive(\Ee z))).
\end{aligned}
\end{equation*}
Let \(\varphi:\L^2(\Omega;\Erdev)
\rightarrow (- \infty,+\infty]\) defined by
\begin{equation*}
z\mapsto 
\begin{cases}
\tfrac{\nu}2\norm[\L^2(\Omega)]{\nabla z}^2 \ \text{ if } \ z\in
\H^1(\Omega;\Erdev),\\
+\infty\ \text{ otherwise},
\end{cases}
\end{equation*}
is a proper, convex lower semicontinuous function.
Clearly, \( \partial\varphi\) is a maximal monotone operator,
for more details on the maximal monotone operators and their properties,
the reader is referred to \cite{Brez73OMMS}.
The resolvant of the subdifferential \(\partial\varphi\) is defined by 
\begin{equation*}
\forall\epsilon>0:\
\calJ_{\epsilon}\eqldef (\bfI{+}\epsilon\partial\varphi)^{-1}.
\end{equation*}
Note that the resolvant \(\calJ_{\epsilon}\)
is a contraction defined on all \(\L^2(\Omega;\Erdev)\). 
Let us introduce also
\begin{equation*} 
\forall \epsilon>0 \ \forall z\in\L^2(\Omega;\Erdev):\
\varphi_{\epsilon} (z)
\eqldef
\min_{\bar z\in\L^2(\Omega;\Erdev)}
\bigl{\{}
\tfrac1{2\varepsilon}\norm[\L^2(\Omega)]{z{-}\bar z}^2{+}\varphi(\bar z)
\bigr{\}},
\end{equation*}
which, is a convex and Fr\'echet differentiable mapping from  
\(\L^2(\Omega;\Erdev)\) to $\Er$.
Using \cite[Prop.~2.11]{Brez73OMMS}, we know that the Yosida 
approximation of $\partial \varphi$ coincide with the Fr\'echet 
differential of $\varphi_{\epsilon}$, i.e., we have
\begin{equation}
\label{eq:101}
\forall\epsilon>0:\
\partial\varphi_{\epsilon}\eqldef
\tfrac1{\epsilon}(\bfI{-}\calJ_{\epsilon}),
\end{equation}
and $\partial \varphi_{\epsilon}$ is a monotone and 
$\tfrac1{\epsilon}$-Lipschitz continuous mapping on \(\L^2(\Omega;\Erdev)\).

We approximate now the problem \eqref{eq:5}--\eqref{eq:9} by 
\begin{equation}
\label{eq:57}
\partial\Psi(\dot z_{\epsilon})+\Me\dot z_{\epsilon}+
\partial\varphi_{\epsilon}(z_{\epsilon})
+\Ee z_{\epsilon}+\D_z H_1(z_{\epsilon})
\ni -(g_1(\theta_{\epsilon}){+}g_2(\calJ_{\epsilon} z_{\epsilon}){+}\theta_{\epsilon}
\D_z H_2(\calJ_{\epsilon} z_{\epsilon})),
\end{equation}
with the initial condition
\begin{equation}
\label{eq:87}
z_{\epsilon}(0,\cdot)=z^0.
\end{equation}
Here \(\theta_{\epsilon} \in \C_0^{\infty}(0,T)\otimes\C_0^{\infty}(\Omega)\) 
and $\partial\Psi(\dot z_{\epsilon})$ is taken in the sense of the 
$\L^2( \Omega)$-extension of the subdifferential of the convex function 
$\Psi$ (see examples 2.1.3 and 2.3.3 in \cite{Brez73OMMS}), 
i.e., \eqref{eq:57} is equivalent to 
\begin{equation*}
\begin{aligned}
&-\Me\dot z_{\epsilon} (t,x) -
\partial\varphi_{\epsilon}(z_{\epsilon} (t,x) )
 -\Ee z_{\epsilon} (t,x) - \D_z H_1(z_{\epsilon} (t,x) )- 
g_1(\theta_{\epsilon}(t,x))
\\&-g_2(\calJ_{\epsilon} z_{\epsilon}(t,x))
-\theta_{\epsilon}(t,x) 
\D_z H_2(\calJ_{\epsilon} z_{\epsilon} (t,x) )) 
 \in \partial\Psi(\dot z_{\epsilon} (t,x) ).
\end{aligned}
\end{equation*}
for all \(t\in [0,T]\) and almost every \(x \in \Omega\).
Since we are looking for a solution with values 
in $\Erdev$, the test-functions should be taken in $\Erdev$. More precisely,
we have 
\begin{equation}
\label{eq:30}
\begin{aligned} 
&\Psi ( \tilde z)-\Psi (\dot z_{\epsilon} (t,x)) 
+(\Me\dot z_{\epsilon} (t,x){+}
\partial\varphi_{\epsilon}(z_{\epsilon} (t,x) ){+}
\Ee z_{\epsilon} (t,x){+}\D_z H_1(z_{\epsilon} (t,x))\\& 
{+} g_1(\theta_{\epsilon}(t,x))   
{+}g_2(\calJ_{\epsilon} z_{\epsilon}(t,x))
{+}\theta_{\epsilon} (t,x) 
\D_z H_2(\calJ_{\epsilon} z_{\epsilon} (t,x) )){:} 
(\tilde z{-} z_{\epsilon}(t,x)) \geq 0,
\end{aligned}
\end{equation}
for all \(\tilde z \in \Erdev\)
and $t \in [0,T]$ and almost every $x \in \Omega$. Observe that 
\eqref{eq:30} is equivalent to
\begin{equation*}
\begin{aligned}
&\Psi_{\textrm{dev}}(\tilde z)-\Psi_{\textrm{dev}}(\dot z_{\epsilon}(t,x))  
+{\textrm{Proj}}_{\Erdev}(\Me\dot z_{\epsilon}(t,x){+}
\partial\varphi_{\epsilon}(z_{\epsilon} (t,x))
{+}\Ee z_{\epsilon}(t,x) 
{+}\D_z H_1(z_{\epsilon}(t,x))\\&
{+} g_1(\theta_{\epsilon}(t,x))
{+}g_2(\calJ_{\epsilon} z_{\epsilon})(t,x) {+}\theta_{\epsilon} (t,x) 
\D_z H_2(\calJ_{\epsilon} z_{\epsilon}(t,x))){:}
(\tilde z{-}z_{\epsilon}(t,x))\ge 0,
\end{aligned}
\end{equation*}
for all \(\tilde z \in \Erdev\)
and $t \in [0,T]$ and almost every $x \in \Omega$. 
Here $\Psi_{\textrm{dev}}$ is the restriction of $\Psi$ to $\Erdev$ and 
${\textrm{Proj}}_{\Erdev}$ denotes the projection on $\Erdev$ 
relatively to the inner product of $\Ers$. Using the definition
of the subdifferential \(\partial\Psi_{\textrm{dev}}(\dot z_{\epsilon}(t,x))\)
leads to the following differential inclusion
\begin{equation*}
\begin{aligned}
&-{\textrm{Proj}}_{\Erdev}
(\partial\varphi_{\epsilon}(z_{\epsilon} (t,x))
{+}\Ee z_{\epsilon}(t,x){+}\D_z H_1(z_{\epsilon}(t,x)) 
+ g_1(\theta_{\epsilon}(t,x))\\&
{+}g_2(\calJ_{\epsilon}z_{\epsilon}(t,x)) 
{+}\theta_{\epsilon}(t,x) 
\D_z H_2(\calJ_{\epsilon} z_{\epsilon}(t,x))) 
\in \partial\Psi_{\textrm{dev}} 
(\dot z_{\epsilon} (t,x) ) + \Me\dot z_{\epsilon}(t,x), 
\end{aligned}
\end{equation*}
for all $t \in [0,T]$ and almost every $x \in \Omega$. 
Note that the operator \(\partial\Psi_{\textrm{dev}}+\Me\) is strongly
monotone on \(\Erdev\) as well as its $\L^2 (\Omega)$-extension, still 
denoted \(\partial\Psi_{\textrm{dev}}+\Me\), on \(\L^2(\Omega;\Erdev)\). This
operator is invertible and
its inverse is $\tfrac1{c^{\Me}}$-Lipschitz continuous. 
Thus \eqref{eq:57} can be rewritten as 
\begin{equation} \label{eq:57bis}
\dot z_{\epsilon}=\Phi(\theta_{\epsilon}, z_{\epsilon}),
\end{equation}
with
\begin{equation*}
\Phi(\theta_{\epsilon}, z)\eqldef(\partial \Psi_{\textrm{dev}}{+}\Me)^{-1} 
\bigl({\textrm{Proj}}_{\Erdev} ({-}\partial\varphi_{\epsilon}(z)
{-}\Ee z{-}\D_z H_1(z)
{-}g_1(\theta_{\epsilon}){-}g_2(\calJ_{\epsilon} z){-}\theta_{\epsilon}
\D_z H_2(\calJ_{\epsilon} z))\bigr), 
\end{equation*}
for all $z \in \L^2(0,T; \L^2(\Omega; \Erdev))$ and we have to solve 
this  differential equation for the unknown function $z_{\epsilon}$.
We may observe that we are not in the classical framework of ODE in 
\(\L^2(\Omega;\Erdev)\) (see \cite{Car90}), since $g_2(\calJ_{\epsilon} z) $ 
can not be defined for $z \in \L^2 (\Omega; \Erdev)$ but only for 
$z \in \L^2 (0,T;\L^2(\Omega;\Erdev))$. Nevertheless, we can solve 
\eqref{eq:57bis} with the Picard's iterations technique. Indeed, 
$\Phi(\theta_{\epsilon}, z) \in \C^0([0,T];\L^2(\Omega))$ for all 
$z \in \C^0([0,T];\L^2(\Omega))$ and \(z_{\epsilon} \in
\C^1([0,T];\L^2(\Omega))\) 
is a solution of  \eqref{eq:57}--\eqref{eq:87} if and only if 
\(z_{\epsilon}\) is a fixed point of the mapping
\begin{equation*}
\begin{aligned}
\Lambda_{\epsilon}:\ \C^0([0,T];\L^2(\Omega)) &\to \C^1([0,T];\L^2(\Omega)) \\
z &\mapsto  \Lambda_{\epsilon} (z) : t \mapsto  z^0 + \int_0^t 
\Phi(\theta_{\epsilon}(\cdot,s), z(\cdot,s))\dd s.
\end{aligned}
\end{equation*}
Let us assume that $z_1, z_2 \in \C^0([0,T];\L^2(\Omega))$. 
By using \eqref{eq:H3} and \eqref{eq:46}, we find that
\begin{equation*}
\begin{aligned}
& \norm[\L^2(\Omega)]{\Lambda_{\epsilon} (z_1(\cdot,t)){-}\Lambda_{\epsilon} 
(z_2(\cdot,t))}
\leq \tfrac{1}{c^{\Me}}\int_0^t \bigl( \norm[\L^2(\Omega)]{\partial 
\varphi_{\epsilon} (z_1(\cdot,s)){-}\partial \varphi_{\epsilon} (z_2(\cdot,s))}\\
&{+}(\norm[\L^{\infty}(\Omega)]{\Ee}{+}C_{zz}^{H_1}{+} 
C_{zz}^{H_2} \norm[\L^{\infty} (\calQ_T)]{\theta_{\epsilon}}) 
\norm[\L^2(\Omega)]{z_1(\cdot,s){-}z_2(\cdot,s)}{+} 
\norm[\L^2(\Omega)]{g_2(\calJ_{\epsilon} z_1(\cdot,s)) 
{-}g_2(\calJ_{\epsilon} z_2(\cdot,s))}\bigr)\dd s,
\end{aligned}
\end{equation*}
for all $t\in[0,T]$.
Since $\partial \varphi_{\epsilon}$ is 
$\tfrac1{\epsilon}$-Lipschitz continuous on $\L^2 (\Omega)$ and $g_2$ 
is linear, it follows that
\begin{equation} 
\label{eqlp:2}
\begin{aligned}
\forall t \in[0,T]: \
&\norm[\L^2(\Omega)]{\Lambda_{\epsilon} (z_1(\cdot,t)){-} 
\Lambda_{\epsilon} (z_2(\cdot,t))}\le C^{\epsilon} 
\int_0^t  \norm[\L^2(\Omega)]{z_1(\cdot,s){-}z_2(\cdot,s)}\dd s\\&
+ \tfrac{1}{c^{\Me}}  \int_0^t  \norm[\L^2(\Omega)]{g_2(\calJ_{\epsilon} 
z_1(\cdot,s){-}\calJ_{\epsilon} z_2(\cdot,s))}\dd s,
\end{aligned}
\end{equation}
where 
\(
C^{\epsilon}\eqldef\frac{1}{c^{\Me}}\bigl(\tfrac{1}{\epsilon}{+} 
\norm[\L^{\infty}(\Omega)]{\Ee}{+}C_{zz}^{H_1}{+}C_{zz}^{H_2} 
\norm[\L^{\infty}(\calQ_T)]{\theta_{\epsilon}}\bigr)\).
Thus the energy estimate \eqref{eqlp:3} allows us to infer that
\begin{equation*}
\begin{aligned}
\forall s \in [0,T]\ \forall z \in \L^2 (0,T; \L^2(\Omega)):\
\norm[\L^2(\Omega)]{g_2 (z(\cdot,s))}&\leq 
\tfrac{\norm[\L^{\infty}(\Omega)]{\Ee}}{\sqrt{c^{\Ee} 
c^{\Le}\cK}} \Bigl(\int_0^s \norm[(\H^1_0(\Omega))']{\dive (\Ee
z(\cdot,\sigma))}^2 
\dd\sigma  
\Bigr)^{\frac12}\\
& \leq\tfrac{\norm[\L^{\infty}(\Omega)]{\Ee}^2}{\sqrt{c^{\Ee}c^{\Le}\cK}} 
\Bigl(\int_0^s \norm[\L^2(\Omega)]{z(\cdot,\sigma)}^2\dd\sigma\Bigr)^{\frac12},
\end{aligned}
\end{equation*}
which, leads to
\begin{equation*}
\begin{aligned}
\forall t \in [0,T]:\ &
\int_0^t  \norm[\L^2(\Omega)]{g_2(\calJ_{\epsilon} 
z_1(\cdot,s){-}\calJ_{\epsilon} z_2(\cdot,s))}\dd s \\&\le
 \tfrac{\norm[\L^{\infty}(\Omega)]{\Ee}^2}{\sqrt{c^{\Ee}c^{\Le}\cK}} 
\int_0^t \Bigl( \int_0^s \norm[\L^2(\Omega)]{\calJ_{\epsilon} z_1
(\cdot,\sigma){-} 
\calJ_{\epsilon} z_2(\cdot,\sigma)}^2\dd \sigma\Bigr)^{\frac12} \dd s.
\end{aligned}
\end{equation*}
Since  $\calJ_{\epsilon}$ is a contraction on $\L^2(\Omega)$ we infer
\begin{equation*}
\begin{aligned}
\forall t \in [0,T]:\
\int_0^t  \norm[\L^2(\Omega)]{g_2(\calJ_{\epsilon} z_1(\cdot,s){-} 
\calJ_{\epsilon} z_2(\cdot,s))}\dd s &\leq
 \tfrac{\norm[\L^{\infty}(\Omega)]{\Ee}^2}{\sqrt{c^{\Ee} c^{\Le}\cK}} \int_0^t 
\Bigl(\int_0^s \norm[\L^2(\Omega)]{z_1(\cdot,\sigma){-}z_2(\cdot,\sigma)}^2 
\dd \sigma   
\Bigr)^{\frac12} \dd s\\&
\leq  
\tfrac{\norm[\L^{\infty}(\Omega)]{\Ee}^2 \sqrt{T}}{\sqrt{c^{\Ee}  
c^{\Le}\cK}} t \norm[\C^0(0,T; \L^2(\Omega))]{z_1{-}z_2}.
\end{aligned}
\end{equation*}
Finally, we find
\begin{equation*}
\forall t \in [0,T]:\
\norm[\L^2(\Omega)]{\Lambda_{\epsilon} (z_1(\cdot,t)){-}\Lambda_{\epsilon} 
(z_2(\cdot,t))} 
\le \bigl(C^{\epsilon}{+}\tfrac{\norm[\L^{\infty}(\Omega)]{\Ee}^2
    \sqrt{T}}{c^{\Me} 
\sqrt{c^{\Ee}c^{\Le}\cK}}\bigr)t\norm[\C^0(0,T; \L^2(\Omega))]{z_1{-}z_2}.
\end{equation*}
We iterate the previous computation to get 
\begin{equation*}
\forall t \in [0,T]\ \forall m\geq 1:\
\norm[\L^2(\Omega)]{\Lambda_{\epsilon}^m (z_1(\cdot,t)){-}\Lambda_{\epsilon}^m 
(z_2(\cdot,t))} 
\le \bigl(C^{\epsilon}{+}\tfrac{\norm[\L^{\infty}(\Omega)]{\Ee}^2
\sqrt{T}}{c^{\Me} 
\sqrt{c^{\Ee}c^{\Le}\cK}}\bigr)^m \tfrac{t^m}{m!} 
\norm[\C^0(0,T; \L^2(\Omega))]{z_1{-}z_2}.
\end{equation*}
It follows that there exists $m_0 \ge 1$ such that 
\begin{equation*}
\bigl(C^{\epsilon}{+}\tfrac{\norm[\L^{\infty}(\Omega)]{\Ee}^2 
\sqrt{T}}{c^{\Me} \sqrt{c^{\Ee}c^{\Le}\cK}}\bigr)^{m_0}\tfrac{T^{m_0}}{m_0!} < 1,
\end{equation*}
and $\Lambda_{\epsilon}^{m_0}$ possesses a unique fixed point in 
$\C^0 ([0,T]; \L^2(\Omega))$, which is also the unique fixed point of 
$\Lambda_{\epsilon}$. We may conclude that
there exists a unique solution \(z_{\epsilon}\in\C^1([0,T];\L^2(\Omega; 
\Erdev))\) 
to the problem \eqref{eq:57}--\eqref{eq:87}.

Let us choose now a sequence \((\theta_{\epsilon})_{\epsilon>0}\) such that
\begin{equation*}
\theta_{\epsilon}
\rightarrow
\theta\ \text{ in } \ \L^q(0,T;\L^p(\Omega)).
\end{equation*}
Since \(p, q\geq 2\), it follows that
\begin{equation*}
\theta_{\epsilon}
\rightarrow
\theta\ \text{ in } \ \L^2(0,T;\L^2(\Omega)).
\end{equation*}
We will establish that there exists a subsequence of
\((z_{\epsilon})_{\epsilon>0}\), still denoted by \((z_{\epsilon})_{\epsilon>0}\),
such that 
\begin{equation*}
\begin{aligned}
z_{\epsilon}
&\rightharpoonup
z \ \text{ in } \ \H^1(0,T;\L^2(\Omega;\Erdev)) \ \text{ weak, }\\
z_{\epsilon}
&\rightharpoonup
z \ \text{ in } \ \L^{\infty}(0,T;\H^1(\Omega;\Erdev)) \ \text{ weak }\ *,
\end{aligned}
\end{equation*}
where \(z\) is a solution of the problem \eqref{eq:5}--\eqref{eq:9}.
To do so, we notice that
\begin{equation*}
\forall t\in[0,T]:\
\dot{z}_{\epsilon}(\cdot,t)=
(\partial\Psi_{\textrm{dev}}{+}\Me)^{-1}\bigl({\textrm{Proj}}_{\Erdev} 
(g_{\epsilon}(t,\calJ_{\epsilon}z_{\epsilon}(\cdot,t))
{-}\partial\varphi_{\epsilon}(z_{\epsilon}(\cdot,t))
{-}\Ee z_{\epsilon}(\cdot,t){-}\D_z H_1(z_{\epsilon}(\cdot,t)))\bigr),
\end{equation*}
where \(g_{\epsilon}(t,\calJ_{\epsilon}z_{\epsilon}(\cdot,t))\eqldef
-(g_1(\theta_{\epsilon}(\cdot,t)){+}
g_2(\calJ_{\epsilon} z_{\epsilon}(\cdot,t)){+}\theta_{\epsilon}(\cdot,t)
\D_z H_2(\calJ_{\epsilon} z_{\epsilon}(\cdot,t)))\).
Define  
\begin{equation*}
\forall t\in [0,T]:\
w_{\epsilon}(\cdot,t)\eqldef g_{\epsilon}(t,\calJ_{\epsilon}z_{\epsilon}(\cdot,t))
-\partial\varphi_{\epsilon}(z_{\epsilon}(\cdot,t))
-\Ee z_{\epsilon}(\cdot,t)-\D_z H_1(z_{\epsilon}(\cdot,t)).
\end{equation*}
We have \(w_{\epsilon}\in\C^0([0,T];\L^2(\Omega;\Ers))\) and 
\begin{subequations}
\label{eq:90}
\begin{align}
&\forall t\in[0,T]:\ 
w_{\epsilon}(\cdot,t)+\partial\varphi_{\epsilon}(z_{\epsilon}(\cdot,t))
+\Ee z_{\epsilon}(\cdot,t)+\D_z H_1(z_{\epsilon}(\cdot,t))=
g_{\epsilon}(t,\calJ_{\epsilon}z_{\epsilon}(\cdot,t)),\label{eq:92}\\
&\forall t\in[0,T]:\  \dot z_{\epsilon}(\cdot,t) =
(\partial\Psi_{\textrm{dev}}{+}\Me)^{-1}\bigl({\textrm{Proj}}_{\Erdev} 
(w_{\epsilon}(\cdot,t))\bigr).
\label{eq:93}
\end{align}
\end{subequations}
In order to obtain a priori estimates, we multiply \eqref{eq:92} by 
\(\dot z_{\epsilon}\), and we integrate this expression over 
\(\calQ_{\tau}\) to get
\begin{equation*}
\int_{\calQ_{\tau}}w_{\epsilon}{:}\dot z_{\epsilon}\dd x\dd t
+\int_{\calQ_{\tau}}\partial\varphi_{\epsilon}(z_{\epsilon}){:}
\dot z_{\epsilon}\dd x\dd t
+\int_{\calQ_{\tau}}\Ee z_{\epsilon}{:}\dot z_{\epsilon}\dd x\dd t
+\int_{\calQ_{\tau}}\D_z H_1(z_{\epsilon}){:}\dot z_{\epsilon}\dd x\dd t=
\int_{\calQ_{\tau}}g_{\epsilon}(t,\calJ_{\epsilon}z_{\epsilon}){:}\dot 
z_{\epsilon}\dd x\dd t.
\end{equation*}
We observe that
\(\int_{\calQ_{\tau}}\partial\varphi_{\epsilon}(z_{\epsilon}){:}\dot 
z_{\epsilon}\dd x\dd t
=\varphi_{\epsilon}(z_{\epsilon}(\cdot,\tau))-\varphi_{\epsilon}(z^0)\),
\(\int_{\calQ_{\tau}}\D_z H_1(z_{\epsilon}){:}\dot z_{\epsilon}\dd x\dd t=
\int_{\Omega}H_1(z_{\epsilon}(\cdot,\tau))\dd x-
\int_{\Omega}H_1(z^0)\dd x\). Moreover, recalling that $\dot z_{\epsilon} 
\in \C^0( [0,T];\L^2(\Omega;\Erdev))$, we have 
\begin{equation*}
\int_{\calQ_{\tau}}w_{\epsilon}{:}\dot z_{\epsilon}\dd x\dd t
= \int_{\calQ_{\tau}} {\rm Proj}_{\Erdev} (w_{\epsilon}) {:}
\dot z_{\epsilon}\dd x\dd t
\geq c^{\Me}\int_0^{\tau}\norm[\L^2(\Omega)]{\dot z_{\epsilon}}^2\dd t.
\end{equation*}
Since \(z^0\in\H^1(\Omega;\Erdev)=\calD(\varphi)\), we have 
\(\varphi_{\epsilon}(z^0)\leq \varphi(z^0)\)
(for technical details, the reader is referred to \cite{Brez73OMMS}) and 
\begin{equation*}
\begin{aligned}
&\tfrac{c^{\Me}}2
\int_0^{\tau}\norm[\L^2(\Omega)]{\dot z_{\epsilon}}^2\dd t+
\varphi_{\epsilon}(z_{\epsilon}(\cdot,\tau))+
\int_{\Omega}H_1(z_{\epsilon}(\cdot,\tau))\dd x+
\tfrac12\int_{\Omega}\Ee z_{\epsilon}(\cdot,\tau){:} z_{\epsilon}(\cdot,\tau)\dd x
\\&\leq
\varphi(z^0)+
\int_{\Omega}H_1(z^0)\dd x+
\tfrac12\int_{\Omega}\Ee z^0{:} z^0\dd x+
\tfrac1{2c^{\Me}}\int_0^{\tau}\norm[\L^2(\Omega)]{
g_{\epsilon}(t,\calJ_{\epsilon}z_{\epsilon}(\cdot,t))}^2\dd t.
\end{aligned}
\end{equation*}
Furthermore
\begin{equation*}
\forall\epsilon>0\ \forall \bar z\in\L^2(\Omega;\Erdev):\
\varphi_{\epsilon}(\bar z)=
\tfrac{\epsilon}2 \norm[\L^2(\Omega)]{\partial\varphi_{\epsilon}(\bar
  z)}^2+\varphi (\calJ_{\epsilon}\bar z) = \tfrac{1}{2 \epsilon} 
\norm[\L^2(\Omega)]{\bar z {-} \calJ_{\epsilon} \bar z}^2 + 
\varphi (\calJ_{\epsilon}\bar z).
\end{equation*}
Thus \(\calJ_{\epsilon}\bar z\in \C^0 ([0,T]; \H^1(\Omega;\Erdev))\)
for all \(\bar z\in \C^0([0,T]; \L^2(\Omega;\Erdev))\) and
\begin{equation*}
\begin{aligned}
&\tfrac{c^{\Me}}2
\int_0^{\tau}\norm[\L^2(\Omega)]{\dot z_{\epsilon}}^2\dd t+
\tfrac{1}{2 \epsilon} \norm[\L^2(\Omega)]{z_{\epsilon} (\cdot, \tau) 
{-}\calJ_{\epsilon} z_{\epsilon}(\cdot, \tau)}^2 + 
\varphi(\calJ_{\epsilon}z_{\epsilon}(\cdot,\tau))+
\int_{\Omega}H_1(z_{\epsilon}(\cdot,\tau))\dd x\\&+
\tfrac12\int_{\Omega}\Ee z_{\epsilon}(\cdot,\tau){:} z_{\epsilon}(\cdot,\tau)\dd x
\leq
\varphi(z^0)+
\int_{\Omega}H_1(z^0)\dd x+
\tfrac12\int_{\Omega}\Ee z^0{:} z^0\dd x+
\tfrac1{2c^{\Me}}\int_0^{\tau}\norm[\L^2(\Omega)]{
g_{\epsilon}(t,\calJ_{\epsilon}z_{\epsilon}(\cdot,t))}^2\dd t.
\end{aligned}
\end{equation*}
By using \eqref{eq:H1}, \eqref{eq:H2} and \eqref{eq:3}, we infer that 
\begin{equation} \label{eqlp:4}
\begin{aligned}
&\tfrac{c^{\Me}}2\int_0^{\tau}
\norm[\L^2(\Omega)]{\dot z_{\epsilon}}^2\dd t+
\tfrac{1}{2 \epsilon} \norm[\L^2(\Omega)]{z_{\epsilon} (\cdot, \tau) 
{-}\calJ_{\epsilon} z_{\epsilon}(\cdot, \tau)}^2 +
C_1\norm[\L^2(\Omega)]{z_{\epsilon}(\cdot,\tau)}^2\\&+
C_1\norm[\L^2(\Omega)]{\nabla(\calJ_{\epsilon}z_{\epsilon}(\cdot,\tau))}^2 
\leq
C_0+\tfrac1{2c^{\Me}}\int_0^{\tau}\norm[\L^2(\Omega)]{
g_{\epsilon}(t,\calJ_{\epsilon}z_{\epsilon}(\cdot,t))}^2\dd t,
\end{aligned}
\end{equation}
where 
\(C_0\eqldef \varphi(z^0)+\int_{\Omega}H_1(z^0)\dd x+ \tilde c^{H_1}\abs{\Omega}+
\tfrac12\int_{\Omega}\Ee z^0{:}z^0\dd x\) and 
\(C_1\eqldef
\min\bigl(c^{H_1}{+}\tfrac{c^{\Ee}}2,\tfrac{\nu}2\bigr)
\).
Since \(\calJ_{\epsilon}\) is a contraction and
\(\calJ_{\epsilon}0=0\), we have
\begin{equation*}
\forall \bar z\in\L^2(\Omega;\Erdev):\
\norm[\L^2(\Omega)]{\calJ_{\epsilon}\bar z}=
\norm[\L^2(\Omega)]{\calJ_{\epsilon}\bar z{-}\calJ_{\epsilon}0}
\leq \norm[\L^2(\Omega)]{\bar z}.
\end{equation*}
Then we find
\begin{equation}
\label{eq:94bis}
\begin{aligned}
\tfrac{c^{\Me}}2\int_0^{\tau}\norm[\L^2(\Omega)]{\dot z_{\epsilon}}^2\dd t+
C_1\norm[\H^1(\Omega)]{\calJ_{\epsilon}z_{\epsilon}(\cdot,\tau)}^2
\leq
C_0+\tfrac1{2c^{\Me}}\int_0^{\tau}\norm[\L^2(\Omega)]{
g_{\epsilon}(t,\calJ_{\epsilon}z_{\epsilon}(\cdot,t))}^2\dd t.
\end{aligned}
\end{equation}
For the bound on \(g_{\epsilon}(t,\calJ_{\epsilon}z_{\epsilon}(\cdot,t))\), we use
\eqref{eq:H2} to get
\begin{equation}
\label{eq:95}
\begin{aligned}
&
\int_0^{\tau}\norm[\L^2(\Omega)]{
g_{\epsilon}(t,\calJ_{\epsilon}z_{\epsilon})}^2\dd t
\leq
3\int_0^{\tau}\norm[\L^2(\Omega)]{g_1(\theta_{\epsilon})}^2\dd t
+
3\int_0^{\tau}\norm[\L^2(\Omega)]{g_2(\calJ_{\epsilon}z_{\epsilon})}^2\dd t\\&+
6 (C_z^{H_2})^2 \int_0^{\tau} \norm[\L^2(\Omega)]{\theta_{\epsilon}}^2 \dd t + 
6(C_z^{H_2})^2 \int_{\calQ_{\tau}}\abs{\theta_{\epsilon}}^2
\abs{\calJ_{\epsilon} z_{\epsilon}}^{2}
\dd x\dd t.
\end{aligned}
\end{equation}
By using the definition of mappings \(g_1\) and \(g_2\), the first two
terms on the right hand side of \eqref{eq:95} can be estimated. More
precisely, it is plain that there exists \(C_2>0\) such that
\begin{equation}
\label{eq:96}
\begin{aligned}
&\int_0^{\tau}\norm[\L^2(\Omega)]{g_1(\theta_{\epsilon})}^2\dd t+
\int_0^{\tau}\norm[\L^2(\Omega)]{g_2(\calJ_{\epsilon} z_{\epsilon})}^2\dd t
\\&\leq \tfrac{C_2}3\Bigl( \norm[\H^1(\Omega)]{u^0}^2 {+} 
 \int_0^{\tau}\norm[\L^2(\Omega)]{\ell}^2\dd t{+}
\int_0^{\tau}\norm[\L^2(\Omega)]{\theta_{\epsilon}}^2\dd t{+}
\int_0^{\tau}\norm[\L^2(\Omega)]{\calJ_{\epsilon} z_{\epsilon}}^2\dd t
\Bigr).
\end{aligned}
\end{equation}
The last term on the right hand side of \eqref{eq:95} is estimated by
using  H\"older's inequality and the continuous embedding $\H^1(\Omega) 
\hookrightarrow \L^4(\Omega)$. Namely, it follows that there 
exists \(C_3>0\) such that
\begin{equation}
\label{eq:97}
\begin{aligned}
\int_{\calQ_{\tau}}\abs{\theta_{\epsilon}}^2
\abs{\calJ_{\epsilon} z_{\epsilon}}^{2}\dd x\dd t
&\leq
C_3 \int_0^{\tau} \norm[\L^4(\Omega)]{\theta_{\epsilon}}^2 
\norm[\H^1(\Omega)]{\calJ_{\epsilon}z_{\epsilon}}^2 \dd t.
\end{aligned}
\end{equation}
We insert \eqref{eq:96} and \eqref{eq:97} into \eqref{eq:94bis}, we obtain
\begin{equation}
\label{eq:98}
\begin{aligned}
&\tfrac{c^{\Me}}2\int_0^{\tau}\norm[\L^2(\Omega)]{\dot z_{\epsilon}}^2\dd t+
C_1\norm[\H^1(\Omega)]{\calJ_{\epsilon}z_{\epsilon}(\cdot,\tau)}^2
\leq
C_0+\tfrac{C_2}{2c^{\Me}} \Bigl(\norm[\H^1(\Omega)]{u^0}^2{+} 
\int_0^{\tau}\norm[\L^2(\Omega)]{\ell}^2 \dd t\Bigr)
\\& +\tfrac{C_2{+}6 (C_z^{H_2})^2}{2c^{\Me}} \int_0^{\tau}
\norm[\L^2(\Omega)]{\theta_{\epsilon}}^2\dd t+
\int_0^{\tau} 
\tfrac{C_2{+}6(C_z^{H_2})^2 C_3\norm[\L^4(\Omega)]{\theta_{\epsilon}}^2}{2c^{\Me}}
\norm[\H^1(\Omega)]{\calJ_{\epsilon}z_{\epsilon}}^2
\dd t.
\end{aligned}
\end{equation}
Since \(\theta_{\epsilon}\) is bounded in \(\L^q(0,T;\L^p(\Omega))\), 
with $q \ge 2$ and $p \ge 4$,
we may deduce from  Gr\"onwall's lemma that 
\(\calJ_{\epsilon} z_{\epsilon}\) is bounded in \(\L^{\infty}(0,T;\H^1(\Omega))\),
\(g_{\epsilon}(\cdot ,\calJ_{\epsilon}z_{\epsilon})\) is bounded in
\(\L^2(0,T;\L^2(\Omega))\) and 
\(z_{\epsilon}\) is bounded in \(\H^1(0,T;\L^2(\Omega)) \cap \L^{\infty} 
(0,T; \L^2(\Omega)) \),
independently of \(\epsilon>0\).
Therefore \eqref{eq:H2} implies that $D_z H_1 (z_{\epsilon})$ is bounded 
in \(\L^{\infty}(0,T;\L^2(\Omega)) \),
independently of \(\epsilon>0\). Furthermore 
the definition of the subdifferential 
\(\partial\Psi_{\textrm{dev}}(\dot z_{\epsilon})\) enables 
us to infer from \eqref{eq:93} that
\begin{equation*}
- \int_{\calQ_{\tau}}\Me\dot z_{\epsilon}{:}w^p_{\epsilon}\dd x\dd t+
\int_0^{\tau}\norm[\L^2(\Omega)]{w^p_{\epsilon}}^2\dd t
 \leq \int_{\calQ_{\tau}} 
(\Psi(w^p_{\epsilon}{+}\dot z_{\epsilon}){-}\Psi( \dot z_{\epsilon}) )\dd x\dd t
 \leq  \int_{\calQ_{\tau}} \Psi(w^p_{\epsilon} ) \dd x\dd t ,
\end{equation*}
where $w^p_{\epsilon}={\textrm{Proj}}_{\Erdev}(w_{\epsilon})$, and using
 \eqref{eq:Psi.bdd}, we get
\begin{equation}
\label{eq:100}
\norm[\L^2(0,T;\L^2(\Omega))]{w^p_{\epsilon}}
\leq 
C^{\Psi}\sqrt{T\abs{\Omega}}+
C^{\Me} \norm[\L^2(0,T;\L^2(\Omega))]{\dot z_{\epsilon}}.
\end{equation}
Observe that \(\dot z_{\epsilon}\) is bounded in \(\L^2(0,T;\L^2(\Omega))\),
independently of \(\epsilon>0\), allows us to
conclude that \(w^p_{\epsilon}\) is bounded in \(\L^2(0,T;\L^2(\Omega))\),
independently of \(\epsilon>0\). On the other hand, we multiply 
\eqref{eq:92} by \(\partial\varphi_{\epsilon}(z_{\epsilon})\) and  we
integrate this result over \(\calQ_{\tau}\). Recalling that $\partial
\varphi_{\epsilon}( z_{\epsilon})$ takes its values in $L^2 (\Omega; \Erdev)$,   
the Cauchy-Schwarz's
inequality, \eqref{eq:H2} and \eqref{eq:100} give 
\begin{equation*}
\begin{aligned}
&
\norm[\L^2(0,T;\L^2(\Omega))]{\partial\varphi_{\epsilon}(z_{\epsilon})}
\leq
(C_z^{H_1}{+} C^{\Psi})\sqrt{T\abs{\Omega}}
+(C_z^{H_1}{+} \norm[\L^{\infty}(\Omega)]{\Ee})
\norm[\L^2(0,T;\L^2(\Omega))]{z_{\epsilon}}\\&+
C^{\Me}\norm[\L^2(0,T;\L^2(\Omega))]{\dot z_{\epsilon}}
+ \norm[\L^2(0,T;\L^2(\Omega))]{ g_{\epsilon} ( \cdot, \calJ_{\epsilon} 
z_{\epsilon}) } .
\end{aligned}
\end{equation*}
Thus \(\partial\varphi_{\epsilon}(z_{\epsilon})\) is bounded in
\(\L^2(0,T;\L^2(\Omega))\), independently of \(\epsilon>0\) and finally 
$w_{\epsilon}$ is bounded in
\(\L^2(0,T;\L^2(\Omega))\), independently of \(\epsilon>0\). 
Furthermore,
there exists \(z\in\H^1(0,T;\L^2(\Omega))\cap\L^{\infty}(0,T;\L^2(\Omega))\),
\(\tilde z\in\L^{\infty}(0,T;\H^1(\Omega))\),
\(w, v\in\L^2(0,T;\L^2(\Omega))\) such that it is possible to extract
subsequences, still denoted by \(z_{\epsilon}\),
\(\calJ_{\epsilon}z_{\epsilon}\),
\(w_{\epsilon}\) and \(\partial\varphi_{\epsilon}(z_{\epsilon})\) satisfying the
following convergences
\begin{equation*}
\begin{aligned}
z_{\epsilon}&\rightharpoonup z \ \text{ in } \ \H^1(0,T;\L^2(\Omega))
\ \text{ weak},\\
z_{\epsilon}&\rightharpoonup z \ \text{ in } \ \L^{\infty}(0,T;\L^2(\Omega))
\ \text{ weak }*,\\
\calJ_{\epsilon}z_{\epsilon}&\rightharpoonup 
\tilde z \ \text{ in } \ \L^{\infty}(0,T;\H^1(\Omega))
\ \text{ weak }*,\\
w_{\epsilon}&\rightharpoonup w\ \text{ in } \ \L^2(0,T;\L^2(\Omega))
\ \text{ weak},\\
\partial\varphi_{\epsilon}(z_{\epsilon})
&\rightharpoonup v
\ \text{ in } \ \L^2(0,T;\L^2(\Omega))
\ \text{ weak}.
\end{aligned}
\end{equation*}
Moreover, using \eqref{eqlp:4} there exists $C_4>0$ such that
\begin{equation} 
\label{eq:102}
\forall \tau \in [0, T]:\
\norm[\L^2(\Omega)]{z_{\epsilon} (\cdot, \tau){-}\calJ_{\epsilon} 
z_{\epsilon}(\cdot, \tau)}^2 \le C_4 \epsilon \bigl(1{+}  
\norm[\L^2(0,T;\L^2(\Omega))]{g_{\epsilon} ( \cdot, 
\calJ_{\epsilon} z_{\epsilon})}^2 \bigr),
\end{equation}
which, allows us to infer that \(z=\tilde z\). 
Since $\calJ_{\epsilon}$ is a contraction on $\L^2(\Omega)$, 
it follows also that, possibly extracting another subsequence,   
still denoted by \(\calJ_{\epsilon}z_{\epsilon}\),
we have
\begin{equation*}
\forall r\in[2,6):\
\calJ_{\epsilon}z_{\epsilon}\rightarrow z
\ \text{ in }\ \C^0([0,T];\L^r(\Omega)).
\end{equation*}
We may deduce from \eqref{eq:102} that
\begin{equation*}
z_{\epsilon}\rightarrow z
\ \text{ in } \ \C^0(0,T;\L^2(\Omega)).
\end{equation*}
Our aim now consists to pass to the limit in \eqref{eq:90}. To do so,
we observe that the mapping \(\calL_0\) is linear and continuous, the
mappings \(g_1\) and \(g_2\) are also continuous, which, gives 
\begin{equation*}
\begin{aligned}
g_1(\theta_{\epsilon})&\rightarrow g_1(\theta)
\ \text{ in } \ \L^2(0,T;\L^2(\Omega)),\\
g_2(\calJ_{\epsilon}z_{\epsilon})&\rightarrow g_2(z)
\ \text{ in } \ \L^2(0,T;\L^2(\Omega)).
\end{aligned}
\end{equation*}
Since \(\D_z H_i\), \(i=1,2\), is Lipschitz continuous, it follows that
\begin{equation*}
\begin{aligned}
\D_z H_1(z_{\epsilon})+\Ee z_{\epsilon}
&\rightarrow \D_z H_1(z)+\Ee z \ \text{ in } \ \C^0([0,T];\L^2(\Omega)),\\
\forall r\in[2,6):\
\D_z H_2(\calJ_{\epsilon} z_{\epsilon})&\rightarrow \D_z H_2(z)
\ \text{ in } \ \C^0([0,T];\L^r(\Omega)).
\end{aligned}
\end{equation*}
Hence we consider \(r\) such that \(\tfrac1p+\tfrac1r=\tfrac12\) to get
\begin{equation*}
\theta_{\epsilon}\D_z H_2(\calJ_{\epsilon}z_{\epsilon})
\rightarrow 
\theta\D_z H_2(z)\ \text{ in } \ \L^2(0,T;\L^2(\Omega)).
\end{equation*}
We can pass to the limit in all terms of \eqref{eq:92}, we obtain
\begin{equation}
\label{eq:106}
w+v+\D_z H_1(z)+\Ee z=-(g_1(\theta){+}g_2(z){+}\theta\D_z H_2(z)).
\end{equation}
It remains to pass to the limit in \eqref{eq:93}
and to prove that $z$ solves \eqref{eq:5}--\eqref{eq:9}. 
Using \eqref{eq:87} and the strong convergence of 
$(z_{\epsilon})_{\epsilon>0}$ to $z$ in $\C^0([0,T];\L^2(\Omega))$, 
we infer that $z(0, \cdot) = z^0$. 
Then we prove that \(v (t, \cdot) \in\partial\varphi(z (t, \cdot) )\) 
for almost every $t \in [0,T]$, i.e, \(v\in\partial\varphi(z)\) 
where $\partial \varphi$ is identified to its $\L^2 (0,T)$-extension. 
To do so, 
we use the classical results for maximal monotone operators 
(see \cite{Brez73OMMS}).
More precisely, since 
\(\partial\varphi_{\epsilon}(z_{\epsilon})\in\partial\varphi(\calJ_{\epsilon}
z_{\epsilon})\), it is sufficient to prove that
\begin{equation} \label{eq:104}
\limsup_{\epsilon\rightarrow 0}
\int_{\calQ_{T}}(\partial\varphi_{\epsilon}(z_{\epsilon}))
{:}(\calJ_{\epsilon}z_{\epsilon})\dd x\dd t
\leq
\int_{\calQ_T}v{:}z\dd x\dd t,
\end{equation}
which is an immediate consequence of the convergence results for the
subsequences \((\partial\varphi_{\epsilon}(z_{\epsilon}))_{\epsilon>0}\)
and \((\calJ_{\epsilon} z_{\epsilon})_{\epsilon>0}\).
Hence, using the definitions of $w$ and $\varphi$, we may conclude 
that $z$ is a solution of \eqref{eq:5}--\eqref{eq:9} if 
we can prove that 
${\textrm{Proj}}_{\Erdev} (w (t, \cdot))- \Me \dot z (t, \cdot) \in  
\partial \Psi_{\textrm{dev}}(\dot z (t, \cdot))$ 
for almost every $t \in [0,T]$.
So we may use the same trick as previously and we only need to check that
\begin{equation} \label{eq:105}
\limsup_{\epsilon\rightarrow 0}
\int_{\calQ_{T}} {\rm Proj}_{\Erdev}  (w_{\epsilon}) {:}\dot z_{\epsilon}\dd x\dd t
\leq
\int_{\calQ_{T}}{\rm Proj}_{\Erdev} (w){:}\dot z\dd x\dd t.
\end{equation}
But $\dot z_{\epsilon}$ takes its values in $\Erdev$ and thus $\dot z $ 
takes also its values in $\Erdev$. It follows that \eqref{eq:105} is 
equivalent to 
\begin{equation*}
\limsup_{\epsilon\rightarrow 0}
\int_{\calQ_{T}} w_{\epsilon} {:}\dot z_{\epsilon}\dd x\dd t
\leq
\int_{\calQ_{T}}w{:}\dot z\dd x\dd t.
\end{equation*}
We compute 
\(\int_{\calQ_{T}}w_{\epsilon}{:}\dot z_{\epsilon}\dd x\dd t\)
and \(\int_{\calQ_{T}}w{:}\dot z\dd x\dd t\) from \eqref{eq:92}
and \eqref{eq:106}, then the convergence results obtained above imply
that \eqref{eq:105} holds if and only if
\begin{equation*} 
\liminf_{\epsilon\rightarrow 0}\int_{\calQ_{T}}
\partial\varphi_{\epsilon}(z_{\epsilon}){:}\dot z_{\epsilon}\dd x\dd t
\geq
\int_{\calQ_{T}}
\partial\varphi(z){:}\dot z\dd x\dd t.
\end{equation*}
We observe that
\(
\int_{\calQ_{T}}
\partial\varphi_{\epsilon}(z_{\epsilon}){:}\dot z_{\epsilon}\dd x\dd t
=\varphi_{\epsilon}(z_{\epsilon}(T))-\varphi_{\epsilon}(z^0)\geq
\varphi(\calJ_{\epsilon}
z_{\epsilon}(T))-\varphi_{\epsilon}(z^0)\)
and \(
\int_{\calQ_{T}}\partial\varphi(z){:}\dot z\dd x\dd t=
\varphi(z(T))-\varphi(z^0)\). 
Recalling that $z^0 \in \calD (\varphi)$, 
we get $\lim_{\epsilon \rightarrow 0} \varphi_{\epsilon} (z^0) = \varphi (z^0)$ and 
the lower semicontinuity of \(\varphi\),
allows us to conclude. This proves the existence result.

Finally we observe that 
\begin{equation*}
-{\textrm{Proj}}_{\Erdev} (w)={\textrm{Proj}}_{\Erdev} 
(g_1(\theta){+}g_2(z){+}\Ee z)+ \D_z H_1 (z)+\theta \D_z H_2 (z)+ v,
\end{equation*}
and, using the definition of the mappings $g_1$, $g_2$ and $\varphi$
\begin{equation*}
-{\textrm{Proj}}_{\Erdev} (w)= {\textrm{Proj}}_{\Erdev} 
(-\Ee( \ee (u){-} z))+ \D_z H_1 (z)+ \theta \D_z H_2 (z)-\nu \Delta z.
\end{equation*}
So we may rewrite 
\eqref{eq:ent_eq_2_2}  as follows
\begin{equation}
\label{eq:15bis} 
\dot z-\nu\Me^{-1}\Delta z=\Me^{-1}f^z.
\end{equation}
with
\(f^z\eqldef {\textrm{Proj}}_{\Erdev} (\Ee(\ee(u){-}z))-\D_z H_1(z)-
\theta\D_zH_2(z)-\psi\) and \(\psi= {\textrm{Proj}}_{\Erdev} (w) 
-\Me \dot z \in\partial\Psi_{\textrm{dev}} (\dot z)\).
With  \eqref{eq:Psi} we infer that 
\begin{equation*} 
\forall\psi\in\partial\Psi_{\textrm{dev}} (\dot z):\ 
\abs{\psi}\leq C^{\Psi} \ \text{ a.e. } \ (x,t)\in\Omega\times (0,T).
\end{equation*}
Furthermore, since $z \in \L^{\infty}(0,T; \H^1(\Omega))$, 
we infer with \eqref{eq:H2} that \(\D_z H_i(z)\)  belongs to
\(\L^{\infty}(0,T;\L^{p}(\Omega))\) for $i=1,2$. Then it follows that 
 \(\Me^{-1}f^z\) belongs to \(\L^{q}(0,T;\L^{2}(\Omega))\). 
Since $z^0 \in \X_{q,p} (\Omega)$, 
we may deduce from the maximal regularity result for parabolic systems
that \(z \in \L^{q}(0,T;\H^2(\Omega)) \cap \C^0([0,T]; \H^1(\Omega))\)
and \(\dot z \in \L^{q}(0,T;\L^{2}(\Omega))\). 
We refer to \cite{Dore91RADE,HieReh08QPMB,PruSch01SMRP} and the references 
therein 
for more details on the maximal regularity result for parabolic systems
and its consequences. 
\eproof

Let us observe that here neither $\partial \Psi_{\textrm{dev}}+\Me$ 
nor $\partial 
\varphi+ \D_z H_1+\theta \D_zH_2+{\textrm{Proj}}_{\Erdev}{\circ} 
(\Ee{+}g_1{+}g_2)$ are linear and self-adjoint 
operators of $\L^2 (\Omega; \Erdev)$, so we can not use the ideas proposed in 
\cite{ColVis90CDNE} 
to prove uniqueness. The uniqueness result proved below relies on the 
boundedness assumption \eqref{eq:H3} for the hardening functional $H_1$
combined with Gr\"onwall's lemma. 
More precisely 
let \(h_1 \in \C^2(\Erdev; \Er)\) be defined by 
\begin{equation}
\label{eq:68}
\forall z \in \Erdev: \ h_1(z)\eqldef H_1(z) - C^{H_1}\abs{z}^2, 
\end{equation}
with a  real number $ C^{H_1}>0$.
Assumption \eqref{eq:H3} implies that
there exists  $C^{h_1}>0$ such that 
\begin{equation}
\label{eq:67}  
\forall z_1, z_2 \in \Erdev : \ 
\abs{\D_z h_1(z_1){-}\D_z h_1(z_2)}
\leq C^{h_1}\abs{z_1{-}z_2}.
\end{equation}

\begin{proposition}[Uniqueness for \((\textrm{P}_{uz})\)]
\label{sec:uniqueness}
Assume that $\theta$ is given in 
$\L^q( 0,T;\L^p (\Omega))$, \eqref{eq:Psi}, 
\eqref{eq:H}, \eqref{eq:LM}  
and \eqref{eq:ell} hold and $u^0 \in \H^1( \Omega)$, 
$z^0 \in \X_{q,p}(\Omega)$. Then the problem 
\eqref{eq:ent_eq_1_1}--\eqref{eq:boun_cond2} 
admits a unique solution. 
\end{proposition}

\bproof 
Let \(\xi_1\eqldef(u_1,z_1)\) and \(\xi_2\eqldef(u_2,z_2)\)
be two solutions of \eqref{eq:ent_eq_1_2}--\eqref{eq:ent_eq_2_2} satisfying
\eqref{eq:init_cond2} and \eqref{eq:boun_cond2}. With the results of the
previous theorem we already know that $(u_i, z_i) \in \C^0 ([0,T];  
\H^1_0(\Omega) \times \H^1(\Omega))$  and $\Delta z_i \in \L^q (0,T; 
\L^2(\Omega))$, $i=1,2$.
Define
\begin{equation}
\label{eq:14} 
\begin{aligned} 
&\gamma(t)\eqldef \tfrac12\int_{\Omega}
\Ee((\ee(u_1){-}z_1){-}(\ee(u_2){-}z_2)){:}
((\ee(u_1){-}z_1){-}(\ee(u_2){-}z_2))
\dd x\\&
+\tfrac{\nu}{2} \int_{\Omega} \nabla (z_1{-}z_2){:} \nabla (z_1{-}z_2)
\dd x+ C^{H_1} \int_{\Omega} \abs{ z_1{-}z_2}^2 \dd x 
\end{aligned}
\end{equation}
for all  $t \in [0,T]$.
Since 
\begin{equation*}
\begin{aligned} 
&\gamma(t)=\tfrac12\int_{\Omega}
\Ee((\ee(u_1){-}z_1){-}(\ee(u_2){-}z_2)){:}
((\ee(u_1){-}z_1){-}(\ee(u_2){-}z_2))
\dd x\\&
-\tfrac{\nu}{2} \int_{\Omega} \Delta (z_1{-}z_2){\cdot} (z_1{-}z_2)
\dd x+ C^{H_1} \int_{\Omega} \abs{ z_1{-}z_2}^2 \dd x,
\end{aligned}
\end{equation*}
the mapping  $\gamma$ is continuous on $[0,T]$ and its derivative 
in the sense of distributions belongs to $\L^1(0,T)$. Thus $\gamma$ is  
absolutely 
continuous on $[0,T]$ and with assumptions \eqref{eq:H2} and  \eqref{eq:3} 
combined with  Korn's inequality, we infer that there exists a real number 
$\kappa>0$ such that 
\begin{equation} 
\label{eq:14bis} 
\gamma (t)\ge\kappa\bigl(\norm[\H^1(\Omega)]{u_1(\cdot,t){-}u_2(\cdot,t)}^2  
+ \norm[\H^1(\Omega)]{z_1(\cdot,t){-}z_2(\cdot,t)}^2\bigr) \ \forall  t\in[0,T].
\end{equation}
On the one hand, recalling  that \(\partial\Psi_{\textrm{dev}}(\cdot)\) 
is a monotone 
operator,
the Green's formula and \eqref{eq:boun_cond2} enable us to deduce from
\eqref{eq:ent_eq_2_2} that
\begin{equation}
\label{eq:12}
\begin{aligned}
\forall \tilde z_i \in \L^2 (\Omega; \Erdev): \ &\int_{\Omega}
(\Me\dot z_i{-}\Ee(\ee(u_i){-}z_i){+}\D_z H_1(z_i){+}
\theta\D_zH_2(z_i)){:} ( \dot z_i{-}{\tilde z}_i)\dd x \\&-
\nu\int_{\Omega}\Delta z_i{\cdot} ( \dot z_i{-}{\tilde z}_i)\dd x\leq 
\int_{\Omega} ( \Psi_{\textrm{dev}} ( \tilde z_i ) {-} \Psi_{\textrm{dev}} 
(\dot z_i)) \dd x.
\end{aligned}
\end{equation} 
On the other hand, we multiply \eqref{eq:ent_eq_1_2} by 
\(\dot u_i-{\tilde u}_i\),
we integrate this expression over \(\Omega\), and with the help of 
the Green's formula together with
\eqref{eq:boun_cond2}, we obtain
\begin{equation}
\label{eq:13}
\forall \tilde u_i \in \H^1_0(\Omega) : \  \int_{\Omega}
(\Ee(\mathrm{e}(u_i){-}z_i){+}\alpha\theta\bfI{+}\Le\mathrm{e}(\dot
u_i)){:}(\ee( \dot u_i)- \ee({\tilde u}_i))\dd x
=\int_{\Omega} \ell{\cdot}(\dot u_i - {\tilde u}_i)\dd x.
\end{equation}
Thus we add \eqref{eq:12} and \eqref{eq:13}, we get
\begin{equation}
\label{eq:10}
\begin{aligned}
&\int_{\Omega}  
(\Me\dot z_i{-}\Ee(\ee(u_i){-}z_i){+}\D_z H_1(z_i){+}
\theta\D_zH_2(z_i)){:}(\dot z_i{-}{\tilde z}_i)\dd x-
\nu\int_{\Omega}\Delta z_i{\cdot} (\dot{z}_i{-}\tilde z_i)\dd
x\\&
+\int_{\Omega}
(\Ee(\mathrm{e}(u_i){-}z_i){+}\alpha\theta\bfI{+}\Le\mathrm{e}(\dot
u_i)){:}(\ee( \dot u_i){-}\ee({\tilde u}_i))\dd x
\leq 
\int_{\Omega}(\ell{\cdot}(\dot u_i{-}{\tilde u}_i)
{+}\Psi_{\textrm{dev}}( \tilde z_i){-}\Psi_{\textrm{dev}}(\dot z_i)) \dd x.
\end{aligned}
\end{equation}
We choose now \((\tilde u_i,\tilde z_i)=(\dot u_{3-i}, \dot z_{3-i})\) 
for \(i=1,2\) in
\eqref{eq:10}, we add these two inequalities, and with the help of 
\eqref{eq:LM}, we find
\begin{equation*}
\begin{aligned}
&\dot \gamma (t)+c^{\Me}
\int_{\Omega}\abs{\dot z_1{-}\dot z_2}^2\dd x+c^{\Le}
\int_{\Omega}\abs{\ee(\dot u_1){-}\ee(\dot u_2)}^2\dd x
\\&
\leq 
- \int_{\Omega}(\D_z h_1(z_1){-}\D_z h_1(z_2)){:}(\dot  z_1{-}\dot z_2)\dd x 
- \int_{\Omega}\theta(\D_z H_2(z_1){-}\D_z H_2(z_2)){:}(\dot z_1{-} \dot z_2)
\dd x \ \text{ a.e. } \ t\in[0,T].
\end{aligned}
\end{equation*}
We estimate the terms of the right hand side of the previous 
inequality as follows:
\begin{equation*}
\begin{aligned}
&
\Bigl|\int_{\Omega}(\D_z h_1(z_1){-}\D_z h_1(z_2)){:}(\dot  z_1{-} \dot z_2)
\dd x 
\Bigr|\le C^{h_1}\int_{\Omega} \abs{z_1{-}z_2} \abs{ \dot z_1{-}\dot z_2}\dd
x 
\\&
\leq \tfrac{c^{\Me}}{4} \int_{\Omega} \abs{\dot z_1{-}\dot z_2}^2 \dd x + 
\tfrac{(C^{h_1})^2}{ c^{\Me}}  \int_{\Omega} \abs{z_1{-}z_2}^2 \dd x,
\end{aligned}
\end{equation*}
and 
\begin{equation*}
\begin{aligned}
&\Bigl| \int_{\Omega}\theta(\D_z H_2(z_1){-}\D_z H_2(z_2)){:}(\dot z_1{-}\dot
  z_2)\dd x\Bigr|
\leq C^{H_2}_{zz}  \int_{\Omega} \abs{\theta}\abs{z_1{-}z_2} 
\abs{\dot z_1{-}\dot z_2}\dd x \\&
\leq \tfrac{c^{\Me}}{4} \int_{\Omega} \abs{\dot z_1{-}\dot z_2}^2 \dd x 
+\tfrac{(C^{H_2}_{zz})^2}{ c^{\Me}}
\int_{\Omega} \abs{\theta}^2 \abs{z_1{-}z_2}^2 \dd x 
\leq  
\tfrac{c^{\Me}}{4} \int_{\Omega}\abs{\dot z_1{-}\dot z_2}^2 \dd x  
+\tfrac{(C^{H_2}_{zz})^2}{c^{\Me}} 
\norm[\L^4(\Omega)]{\theta}^2 \norm[\L^4(\Omega)]{z_1{-}z_2}^2.
\end{aligned}
\end{equation*}
Since $\H^1 (\Omega)\hookrightarrow\L^4 (\Omega)$ with continuous embedding, 
we infer that there exists $C>0$ such that 
\begin{equation*}
\dot \gamma (t) \leq C \bigl(1{+}\norm[\L^4(\Omega)]{\theta (\cdot,t)}^2\bigr) 
\norm[\H^1 (\Omega)]{z_1(\cdot,t){-}z_2(\cdot,t)}^2 \leq \tfrac{C}{\kappa} 
\bigl(1{+}\norm[\L^4(\Omega)]{\theta(\cdot,t)}^2 \bigr) 
\gamma(t) \ \text{ a.e. } \ t\in[0,T].
\end{equation*}
But $\theta \in \L^q(0,T;\L^p (\Omega))$ with $q \ge 4 $ and $p \ge 4$, 
thus we get
\begin{equation*}
\displaystyle \forall t \in [0,T]: \
\gamma (t) \leq \gamma (0) 
\ee^{\frac{C}{\kappa}\int_0^t  
(1{+}\norm[\L^4(\Omega)]{\theta(\cdot,s)}^2 )\dd s},
\end{equation*}
which, allows us to conclude.
\eproof

We provide that \(\tilde \vartheta\mapsto(u,z)\) is a continuous mapping from
\(\L^{\bar q}(0,T;\L^{\bar p}(\Omega))\) into
\(\H^1(0,T;\H^1_0(\Omega){\times}\L^2(\Omega))
\cap\L^{\infty}(0,T;\H^1_0(\Omega){\times}\H^1(\Omega))\) 
where $(u,z)$ is the unique solution of $(\textrm{P}_{uz})$ 
when $\theta= \zeta (\tilde \vartheta)$.

\begin{lemma}
\label{lm:cont_uz}
Assume that \eqref{eq:Psi}, \eqref{eq:H}, \eqref{eq:LM} and  
\eqref{eq:ell} 
hold and that $u^0 \in \H^1(\Omega)$, $z^0 \in \X_{q,p}(\Omega)$. 
Then \(\vartheta\mapsto(u,z)\) 
is continuous from
\(\L^{\bar q}(0,T;\L^{\bar p}(\Omega))\) into
\(\H^1(0,T;\H^1_0(\Omega){\times}\L^2(\Omega))
{\cap}\L^{\infty}(0,T;\H^1_0(\Omega){\times}\H^1(\Omega))\).
\end{lemma}

\bproof
We consider \(\vartheta_i\in \L^{\bar q}(0,T;\L^{\bar p}(\Omega))\) and 
for $i=1,2$, we
denote by \(\theta_i\eqldef\zeta (\vartheta_i) \in \L^q(0,T;\L^p(\Omega))\) and 
\((u_i,z_i)\) the solution of the
following system:
\begin{subequations}
\label{eq:61}
\begin{align}
&-\dive(\Ee(\mathrm{e}(u_i){-}z_i){+}\alpha\theta_i\bfI{+}\Le\mathrm{e}(\dot
u_i))
=\ell,\label{eq:62}\\&
\partial\Psi(\dot z_i)+\Me\dot z_i-\Ee(\ee(u_i){-}z_i)+\D_z H_1(z_i)+
\theta\D_zH_2(z_i)-\nu\Delta z_i\ni 0,\label{eq:63}
\end{align}
\end{subequations}
together with initial conditions
\begin{equation}
\label{eq:64}
u_i(0,\cdot)=u^0,\quad
z_i(0,\cdot)=z^0,
\end{equation}
and boundary conditions
\begin{equation}
\label{eq:65}
u_{i|_{\partial\Omega}}=0,\quad
\nabla z_i{\cdot}\eta_{|_{\partial\Omega}}=0.
\end{equation}
Since the mapping $\phi_1: 
\vartheta \mapsto \theta = \zeta (\vartheta)$ 
is continuous from $\L^{\bar q}(0,T;\L^{\bar p}(\Omega))$ to 
$\L^{\bar q}(0,T;\L^{\bar p}(\Omega))$, it is enough to check that 
the mapping $\theta = \zeta (\vartheta) \rightarrow (u,z)$ is continuous
from $\L^{\bar q}(0,T;\L^{\bar p}(\Omega))$ to 
$\L^{\bar q}(0,T;\L^{\bar p}(\Omega))$.

Once again the key tool is the structural decomposition 
\eqref{eq:68}  combined with  Gr\"onwall's lemma. 
We reproduce the same kind of computations as in   
Proposition \ref{sec:uniqueness}. 
More precisely, the total energy inequality associated to 
\eqref{eq:61}--\eqref{eq:65} is obtained by
multiplying \eqref{eq:62} and \eqref{eq:63} by \(\dot u_{3{-}i}{-}\dot u_i\)
and \(\dot z_{3{-}i}{-}\dot z_i\) respectively, and integrating 
over $\Omega$. Then, we add these two inequalities, we find
\begin{equation}
\label{eq:66}
\begin{aligned}
&
\int_{\Omega}
\Ee(\ee(u_i){-}z_i){:}((\ee(\dot u_{3{-}i}){-}\dot z_{3{-}i}){-}
(\ee(\dot u_i){-}\dot z_i))\dd x
+
\int_{\Omega}
\Me\dot z_i{:}(\dot z_{3{-}i}{-}\dot z_i)\dd x \\&
+
\int_{\Omega}
\alpha\theta_i
\tr(\ee(\dot u_{3{-}i}){-}\ee(\dot u_i))\dd x
+
\int_{\Omega}
\Le\ee(\dot u_i){:}
(\ee(\dot u_{3{-}i}){-}\ee(\dot u_i))\dd x \\&
-
\int_{\Omega}
\nu\Delta z_i{\cdot}(\dot z_{3{-}i}{-}\dot
z_i)\dd x
+
\int_{\Omega}
\D_z H_1(z_i){:}(\dot z_{3{-}i}{-}\dot z_i)\dd x \\&
+
\int_{\Omega} 
\theta_i\D_zH_2(z_i){:}(\dot z_{3{-i}}{-}\dot
z_i)\dd x
-
\int_{\Omega}
\ell{\cdot}(\dot u_{3{-}i}{-}\dot u_i)\dd x+
\Psi_{\textrm{dev}} (\dot z_{3{-i}})-\Psi_{\textrm{dev}} 
(\dot z_i)\geq 0 \ \text{ a.e. } 
\ t\in[0,T].
\end{aligned}
\end{equation}
It is convenient to introduce the notations:
\(\bar u\eqldef u_1{-}u_2\), \(\bar z\eqldef z_1{-}z_2\)
and \(\bar\theta\eqldef \theta_1{-}\theta_2\). Therefore, we 
take \(i=1,2\) in \eqref{eq:66} and we add these two inequalities, we get
\begin{equation*}
\begin{aligned}
&
\int_{\Omega}
\Ee(\ee(\bar u){-}\bar z){:}(\ee(\dot{\bar u}){-}\dot{\bar z})\dd x
+
\int_{\Omega}
\Me\dot{\bar z}{:}\dot{\bar z}\dd x
+
\int_{\Omega}
\Le\ee(\dot{\bar u}){:}\ee(\dot{\bar u})\dd x
- 
\int_{\Omega}
\nu\Delta {\bar z}{\cdot}\dot{\bar z}\dd x\\&
+
\int_{\Omega}
(\D_z H_1(z_1){-}\D_z H_1(z_2)){:}\dot{\bar z}\dd x
\leq -
\int_{\Omega}
\alpha\bar\theta
\tr(\ee(\dot{\bar u}))\dd x
-
\int_{\Omega}
(\theta_1\D_zH_2(z_1){-}\theta_2\D_zH_2(z_2)){:}
\dot{\bar z}\dd x.
\end{aligned}
\end{equation*}
Define  
\begin{equation*}
\begin{aligned}
\forall t \in [0,T]: \
& \gamma(t)\eqldef\tfrac12 \int_{\Omega}
\Ee((\ee(u_1){-}z_1){-}(\ee(u_2){-}z_2)){:}
((\ee(u_1){-}z_1){-}(\ee(u_2){-}z_2))
\dd x\\&- 
\tfrac{\nu}{2} \int_{\Omega}\nabla (z_1{-}z_2){:} \nabla (z_1{-}z_2)
\dd x+ C^{H_1} \int_{\Omega} \abs{ z_1{-}z_2}^2 \dd x.
\end{aligned}
\end{equation*}
As in Proposition \ref{sec:uniqueness}, we can check that 
$\gamma$ is absolutely continuous on $[0,T]$ and 
 \eqref{eq:LM} and \eqref{eq:68} 
imply that 
\begin{equation*}
\begin{aligned}
&
\dot \gamma (t) 
+c^{\Me}
\norm[\L^2(\Omega)]{\dot{\bar z}}^2
+c^{\Le}
\norm[\L^2(\Omega)]{\ee(\dot{\bar u})}^2 
\le 
- \int_{\Omega} (\D_z h_1(z_1){-}\D_z h_1(z_2)){:}\dot{\bar z}\dd x\\& 
-\int_{\Omega}
\alpha\bar\theta
\tr(\ee(\dot{\bar u}))\dd x
-\int_{\Omega}(\theta_1\D_zH_2(z_1){-}\theta_2\D_zH_2(z_2)){:}\dot{\bar z}\dd x,
\end{aligned}
\end{equation*}
for almost every \(t\in [0,T]\).
Clearly, it follows from \eqref{eq:67}  that
\begin{equation*}
\begin{aligned}
&
\dot \gamma (t)+c^{\Me}
\norm[\L^2(\Omega)]{\dot{\bar z}}^2
+c^{\Le}
\norm[\L^2(\Omega)]{\ee(\dot{\bar u})}^2
\leq-\int_{\Omega}
\alpha\bar\theta
\tr(\ee(\dot{\bar u}))\dd x
\\&
-\int_{\Omega}
(\theta_1\D_zH_2(z_1){-}\theta_2\D_zH_2(z_2)){:}\dot{\bar z}\dd
x+C^{h_1}\int_{\Omega} \abs{\bar z}\abs{\dot{\bar z}}\dd x.
\end{aligned}
\end{equation*}
We estimate the first and third term on the right hand side
with the help of  Cauchy-Schwarz's inequality, while for 
the second term, we use the following decomposition 
\begin{equation*}
(\theta_1\D_zH_2(z_1){-}\theta_2\D_zH_2(z_2)){:}\dot{\bar z}=
(\bar\theta\D_zH_2(z_1){+}\theta_2(\D_zH_2(z_1){-}\D_zH_2(z_2))){:}\dot{\bar z}.
\end{equation*}
Hence, we obtain
\begin{equation}
\label{eq:71}
\begin{aligned}
&\dot \gamma (t) 
+\tfrac{3c^{\Me}}4
\norm[\L^2(\Omega)]{\dot{\bar z}}^2
+\tfrac{c^{\Le}}2
\norm[\L^2(\Omega)]{\ee(\dot{\bar u})}^2 
\leq 
\tfrac{(C^{h_1})^2}{c^{\Me}}
\norm[\L^2(\Omega)]{\bar z}^2+
\tfrac{3 \alpha^2}{2c^{\Le}}
\norm[\L^2(\Omega)]{\bar \theta}^2
\\&+\int_{\Omega}
\big(\abs{\bar\theta}\abs{\D_zH_2(z_1)}\abs{\dot{\bar z}}{+}
\abs{\theta_2}\abs{\D_zH_2(z_1){-}\D_zH_2(z_2)}\abs{\dot{\bar z}}\bigr)\dd x.
\end{aligned}
\end{equation}
It remains to estimate the last term on the right hand side of
\eqref{eq:71}. We use \eqref{eq:H3} and  \eqref{eq:H2} to get
\begin{equation*}
\begin{aligned}
&\int_{\Omega} 
\big(\abs{\bar\theta}\abs{\D_zH_2(z_1)}\abs{\dot{\bar z}}{+}
\abs{\theta_2}\abs{\D_zH_2(z_1){-}\D_zH_2(z_2)}\abs{\dot{\bar z}}\bigr)\dd
x
\\&\leq 
C_z^{H_2}
\int_{\Omega}
\bigl(1{+}\abs{z_1}\bigr)\abs{\bar\theta}\abs{\dot{\bar z}}
\dd x
+C_{zz}^{H_2}
\int_{\Omega} 
\abs{\theta_2}\abs{\bar z}\abs{\dot{\bar z}}\dd x.
\end{aligned}
\end{equation*}
The Young's inequality implies that there exists \(\gamma_i>0\),
\(i=1,2,3\), such that
\begin{equation*}
\begin{aligned}
&\int_{\Omega}
\big(\abs{\bar\theta}\abs{\D_zH_2(z_1)}\abs{\dot{\bar z}}{+}
\abs{\theta_2}\abs{\D_zH_2(z_1){-}\D_zH_2(z_2)}\abs{\dot{\bar z}}\bigr)\dd x
\leq 
\tfrac{C_z^{H_2}}{2\gamma_1}
\norm[\L^2(\Omega)]{\bar\theta}^2 \\&+
\tfrac{C_z^{H_2}}{2\gamma_2}
\int_{\Omega}
\abs{\bar\theta}^2\abs{z_1}^2\dd x +
\tfrac{C_{zz}^{H_2}}{2\gamma_3}
\int_{\Omega}
\abs{\theta_2}^2\abs{\bar z}^2\dd x
+\tfrac{C_z^{H_2}(\gamma_1{+}\gamma_2){+}C_{zz}^{H_2}\gamma_3}2
\norm[\L^2(\Omega)]{\dot{\bar z}}^2.
\end{aligned}
\end{equation*}
Since $z_1 \in \L^q (0,T; \H^2 (\Omega))$ and $\H^2(\Omega) 
\hookrightarrow \L^{\infty}(\Omega)$ with continuous embedding, we get
\begin{equation}
\label{eq:73}
\begin{aligned}
&\int_{\Omega} 
\big(\abs{\bar\theta}\abs{\D_zH_2(z_1)}\abs{\dot{\bar z}}{+}
\abs{\theta_2}\abs{\D_zH_2(z_1){-}\D_zH_2(z_2)}\abs{\dot{\bar z}}\bigr)\dd x
\leq \tfrac{C_z^{H_2}}{2\gamma_1}
\norm[\L^2(\Omega)]{\bar\theta}^2 \\&
+
\tfrac{C_z^{H_2}}{2 \gamma_2}
\norm[\L^{\infty}(\Omega)]{z_1}^2
\norm[\L^2(\Omega)]{\bar\theta}^2+
\tfrac{C_{zz}^{H_2}}{2\gamma_3}
\norm[\L^4(\Omega)]{\theta_2}^2
\norm[\L^4(\Omega)]{\bar z}^2
+\tfrac{ C_z^{H_2}(\gamma_1{+}\gamma_2){+}C_{zz}^{H_2}\gamma_3}2
\norm[\L^2(\Omega)]{\dot{\bar z}}^2.
\end{aligned}
\end{equation}
We insert \eqref{eq:73} in \eqref{eq:71} and we choose
\(\gamma_1=\gamma_2=\tfrac{c^{\Me}}{4C_z^{H_2}}\) and 
\(\gamma_3=\tfrac{c^{\Me}}{2C_{zz}^{H_2}}\). Therefore the continuous embedding 
\(\H^1(\Omega)\hookrightarrow\L^4(\Omega)\) and \eqref{eq:14bis} give  
\begin{equation}
\label{eq:74}
\begin{aligned}
\dot \gamma (t) 
+\tfrac{c^{\Me}}4
\norm[\L^2(\Omega)]{\dot{\bar z}}^2
+\tfrac{c^{\Le}}2
\norm[\L^2(\Omega)]{\ee(\dot{\bar u})}^2 
\leq C(\bar\theta,z_1)+
C(\theta_2)
\norm[\H^1(\Omega)]{\bar z}^2 \leq  C(\bar\theta,z_1)+
\tfrac{C(\theta_2)}{\kappa} \gamma (t),
\end{aligned}
\end{equation}
for almost every $t \in [0,T]$. Here \(C(\bar\theta,z_1)\eqldef
\bigl(\tfrac{(3\alpha)^2}{2c^{\Le}}{+}\tfrac{2(C_z^{H_2})^2}{c^{\Me}}\bigr)
\norm[\L^2(\Omega)]{\bar \theta}^2 +
\tfrac{2 (C_z^{H_2})^2 }{c^{\Me}}
\norm[\L^{\infty}(\Omega)]{z_1}^2
\norm[\L^2(\Omega)]{\bar\theta}^2\)
and \(C(\theta_2)\eqldef
\tfrac{(C^{h_1})^2}{c^{\Me}}{+}\tfrac{C^2(C_{zz}^{H_2})^2}{c^{\Me}}
\norm[\L^4(\Omega)]{\theta_2}^2\) where $C>0$ 
is the generic constant involved in the continuous embedding  
$\H^1(\Omega)\hookrightarrow\L^4(\Omega)$.
Since $\bar q \ge 3$ we can check that $C(\bar \theta, z_1) \in \L^1(0,T)$ and
\begin{equation*}
\int_0^T C(\bar \theta, z_1) \dd t \le 
\bigl(\tfrac{(3\alpha)^2}{2c^{\Le}}{+}\tfrac{2(C_z^{H_2})^2}{c^{\Me}}\bigr)
\norm[\L^2(0,T; \L^2(\Omega))]{\bar \theta}^2
+ \tfrac{2 (C_z^{H_2})^2 }{c^{\Me}}
\norm[\L^{\frac{2 \bar q}{\bar q -2}}(0,T; \L^{\infty}(\Omega))]{z_1}^2
\norm[\L^{\bar q} (0,T; \L^2(\Omega))]{\bar\theta}^2.
\end{equation*}
Note that $C(\theta_2)\in\L^1(0,T)$, which, thanks to  Gr\"onwall's lemma, gives
\begin{equation*}
\forall t \in [0,T]: \
\gamma (t) \leq \int_0^t C(\bar \theta(\cdot,s), z_1(\cdot,s)) 
\ee^{\frac{1}{\kappa} \int_s^t C(\theta_1(\cdot,\sigma){-}\bar \theta(\cdot,\sigma)) 
\dd \sigma}\dd s.
\end{equation*}
Recalling that the mapping $\tilde \vartheta \mapsto \theta$ is 
Lipschitz continuous from $\L^{\bar q}(0,T; \L^{\bar p}(\Omega))$ 
to $\L^{\bar q}(0,T; \L^{\bar p}(\Omega))$ and maps any bounded subset of 
$\L^{\bar q}(0,T; \L^{\bar p}(\Omega))$ into a bounded subset of 
$\L^q(0,T; \L^p(\Omega))$ with $\bar p=2$ and $\bar q>4$, $p \ge 4$, 
$q = \beta_1 \bar q>8$, the last estimate 
 allows us to conclude.
\eproof

Let us observe furthermore that the image of a bounded set of 
\(\L^{\bar q}(0,T;\L^{\bar p}(\Omega))\) by the mapping
$\tilde \vartheta \mapsto \zeta (\tilde \vartheta) = \theta \mapsto (u,z)$ 
is a bounded subset of 
\(\H^1 (0,T ; \H^1_0(\Omega){\times}\L^2 (\Omega))\cap
\L^{\infty}(0,T;\H^1_0(\Omega){\times}\H^1(\Omega))\).

Let us introduce now some new notations. For any $r>1$, let 
\begin{equation*}
 \V^{r}(\Omega; \Er^3)\eqldef \bigl\{u\in\L^2(\Omega; \Er^3): \ \nabla u\in
\L^{r}(\Omega; \Er^{3 \times 3}) \bigr\},
\end{equation*}
and for any $r \ge 2$, let
\begin{equation*}
\V_0^{r}(\Omega; \Er^3)
\eqldef \bigl\{ u \in \V^r(\Omega; \Er^3): \  u_{|_{\partial\Omega}}=0 \bigr\}.
\end{equation*}
We endowed $\V^r (\Omega; \Er^3)$ with the following norm
\begin{equation*}
\forall u \in \V^r(\Omega; \Er^3): \
\norm[\V^{r}(\Omega)]{ u}\eqldef\norm[\L^2 (\Omega)]{u} 
+ \norm[\L^r (\Omega)]{\nabla u}.
\end{equation*}

The aim of the next two lemmas is to prove further regularity results for 
the solutions $(u,z)$ of the system composed by the momentum equilibrium
equation and the flow rule. More precisely, 
assuming that $\theta$ remains in a bounded subset of  
\(\L^{q}(0,T;\L^{p}(\Omega))\), 
we will prove that $\ee (\dot u)$, $\dot z$ and $z$ remain 
in a bounded subset of  \(\L^{q}(0,T;\L^{p}(\Omega))\), 
\(\L^{q}(0,T;\L^{ p/2}(\Omega)) \cap \L^{q/2}(0,T;\L^{p}(\Omega))\) 
and \(\L^{q}(0,T;\H^2(\Omega))\), respectively.

\begin{lemma}
\label{lm:reg_e}
Assume that \eqref{eq:Psi}, \eqref{eq:H}, \eqref{eq:46}, \eqref{eq:LM}
and \eqref{eq:ell} hold. Assume moreover that $z^0 \in \X_{q,p}(\Omega)$ and 
$u^0 \in \V^p_0(\Omega;\Er^3)$.
Then \(\mathrm{e}(u)\) 
belongs to \(\W^{1,q}(0,T;\L^p(\Omega))\)
and \(\theta\mapsto \ee( u)\)
maps  
any bounded subset of \(\L^{{q}}(0,T;\L^{ p}(\Omega))\) into a 
bounded subset of \(\W^{1,{q}}(0,T;\L^p(\Omega))\).
\end{lemma}

\bproof
The idea of the proof is to interpret \eqref{eq:ent_eq_1_2} as an
ODE for $u$ in an appropriate Banach space. 
More precisely, let ${\mathcal F}_p\eqldef\L^2(\Omega;\Er^3) 
{\times}\L^p (\Omega; \Ers)$ endowed with the norm
\begin{equation*}
\forall \varphi =( \varphi_1, \varphi_2) \in {\mathcal F}_p : \ 
\norm[{\mathcal F}_p]{\varphi} 
\eqldef \norm[{\mathcal F}_p]{(\varphi_1, \varphi_2)} = 
\norm[\L^2(\Omega)]{\varphi_1} + \norm[\L^p(\Omega)]{\varphi_2}.
\end{equation*}
It follows that ${\mathcal F}_p$ is a Banach space. Let us introduce now the 
mapping $\calA_{\Ee}$
\begin{equation*}
 \begin{aligned}
\calA_{\Ee}:\ 
\V^p_0 (\Omega;\Er^3) &\rightarrow {\mathcal F}_p,\\
u&\mapsto \varphi\eqldef (0,\Ee\ee(u)).
\end{aligned} 
\end{equation*}
Since $\Ee\in\L^{\infty}(\Omega)$, we infer that $\calA_{\Ee}$
is a linear continuous mapping from $\V^p_0(\Omega)$ 
to ${\mathcal F}_p$. Besides since $\Le$ is a symmetric, 
positive definite  tensor, classical results about PDE in 
Banach spaces imply that, 
for all $\varphi =(\varphi_1,\varphi_2)\in{\mathcal F}_p$, there exists a 
unique $u\in\V_0^p(\Omega;\Er^3)$, denoted $u= \Lambda_p (\varphi)$, such 
that 
\begin{equation*}
\forall v\in {\calD} (\Omega;\Er^3):\
\int_{\Omega}\Le\ee(u){:}\ee(v)\dd x=\int_{\Omega} 
\varphi_1{\cdot}v\dd x+\int_{\Omega} \varphi_2{:}\ee(v)\dd x,
\end{equation*}
where $p^*$ is the conjugate of $p$, i.e $\frac{1}{p^*}+\frac{1}{p}=1$.
Furthermore there exists a real number $C>0$, independent of 
$\varphi$, such that
\begin{equation*}
\norm[\V^p (\Omega; \Er^3)]{u}
\leq C\bigl(\norm[\L^2 (\Omega)]{\varphi_1}  
{+}\norm[\L^p(\Omega)]{\varphi_2}\bigr)=C\norm[{\mathcal F}_p]{\varphi},
\end{equation*} 
and $\Lambda_p$ is linear continuous from ${\mathcal F}_p$ 
to $\V^p_0(\Omega; \Er^3)$
(for more details, the reader is referred to 
\cite{Vale88LTEU}).
It follows that 
\eqref{eq:ent_eq_1_2} can be rewritten as 
\begin{equation}
\label{eq:11}
\dot u=\calG_p(\varphi_{z\theta},u),
\end{equation}
with $\varphi_{z \theta}\eqldef(\ell,\Ee z{-}\alpha\theta \bfI)$ 
and \(\calG_p(\varphi_{z\theta},u)\eqldef
\Lambda_p (\varphi_{z\theta}{-}\calA_{\Ee}u) = 
\Lambda_p (\varphi_{z\theta}) 
{-}\Lambda_p ( \calA_{\Ee}u)  \).

With the assumption \eqref{eq:ell} and the previous results, we already know 
that $\varphi_{z \theta}\in\L^{q }(0,T;{\mathcal F}_p)$ and we can apply 
classical results for ODE in Banach spaces to conclude that 
$u\in\W^{1,q}(0,T;\V^p_0(\Omega; \Er^3))$, the reader is referred to
\cite{Car90} for more details.

We can also obtain estimates for $u$ and $\dot u$ in $\V^p_0(\Omega;\Er^3)$.
To this aim, we introduce the
following notations:
\(
C_{\Lambda_p}\eqldef 
\norm[\calL({\mathcal F}_p,\V^p_0(\Omega))]{\Lambda_p}\)
and \(C_{\calA_{\Ee}}\eqldef 
\norm[\calL(\V^p_0(\Omega),{\mathcal F}_p)]{\calA_{\Ee}}\).
Then, we observe that \eqref{eq:11} gives
\begin{equation}
\label{eq:48}
\norm[\V^p(\Omega)]{\dot u(\cdot,t)}
\leq
C_{\Lambda_p}
\bigl(
\norm[{\mathcal F}_p]{\varphi_{z\theta}(\cdot,t){-}\calA_{\Ee}u^0}
{+}
C_{\calA_{\Ee}}
\norm[\V^p(\Omega)]{u(\cdot,t){-}u^0}
\bigr) \ \text{ a.e. } \ t \in [0,T].
\end{equation}
Let us turn now to the term \(\norm[\V^p(\Omega)]{u(\cdot,t){-}u^0}\).
It is clear that
\begin{equation*}
\begin{aligned}
&\norm[\V^p(\Omega)]{u(\cdot,t){-}u^0}
\leq
\int_0^t \norm[\V^p(\Omega)]{\calG_p(\varphi_{z\theta}(\cdot,s),u(\cdot,s))}\dd s
\\&\leq
C_{\Lambda_p}
\int_0^t
\norm[{\mathcal F}_p]{\varphi_{z\theta}(\cdot,s){-}\calA_{\Ee}u^0}\dd s+
C_{\Lambda_p} C_{\calA_{\Ee}}
\int_0^t\norm[\V^p(\Omega)]{u(\cdot,s){-}u^0}\dd s.
\end{aligned}
\end{equation*}
Therefore, we may infer from  Gr\"onwall's lemma and 
the continuous embedding $\H^1(\Omega)\hookrightarrow\L^p (\Omega)$
that there exists a generic constant \(C_1>0\) such that
\begin{equation}
\label{eq:49}
\begin{aligned}
&\norm[\V^p(\Omega)]{u(\cdot,t){-}u^0}
\leq \ee^{C_{\Lambda_p} C_{\calA_{\Ee}}t} C_{\Lambda_p}
\int_0^t
\norm[{\mathcal F}_p]{\varphi_{z\theta}(s){-}\calA_{\Ee}u^0}\dd s 
\\&\leq \ee^{C_{\Lambda_p} C_{\calA_{\Ee}}t} C_{\Lambda_p}C_1
\int_0^t
\bigl(\norm[\L^2 (\Omega)]{\ell(s)}{+}\norm[\L^{\infty}(\Omega)]{\Ee} 
\norm[\H^1 (\Omega)]{z(\cdot,s)}{+}\alpha 
\norm[\L^p(\Omega)]{\theta(\cdot,s)}   
{+}C_{\calA_{\Ee}}\norm[\V^p(\Omega)]{u^0}\bigr)\dd s.
\end{aligned}
\end{equation}
We insert \eqref{eq:49} in \eqref{eq:48}, 
we find that there exists \(C_2>0\)
such that
\begin{equation*}
\norm[\L^q(0,T;\V^p(\Omega))]{\dot u}
\leq
C_2\bigl( \norm[\L^{q}(0,T;\L^2 (\Omega))]{\ell} 
{+}\norm[\L^{\infty}(\Omega)]{\Ee} 
\norm[\L^{q} (0,T;\H^1 (\Omega))]{z}{+}\alpha 
\norm[\L^{q} (0,T;\L^p(\Omega))]{\theta}{+}
\norm[\V^p(\Omega)]{u^0}\bigr).
\end{equation*}
Recalling Lemma \ref{lm:cont_uz} we infer that 
 the image of any bounded set of 
\(\L^{\bar{q}}(0,T;\L^{\bar p} (\Omega))\) by the mappings 
$\tilde \vartheta \mapsto \zeta(\tilde \vartheta) = \theta 
\mapsto \mathrm{e}(u)$ and 
$\tilde \vartheta \mapsto \zeta(\tilde \vartheta) = \theta 
\mapsto \mathrm{e}(\dot u)$ 
are still bounded sets in $\L^{\infty} (0,T; \L^p (\Omega))$ and 
\(\L^q(0,T;\L^p(\Omega))\), respectively.
\eproof

Let us conclude this section with some regularity results for \(z\) and 
$\dot z$. 

\begin{lemma}
\label{lm:reg_z}
Assume that \eqref{eq:Psi}, \eqref{eq:H}, \eqref{eq:46}, \eqref{eq:LM}
and \eqref{eq:ell} hold. Assume moreover that $z^0 \in \X_{q,p}(\Omega)$ 
and $u^0 \in \V^p_0(\Omega; \Er^3)$.
Then $\dot z$ and $\Delta z$ belong to  
\(\L^{{q/2}}(0,T;\L^p(\Omega)) {\cap}\L^{{q}}(0,T;\L^{p/2}(\Omega))
\) and $z \in \C^0([0,T], \X_{q,p}(\Omega)) \cap \L^q(0,T; \H^2(\Omega))$.
Moreover \(\theta\mapsto (\dot z, \Delta z, z) \)  maps  any bounded subset of 
\(\L^{{q}}(0,T;\L^p(\Omega))\) into a bounded subset of 
\((\L^{{q/2}}(0,T;\L^p(\Omega)) 
{\cap}\L^{{q}}(0,T;\L^{p/2}(\Omega)))^2
\times(\C^0([0,T]; \X_{q,p}(\Omega)){\cap}\L^q(0,T; \H^2(\Omega)))\).
\end{lemma}

\bproof
Let us rewrite again 
\eqref{eq:ent_eq_2_2} as
\begin{equation}
\label{eq:15}
\dot z-\nu\Me^{-1}\Delta z=\Me^{-1}f^z.
\end{equation}
with
\(f^z\eqldef {\textrm{Proj}}_{\Erdev} (\Ee(\ee(u){-}z)) -\D_z H_1(z)-
\theta\D_zH_2(z)-\psi\) and \(\psi\eqldef {\textrm{Proj}}_{\Erdev} (w ) 
- \Me \dot z \in\partial\Psi_{\textrm{dev}}(\dot z)\)
where \(w\) has been defined in the proof of Theorem \ref{sec:existence}. 
With the assumption \eqref{eq:Psi.bdd}, we have
\begin{equation*}
\forall \psi\in\partial\Psi_{\textrm{dev}} (\dot z):\
\abs{\psi}\leq C^{\Psi} \ \text{ a.e. } \ (x,t)\in\Omega\times (0,T).
\end{equation*}
Then Lemma \ref{lm:reg_e} enables us to infer that
\({\textrm{Proj}}_{\Erdev} (\Ee\ee(u))-\psi\) 
remains bounded in 
\(\L^{\infty}(0,T;\L^p(\Omega))\). Furthermore, we know with Lemma 
\ref{lm:cont_uz} that \(z\) is bounded in \(\L^{\infty}(0,T;\H^1(\Omega))\), 
so, using \eqref{eq:H2}, we infer that \(\D_z H_i(z)\) is bounded in 
\(\L^{\infty}(0,T;\L^p(\Omega))\) for $i=1,2$. Then it follows that 
 \(\Me^{-1}f^z\) is bounded in \(\L^{q}(0,T;\L^{p/2}(\Omega))\)
if $\theta$ belongs to a bounded subset of \(\L^{q}(0,T;\L^{p}(\Omega))\). 
We may deduce from the maximal regularity result for parabolic systems
that \(z\) is bounded in \(\L^{q}(0,T;\H^2(\Omega))\)
and \(\dot z\) is bounded in \(\L^{q}(0,T;\L^{p/2}(\Omega))\)
(see \cite{Dore91RADE,HieReh08QPMB,PruSch01SMRP}).

Since \(\H^2(\Omega)\hookrightarrow\L^{\infty}(\Omega)\) 
with continuous embedding, 
\(z\) is bounded in \(\L^{q}(0,T;\L^{\infty}(\Omega))\) and 
thus \(\theta \D_z H_2(z)\) is bounded in \(\L^{ q/2}(0,T;\L^{p}(\Omega))\). 
We may deduce that 
\(f^z\) is bounded in \(\L^{ q/2}(0,T;\L^p(\Omega))\) and, the 
maximal regularity result for parabolic systems allows us to infer that
\(\dot z\) and \(\Delta z\) belong to a bounded subset of 
\(\L^{q}(0,T;\L^{p/2}(\Omega))\cap \L^{q/2}(0,T;\L^p(\Omega))\) 
and \(z\) belongs to a bounded subset of \(\C^0([0,T]; \X_{q,p}(\Omega))\)
whenever $\theta$ belongs to a bounded subset of 
\(\L^{q}(0,T;\L^{p}(\Omega))\).
\eproof

\section{Existence and regularity results for the 
enthalpy equation}
\label{sec:enthalpy-equation}

In this section existence and uniqueness results for the enthalpy equation
are recalled and some regularity results are obtained.
More precisely, let us consider the
enthalpy equation $(\textrm{P}_{\vartheta})$: 
\begin{equation}
\label{eq:19}
\dot\vartheta-\dive(\tilde \kappa^c\nabla\vartheta)= f
\end{equation}
together with initial conditions
\begin{equation}
\label{eq:20}
\vartheta(0)=\vartheta^0,
\end{equation}
and boundary conditions
\begin{equation}
\label{eq:21}
\tilde \kappa^c\nabla\vartheta{\cdot}\eta_{|_{\partial\Omega}}=0.
\end{equation}
We assume that the initial enthalpy \(\vartheta^0\)
belongs to \(\L^2(\Omega)\) and $f$
belongs to \(\L^2(0,T;\L^2(\Omega))\). 
Furthermore, we assume that $\tilde \kappa^c\in\L^{\infty}(\calQ_{T}; \Ers)$
and satisfies
\begin{subequations}
\label{eq:41bis}
\begin{align}
\exists c^{\kappa^c}>0\
\forall v \in \Er^3 :\ 
\tilde \kappa^c(x,t) v {\cdot} v \ge c^{\kappa^c} \abs{v}^2 \ 
\text{ a.e. } \ (x,t)\in\calQ_{T}, \\
\exists C^{\kappa^c}>0:\
\abs{\tilde \kappa^c (x,t)} \le C^{\kappa^c} \ \text{ a.e. } \ (x,t)\in\calQ_{T}.
\end{align}
\end{subequations}
The weak formulation of the problem is given by
\begin{equation}
\label{eq:28}
\begin{cases}
\text{Find }\vartheta: [0,T] \rightarrow  
\H^1(\Omega) \text{ such that } \vartheta (0) = \vartheta^0 \text{ and 
for all } \xi \in \H^1(\Omega),\\ 
\displaystyle{\int_{\Omega}\dot\vartheta \xi \dd x
+\int_{\Omega} \tilde \kappa^c\nabla\vartheta{\cdot}\nabla\xi \dd x =
\int_{\Omega}
f  \xi \dd x} \ \text{ in the sense of distributions.}
\end{cases}
\end{equation}

\begin{theorem}[Existence and uniqueness for $(\textrm{P}_{\vartheta})$]
\label{sec:exist-uniq-results}
Under the previous assumptions, the problem 
\eqref{eq:19}--\eqref{eq:21} possesses  a unique
solution \(\vartheta\in \C^0([0,T];\L^2(\Omega))\cap 
\L^2(0,T;\H^1(\Omega))\) with $\dot \vartheta \in \L^2(0,T; (\H^1(\Omega)')$.
Moreover we have
\begin{equation*}
\forall \tau \in [0,T]:\
\norm[\L^2(\Omega)]{\vartheta (\tau)}^2 
+ 2 c^{\kappa^c} \int_0^{\tau} \norm[\L^2(\Omega)]{\nabla\vartheta (t)}^2 \dd
t \leq \ee^{\tau} 
\bigl(\norm[\L^2(\Omega)]{\vartheta^0}^2{+}\norm[\L^2(0,T; \L^2(\Omega))]{f
}^2\bigr).
\end{equation*}
\end{theorem}

\bproof
The proof of existence and uniqueness of 
a solution is quite classical
and can be found in \cite{Brez83AFTA,RavTho83}. 
The estimate is straightforward and its verification is left to the reader.
\eproof

Let us introduce the following functional space
\begin{equation*}
{\calW}\eqldef\bigl\{ \vartheta \in  \L^2(0,T;\H^1(\Omega)) 
\cap \L^{\infty} (0,T; \L^2(\Omega)): \ 
\dot\vartheta \in \L^2(0,T; (\H^1(\Omega)') \bigr\},
\end{equation*}
endowed with the norm
\begin{equation*}
\forall \vartheta \in {\calW} : \ 
\norm[{\mathcal W}]{\vartheta}  \eqldef  
 \norm[\L^2(0,T;\H^1(\Omega))]{\vartheta}+ \norm[\L^{\infty}(0,T; 
\L^2(\Omega))]{\vartheta}
+ \norm[\L^2(0,T; (\H^1(\Omega)')]{\dot \vartheta} .
\end{equation*}
Due to \cite{Sim87CSLP}, we know that ${\mathcal W}$ is compactly 
embedded in $\L^{\bar q} (0,T; \L^{\bar p} (\Omega))$.
Note that the previous estimate implies that there exists a 
generic constant $C>0$ such that the solution of problem 
$(\textrm{P}_{\vartheta})$ satisfies
\begin{equation*}
\norm[{\mathcal W}]{\vartheta}\leq  
C
\bigl(\norm[\L^2(\Omega)]{\vartheta^0} 
{+} \norm[\L^2(0,T; \L^2(\Omega))]{f}\bigr). 
\end{equation*}

\section{Local existence result}
\label{sec:existence-result}

This section is dedicated to the proof of a local existence result for 
\eqref{eq:ent_eq1}--\eqref{eq:init_cond1}
by using a  fixed-point argument.
To this aim, for any given $\tilde \vartheta 
\in\L^{\bar q}(0,T; \L^{\bar p} (\Omega))$, we consider  
$\tilde \kappa^c\eqldef \kappa^c (\ee(u), z, \theta)$ and
 \(f= f^{\tilde \vartheta}\eqldef 
\Le\ee(\dot u){:}\ee(\dot u)
+ \theta ( \alpha\tr(\ee(\dot u))+\D_zH_2(z){:}\dot z) +\Psi(\dot z)+
\Me \dot z{:}\dot z\) in $(\textrm{P}_{\vartheta})$, where 
$(u,z)$ are the solutions of $(\textrm{P}_{uz})$ with 
$\theta = \zeta (\tilde \vartheta)$. With the results obtained 
in the Section \ref{sec:const-equat}, we already know that 
$f^{\tilde \vartheta}$ belongs to $\L^{q/4} (0,T; \L^{p/2} (\Omega))$. 
Since $p \ge 4$ and $q> 8$, we infer that $f^{\tilde \vartheta}$ 
belongs to $\L^{2} (0,T; \L^2(\Omega))$ and we can define 
$\vartheta\in\C^0 ([0,T];\L^2 (\Omega))\cap{\mathcal W}$ 
as the unique solution of $(\textrm{P}_{\vartheta})$. 
Thus we can introduce the fixed point mapping 
$\phi: \tilde \vartheta \mapsto \vartheta$ from 
$\L^{\bar q} (0,T; \L^{\bar p} (\Omega))$ to $\L^{\bar q} (0,T; \L^{\bar p} 
(\Omega))$. 

\begin{proposition}
\label{cont-phi}
The  mapping $\phi$ is continuous 
from $\L^{\bar q} (0,T; \L^{\bar p} (\Omega))$ to 
$\L^{\bar q} (0,T; \L^{\bar p} (\Omega))$.
\end{proposition}

\bproof
Let $(\tilde \vartheta_n)_{n \in \en}$ be  a converging sequence 
of $\L^{\bar q} (0,T; \L^{\bar p} (\Omega) )$ and let $\tilde 
\vartheta_*$ be its limit. We denote by 
$\vartheta_n\eqldef\phi(\tilde \vartheta_n)$ for all 
$n\ge 0$ and $\vartheta_* \eqldef\phi(\tilde \vartheta_*)$. 
Since $(\tilde \vartheta_n)_{n \in \en}$ is a bounded family 
of $\L^{\bar q} (0,T; \L^{\bar p} (\Omega))$, 
we infer from the previous results that $(\vartheta_n)_{n \in \en}$ 
is bounded in $\C^0( [0,T];\L^2 (\Omega))\cap{\mathcal W}$. 
We may deduce that $(\vartheta_n)_{n \in \en}$ is relatively compact 
in $\L^{\bar q} (0,T;\L^{\bar p} (\Omega))$ (see \cite{Sim87CSLP}) 
and there exists a subsequence, still denoted $(\vartheta_n)_{n \in \en}$, 
such that
\begin{equation*}
\begin{aligned}
\vartheta_n &\rightharpoonup \vartheta \ \text{ in } \ \L^2(0,T;\H^1(\Omega))
\ \text{ weak},\\
\vartheta_n &\rightarrow \vartheta \ \text{ in } \ 
\L^{\bar q}(0,T;\L^{\bar p}(\Omega)).
\end{aligned}
\end{equation*}
Let us define \(\calV_T\eqldef\bigl{\{}w\in {\C^{\infty}}([0,T]):\ w(T) =0\bigr{\}}\).
Hence we observe that for all $n \ge 0$, we have
\begin{equation}
\label{cont-phi-1}
\begin{aligned}
\forall\xi\in\H^1(\Omega) \ \forall w \in\calV_T: \
&-\int_{\calQ_T} \vartheta_n (x,t)  \xi (x) \dot w  (t) \dd x \dd t 
+ \int_{\calQ_T} \tilde \kappa^c_n \nabla \vartheta_n (x,t) \nabla \xi (x) 
w(t) \dd x
\dd t \\
&  
= \int_{\calQ_T} f^{\tilde \vartheta_n} (x,t) \xi (x) w (t) \dd x \dd t
+ \int_{\Omega} \vartheta^0 (x) \xi(x) w(0) \dd x,
\end{aligned}
\end{equation}
with $\tilde \kappa^c_n\eqldef\kappa^c \bigl( \ee(u_n),z_n,\theta_n)$ and 
$(u_n,z_n)$ solutions of  $(\textrm{P}_{uz})$ with 
$\theta_n\eqldef\zeta (\tilde \vartheta_n)$. Since 
$(\tilde \vartheta_n)_{n \in \en}$ converges to 
$\tilde \vartheta_*$ in $\L^{\bar q}(0,T;\L^{\bar p}(\Omega))$, 
we infer from Lemma \ref{lm:cont_uz} that $(u_n,z_n)_{n \in \en}$ 
converges to the solution $(u_*, z_*)$ of $(\textrm{P}_{uz})$ 
with $\theta_* = \zeta (\tilde \vartheta_*)$ in 
$\H^1(0,T;\H^1_0(\Omega){\times}\L^2 (\Omega)) 
\cap\L^{\infty}(0,T;\H^1_0(\Omega){\times}\H^1(\Omega))$. 
Let us recall that the mapping 
$\phi_1:\tilde \vartheta\mapsto\theta=\zeta(\tilde \vartheta)$ 
is Lipschitz continuous from $\L^{\bar q} (0,T;\L^{\bar p}(\Omega))$ to 
$\L^{\bar q}(0,T;\L^{\bar p}(\Omega))$, which implies, 
possibly after extracting another subsequence, that
\begin{equation*}
\theta_n, u_n, z_n \rightarrow 
\theta_*, u_*, z_* \ \text{ a.e. } \ (x,t)\in\calQ_T.
\end{equation*}
Note that the continuity of the mapping $\kappa^c$ gives
\begin{equation*}
\tilde \kappa^c_n=\kappa^c(\ee(u_n),z_n,\theta_n) 
\rightarrow \tilde \kappa^c_*=\kappa^c( \ee(u_*),z_*,\theta_*) 
\ \text{ a.e. } \ (x,t)\in\calQ_T,
\end{equation*}
and due to the boundedness assumption on $\kappa^c$, we obtain 
with the Lebesgue's theorem that
\begin{equation*}
\tilde \kappa^c_n \nabla \xi  w \rightarrow \tilde \kappa^c_* \nabla \xi w 
\ \text{ in } \ \L^2(0,T; \L^2(\Omega)).
\end{equation*}
Therefore it is possible to pass to the limit in all the terms of the left 
hand side of \eqref{cont-phi-1} to get
\begin{equation*}
\begin{aligned}
\forall\xi\in\H^1(\Omega)\ \forall w \in \calV_T:\ 
&-\int_{\calQ_T} \vartheta (x,t)  \xi (x) \dot w  (t) \dd x \dd t 
+ \int_{\calQ_T} \tilde \kappa^c_* \nabla \vartheta (x,t) \nabla \xi (x) 
w(t) \dd x
\dd t \\
& = \lim_{n \rightarrow + \infty} 
\int_{\calQ_T} f^{\tilde \vartheta_n} (x,t)\xi(x)w(t)\dd x\dd t
+ \int_{\Omega} \vartheta^0(x) \xi(x) w(0) \dd x.
\end{aligned}
\end{equation*}
Recalling that $p \in [4,6]$ and $q > 8$,  we infer that the 
mapping $\tilde \vartheta \mapsto f^{\tilde \vartheta}$ 
is continuous from $\L^{\bar q}(0,T;\L^{\bar p} (\Omega))$ 
to $\L^{r_1} (0,T;\L^{r_2} (\Omega))$ with 
$\tfrac{1}{r_1}= \frac{3}{q} + \tfrac12$
and $\tfrac1{r_2}= \frac{1}{p} + \frac{1}{2}$.
Indeed, for any $\tilde\vartheta_i$ in $\L^{\bar q} (0,T;\L^{\bar p}(\Omega))$, 
let $(u_i,z_i)$ be the solution of $(\textrm{P}_{uz})$ with 
$\theta_i = \zeta (\tilde \vartheta_i)$, $i= 1,2$, we find
\begin{equation*}
\begin{aligned}
& f^{\tilde \vartheta_1}-f^{\tilde \vartheta_2}= 
\Le\ee(\dot u_1{+}\dot u_2){:}\ee(\dot u_1{-}\dot u_2)+
(\theta_1{-}\theta_2)(\alpha\tr(\ee(\dot u_1)){+}\D_z H_2(z_1){:} 
\dot z_1)
\\&
+ \theta_2(\alpha 
\tr(\ee(\dot u_1{-}\dot u_2)){+}\D_zH_2(z_1){:}\dot z_1{-}\D_zH_2(z_2){:}
\dot z_2) 
+\Psi(\dot z_1)-\Psi(\dot z_2)+\Me (\dot z_1{+}\dot z_2){:}(\dot z_1{-}
\dot z_2).
\end{aligned}
\end{equation*}
On the other hand, \eqref{eq:4} and \eqref{eq:Psi.bdd} imply
\begin{equation*}
\abs{\Psi(\dot z_1){-}\Psi(\dot z_2)} \le C^{\Psi} \abs{\dot z_1{-}\dot z_2},
\end{equation*}
and  \eqref{eq:H2} and \eqref{eq:H3} 
give
\begin{equation*}
\begin{aligned}
\abs{\D_zH_2(z_1){:}\dot z_1{-}\D_zH_2(z_2){:}\dot z_2 }
&\leq\abs{\D_zH_2(z_1)}\abs{\dot z_1{-}\dot z_2}+
\abs{\D_zH_2(z_1){-}\D_zH_2(z_2) } 
\abs{\dot z_2}
\\&
\leq C_z^{H_2} (1{+}\abs{z_1}) 
\abs{\dot z_1{-}\dot z_2}+C^{H_2}_{zz}\abs{z_1{-}z_2}\abs{\dot z_2}. 
\end{aligned}
\end{equation*}
The boundedness and the continuity properties proved
in Lemma \ref{lm:reg_e} and Lemma \ref{lm:cont_uz}, respectively, 
allow us to deduce the desired result. Therefore, we may infer that 
\begin{equation*}
\begin{aligned}
\forall\xi\in\H^1(\Omega) \ \forall w \in \calV_T:\
&\lim_{n \rightarrow + \infty}\int_{\calQ_T}f^{\tilde \vartheta_n}(x,t) 
\xi(x)w(t)\dd x\dd t=\int_{\calQ_T}f^{\tilde \vartheta_*}(x,t) 
\xi(x)w (t)\dd x\dd t.
\end{aligned}
\end{equation*}
We conclude that $\vartheta$ is solution of problem $(\textrm{P}_{\vartheta})$ 
with 
the data $\tilde \kappa^c_*$ and $f^{\tilde \vartheta_*}$. Moreover by 
uniqueness of 
the solution, it follows that $\vartheta=\vartheta_*$ and the 
whole sequence $(\vartheta_n)_{n \in \en}$ converges to 
$\vartheta_* = \zeta (\tilde \vartheta_*)$.
\eproof

We establish now that the mapping $\phi$ fulfills 
the other assumptions of the Schauder's fixed point theorem. To this aim, we
introduce some notations: 
let $R^0, R^{\vartheta}>0$ be any given positive real numbers 
such that
$\max \bigl(\norm[\V^p (\Omega)]{u^0} ,  
\norm[\X_{q,p}(\Omega)]{z^0} \bigr) \leq R^0$ and 
$\norm[\L^{\bar q}(0,T; \L^{\bar p}(\Omega))]{\tilde \vartheta} 
\leq R^{\vartheta}$. Clearly, we have
\begin{equation*}
\begin{aligned}
& \norm[\L^q(0,T;\L^p(\Omega))]{\zeta(\tilde \vartheta)} = 
\norm[\L^q(0,T;\L^p(\Omega))]{\theta} 
\leq 
\bigl(\tfrac{\beta_1}{c^c}\bigr)^{\frac1{\beta_1}} 
\abs{\Omega}^{\frac{\beta_1 \bar p - p}{\beta_1 p \bar p}}
\norm[\L^{\bar q}(0,T;\L^{ \bar p}(\Omega))]{\tilde\vartheta}^{\frac1{\beta_1}} 
\leq R^{\theta}\eqldef\bigl(\tfrac{\beta_1}{c^c} 
R^{\vartheta}\bigr)^{\frac1{\beta_1}}  
\abs{\Omega}^{\frac{\beta_1 \bar p{-}p}{\beta_1 p \bar p}}.
\end{aligned}
\end{equation*}
Therefore once again the results of Section 
\ref{sec:const-equat} are used, which imply that
there exists a constant
$R^{f}\eqldef R^{f} (R^0,R^{\theta},\norm[{\C^0([0,T];\L^2 (\Omega))}]{\ell})>0$, 
depending only on $R^0$, $R^{\theta}$ and 
$\norm[{\C^0([0,T];\L^2 (\Omega))}]{\ell}$, such that
\begin{equation*}
\norm[\L^{q/4}(0,T;\L^{p/2}(\Omega))]{f^{\tilde \vartheta}} 
\leq R^{f}(R^0,R^{\theta},\norm[{\C^0([0,T];\L^2(\Omega))}]{\ell}).
\end{equation*}
The results of Section \ref{sec:enthalpy-equation} imply that there exists a 
generic constant $C>0$ such that
\begin{equation*}
\begin{aligned}
&\norm[\L^{\infty} (0,T; \L^{2} (\Omega))]{ \vartheta} \leq
 C  (\norm[\L^{2} (0,T; \L^{2} (\Omega))]{ f^{\tilde \vartheta}} {+} 
\norm[\L^2(\Omega)]{\vartheta^0}) 
\leq C \bigl(|\Omega|^{\frac{p-4}{2p}}T^{\frac{q-8}{2q}} 
\norm[\L^{q/4} (0,T ;\L^{p/2}(\Omega))]{f^{\tilde \vartheta}}  
{+}\norm[\L^2(\Omega)]{\vartheta^0}\bigr) \\&
 \leq C \bigl( |\Omega|^{\frac{p-4}{2p}}   
T^{\frac{q-8}{2q}} R^{f}(R^0, R^{\theta}, 
\norm[{\C^0([0,T]; \L^2 (\Omega))}]{\ell})
{+}\norm[\L^2(\Omega)]{\vartheta^0} \bigr).
\end{aligned}
\end{equation*}

Now let $0 < \tau \le T$ and let us introduce the following
functional space
\begin{equation*}
{\calW}_{\tau}\eqldef\bigl\{
\vartheta \in  \L^2(0, \tau;\H^1(\Omega)) \cap  \L^{\infty}(0, \tau; 
\L^2(\Omega)): \ \dot \vartheta \in \L^2(0, \tau; (\H^1(\Omega))') \bigr\}.
\end{equation*}
For any $\tilde \vartheta \in 
\L^{\bar q} (0, \tau;\L^{\bar p} (\Omega))$, we define its 
extension 
$\tilde \vartheta_{\text{ext}}$ by $\tilde \vartheta_{\text{ext}} = \tilde 
\vartheta $ on $[0, \tau]$ and $\tilde \vartheta_{\text{ext}} =0$ on $(\tau,T]$. 
It is clear that $\tilde \vartheta_{\text{ext}} \in 
\L^{\bar q} (0, T; \L^{\bar p} (\Omega))$ and the mapping $\tilde \vartheta 
\mapsto \tilde \vartheta_{\text{ext}}$ is a contraction from 
$\L^{\bar q} (0, \tau; \L^{\bar p} (\Omega))$ into $\L^{\bar q} (0, T; 
\L^{\bar p} (\Omega))$. For any $\tilde \vartheta \in 
\L^{\bar q} (0, \tau; \L^{\bar p} (\Omega))$, we define $\phi_{\tau} 
(\tilde \vartheta)$ as the restriction on $[0, \tau]$ of $\phi(\tilde
\vartheta_{\text{ext}})$. 
We infer immediately from Proposition \ref{cont-phi} that $\phi_{\tau}$ 
is continuous from $\L^{\bar q} (0, \tau; \L^{\bar p} (\Omega))$ to 
$\L^{\bar q} (0, \tau; \L^{\bar p} (\Omega))$. 
Furthermore, for any $\tilde \vartheta \in \L^{\bar q} (0, \tau; \L^{\bar p}
(\Omega))$, 
we have 
\begin{equation*}
\norm[\L^{\bar q} (0, \tau; \L^{\bar p} (\Omega))]{\phi_{\tau} 
(\tilde \vartheta)} = 
\norm[\L^{\bar q} (0, \tau; \L^{2} (\Omega))]{\phi_{\tau} (\tilde \vartheta)}
= \Bigl( \int_0^\tau \norm[\L^2(\Omega)]{\phi(\tilde 
\vartheta_{\text{ext}}(\cdot,t))}^{\bar q} 
\dd t\Bigr)^{\frac1{\bar q}}\le \tau^{\frac1{\bar q}}  
\norm[\L^{\infty} (0,T; \L^2(\Omega))]{\phi(\tilde \vartheta_{\text{ext}})},
\end{equation*}
and the previous estimates allow us to show that, for any $R^{\vartheta} >0$, 
there exists $\tau \in (0,T]$ such that
$\phi_{\tau}$ maps the closed  ball 
${\bar B}_{\L^{\bar q} (0,\tau;\L^{\bar p}(\Omega))}(0,R^{\vartheta})$ 
into itself. Note that the image of 
${\bar B}_{\L^{\bar q} (0,T;\L^{\bar p}(\Omega))}(0,R^{\vartheta})$ by $\phi$  is 
a bounded subset of ${\calW}$ and thus it is relatively compact 
in $\L^{\bar q}(0,T;\L^{\bar p}(\Omega))$. It follows that the image of 
${\bar B}_{\L^{\bar q} (0, \tau ;\L^{\bar p}(\Omega))}(0,R^{\vartheta})$ by 
$\phi_{\tau}$  
is also relatively compact 
in $\L^{\bar q}(0,\tau ;\L^{\bar p}(\Omega))$.
Consequently, we may conclude that the problem 
\eqref{eq:ent_eq1}--\eqref{eq:init_cond1} 
possesses a local solution $(u,z, \vartheta)$ defined 
on $[0,\tau]$ such that
\(u\in\W^{1,q}(0,\tau;\V^p_0 (\Omega))\),
\(z\in\L^{\infty}(0,\tau;\H^1(\Omega))\cap\H^1(0,\tau;\L^2(\Omega)) 
\cap \C^0([0,\tau]; \X_{q,p} (\Omega)) \cap \L^q(0, \tau; \H^2(\Omega)) \), 
\(\dot z, \Delta z \in \L^{q/2} (0,\tau;\L^p(\Omega))
\cap\L^{q}(0, \tau; \L^{p/2} (\Omega))\) and
\(\vartheta \in {\calW}_{\tau} \). 

We have to go back to the problem \eqref{eq:ent_eq}--\eqref{eq:boun_cond}. 
First we observe that $g$ and $\zeta$ define a $\C^1$-diffeo\-morphism from 
$(0,\infty)$ to $(0,\infty)$ and  
any solution of \eqref{eq:ent_eq1}--\eqref{eq:init_cond1} provides 
a solution of \eqref{eq:ent_eq}--\eqref{eq:boun_cond}
as soon as the enthalpy $\vartheta$ remains strictly positive.
 So we assume now that the initial enthalpy 
is strictly positive almost everywhere on $\Omega$, 
i.e., there exists $\bar\vartheta>0$ such that
\begin{equation}
\label{init_enthalpy}
g(\theta^0 (x))=\vartheta^0(x) 
\ge\bar\vartheta>0 \ \text{ a.e. }\ x\in\Omega.
\end{equation}
Therefore it is possible to use the Stampacchia's truncation method 
and to prove a  local existence result for the problem 
\eqref{eq:ent_eq}--\eqref{eq:boun_cond}.

\begin{theorem}[Local existence result]
\label{local_existence}
 Assume that \eqref{eq:Psi}, \eqref{eq:H}, \eqref{eq:47}, \eqref{eq:LM},  
\eqref{eq:ell} and \eqref{eq:36}
hold. Then, for any initial data 
$u^0\in\V^p_0(\Omega)$, 
$z^0\in\X_{q,p}(\Omega)$ and $\vartheta^0\in\L^2(\Omega)$ satisfying 
\eqref{init_enthalpy}, there exists $\tau\in(0,T]$ such that the problem 
\eqref{eq:ent_eq}--\eqref{eq:boun_cond} admits a solution on $[0, \tau]$.
\end{theorem}

\bproof
Let $(u,z, \vartheta)$ be a solution 
of \eqref{eq:ent_eq1}--\eqref{eq:init_cond1} on 
$\bigl[0,\tau\bigr]$. 
We prove now that
\begin{equation*}
\vartheta (x,t) >0 \ \text{ a.e. } \ (x,t)\in\calQ_{\tau}.
\end{equation*}
To do so, we introduce some notations: 
let \(C\eqldef\tfrac{(3 \alpha)^2}{2 c^{\Le}}+\tfrac{(C_z^{H_2})^2}{c^{\Me}}\)
and let $\varphi: [0,\tau]\rightarrow\Er$ such that
\begin{equation}
\label{eq:16}
\forall t\in[ 0, \tau]:\
\varphi(t)\eqldef\bar\vartheta \ee^{- \tfrac{\beta_1}{c^c} \int_0^t 
(C{+}\frac{(C_z^{H_2})^2}{c^{\Me}}\norm[\L^{\infty}
  (\Omega)]{z(\cdot,s)}^2)\dd s}.
\end{equation}
We use the classical Stampacchia's truncation method. Let us 
define \(G\in\C^1(\Er)\) satisfying
\begin{enumerate}[(i)]
\item \(\forall \sigma\in\Er\ \exists C^G>0:\ \abs{G'(\sigma)}\leq C^G\),
\item \(G\) is strictly increasing on \((0,\infty)\),
\item \(\forall \sigma\leq 0:\ G(\sigma)=0\).
\end{enumerate}
Let us define also \(H(\sigma)\eqldef\int_0^{\sigma}G(s)\dd s\) 
for all $\sigma\in\Er$, $\vartheta_1\eqldef-\vartheta+\varphi$ and 
\(h(t)\eqldef\int_{\Omega}H(\vartheta_1)\dd x\). Clearly, 
we have $H\in\C^2 (\Er;\Er)$ and $H(\sigma)>0$ for all $\sigma>0$. 
Furthermore, $\vartheta_1(0)=-\vartheta^0 + \bar\vartheta\leq 0$ 
almost everywhere on $\Omega$ implies that $h(0)=0$. Since   
$\vartheta\in{\mathcal W}_{\tau}$ and $\varphi\in\H^1(0,\tau;\Er)$, 
we infer that $h$ is absolutely continuous and 
\begin{equation*}
\begin{aligned}
&\dot h(t)= 
\int_{\Omega} G(\vartheta_1)\dot\vartheta_1\dd x
\\&= {-}\int_{\Omega} G(\vartheta_1)\bigl(\dive(\kappa^c \nabla \vartheta){+}
\Le\ee(\dot u){:}\ee(\dot u)
{+}\theta(\alpha\tr(\ee(\dot u)){+}\D_zH_2(z){:}\dot z){+}\Psi(\dot z){+}
\Me \dot z{:}\dot z{-}\dot \varphi\bigr)\dd x\\&
=-\int_{\Omega} G' (\vartheta_1)\kappa^c \nabla
\vartheta_1 {:}\nabla\vartheta_1\dd x 
-\int_{\Omega} G(\vartheta_1)\bigl( 
\Le\ee(\dot u){:}\ee(\dot u)
{+} \theta (\alpha\tr(\ee(\dot u)){+}\D_zH_2(z){:}\dot z) {+}\Psi(\dot z){+}
\Me \dot z{:}\dot z{-}\dot \varphi \bigr)\dd x,
\end{aligned}
\end{equation*}
for almost every $t\in[0,\tau]$.
It follows from \eqref{eq:H2}, \eqref{eq:LM} and  Cauchy-Schwarz's 
inequality that 
\begin{equation*}
\Le\ee(\dot u){:}\ee(\dot u)
{+}\alpha \theta \tr(\ee(\dot u)) \ge c^{\Le} 
\abs{\ee(\dot u)}^2 -3\alpha \abs{\theta} \abs{\ee(\dot u)} 
\ge \tfrac{c^{\Le}}{2}\abs{\ee(\dot u)}^2 -\tfrac{(3\alpha)^2 
\abs{\theta}^2}{2 c^{\Le}},
\end{equation*}
and 
\begin{equation*}
\Me \dot z{:}\dot z + \theta \D_zH_2(z){:}\dot z \ge c^{\Me} \abs{\dot z}^2 - 
\abs{\theta} \abs{\D_zH_2(z){:}\dot z} \geq c^{\Me} \abs{\dot z}^2 
- C_z^{H_2} \abs{\theta}
(1{+}\abs{z}) \abs{\dot z} 
\geq \tfrac{c^{\Me}}{2} \abs{\dot z}^2 -
\tfrac{(C_z^{H_2})^2 \abs{\theta}^2}{c^{\Me}}(1{+}\abs{z}^{2}) .
\end{equation*}
Since $G'(\vartheta_1) \ge 0$ and $G(\vartheta_1) \geq 0$ 
almost everywhere and \eqref{eq:Psi.bdd} and \eqref{eq:41} hold, we get
\begin{equation*}
\dot h(t)\leq\int_{\Omega}G (\vartheta_1)\bigl(\abs{\theta}^2 \bigl(C{+}
\tfrac{(C_z^{H_2})^2}{c^{\Me}}  
\abs{z}^{2 }\bigr) {+}\dot\varphi\bigr) \dd x  \ \text{ a.e. } \ t\in 
[0,\tau].
\end{equation*}
But $\theta =\zeta (\vartheta)$, and reminding that $\beta_1 \ge 2$, we have 
with \eqref{eqlp:1} that
\begin{equation*}
\abs{\theta} = \abs{\zeta(\vartheta)}
\le 
\sqrt{\tfrac{\beta_1}{c^c}\vartheta^+{+}1} -1 
\le \sqrt{\tfrac{\beta_1}{c^c}\vartheta^+} \ 
\text{ a.e. } \ (x,t) \in \calQ_{\tau}.\end{equation*} 
On the other hand, $G(\vartheta_1)$ vanishes whenever  $\vartheta \ge 
\varphi$, it follows that
\begin{equation*}
\dot h(t)\leq\int_{\Omega}G (\vartheta_1)
\bigl(\tfrac{\beta_1}{c^c} \varphi \bigl(C{+}
\tfrac{(C_z^{H_2})^2}{c^{\Me}}  
\abs{z}^{2 }\bigr){+}\dot\varphi\bigr)\dd x  \le 0 \ \text{ a.e. } \ t\in 
[0,\tau].
\end{equation*}
We may deduce that $h(t)\leq h(0)=0$ for all $t\in[0,\tau]$. Then we infer 
that 
\begin{equation*}
H(\vartheta_1) = 0 \ \text{ a.e. } \ 
(x,t)\in\Omega\times (0,\tau),
\end{equation*}
which implies that
\begin{equation*}
\vartheta_1=-\vartheta+\varphi\leq 0 \ \text{ a.e. } \
(x,t)\in \Omega\times (0,\tau). 
\end{equation*}
This concludes the proof.
\eproof

\section{Global existence result}
\label{sec:global-existence-result}

We begin this section with some a priori estimates for the solutions 
of the problem \eqref{eq:ent_eq1}--\eqref{eq:init_cond1}. 
As usual, the result relies  on an energy balance combined with  
Gr\"onwall's lemma. Then the  global existence result 
is proved by using a contradiction
argument together with the results obtained in the previous sections.

\begin{proposition}[Global energy estimate]
\label{thm:existence-result-global}
Assume that \eqref{eq:Psi}, \eqref{eq:H}, \eqref{eq:47}, \eqref{eq:LM},
\eqref{eq:ell} and \eqref{eq:36} hold. Assume moreover that 
$u^0 \in\V^p_0 (\Omega)$, $z^0 \in \X_{q,p} (\Omega)$, $\vartheta^0 
\in \L^2(\Omega)$ such that  \eqref{init_enthalpy} holds. 
Then, there exists a constant $\tilde C>0$,  
depending only on $\norm[\H^1(\Omega)]{u^0}$, $\norm[\H^1(\Omega)]{z^0}$, 
$\norm[\L^1(\Omega)]{\vartheta^0}$ and the data such that
for any solution $(u,z, \vartheta) $ of problem 
\eqref{eq:ent_eq1}--\eqref{eq:init_cond1} 
defined on $[0,\tau]$, $\tau \in (0,T]$,  we have
\begin{equation*}
\begin{aligned}
\forall \tilde \tau \in [0, \tau]:\
\norm[\H^1(\Omega)]{u(\cdot,\tilde \tau)}^2
+\norm[\H^1(\Omega)]{z(\cdot,\tilde \tau)}^2 
+ \norm[\L^{1}(\Omega)]{\vartheta(\cdot,\tilde \tau)}
\leq \tilde C.
\end{aligned}
\end{equation*}
\end{proposition}

\bproof
On the one hand, we multiply \eqref{eq:ent_eq1_11}
by \(\dot u\) and we integrate this expression over 
\(\calQ_{\tilde \tau}\), with $\tilde \tau \in [0,  \tau]$, to
get
\begin{equation}
\label{eq:52}
\int_{\calQ_{\tilde \tau}}
(\Ee(\mathrm{e}(u){-}z){+}\alpha\theta\bfI{+}\Le\mathrm{e}(\dot
u)){:}\ee(\dot u)\dd x\dd t
=\int_{\calQ_{\tilde \tau}}\ell{\cdot}\dot u\dd x\dd t.
\end{equation}
On the other hand, by using the definition of the subdifferential 
\(\partial\Psi(\dot z)\), we deduce from \eqref{eq:ent_eq1_12} that
\begin{equation}
\label{eq:54}
\int_{\calQ_{\tilde \tau}}
(\Me\dot z{-}\Ee(\ee(u){-}z){+}\D_z H_1(z){+}
\theta\D_zH_2(z){-}\nu\Delta z){:}
\dot z\dd x\dd t
+\int_{\calQ_{\tilde \tau}}\Psi(\dot z)\dd x\dd t=0.
\end{equation}
Adding \eqref{eq:52} and \eqref{eq:54}, we obtain
\begin{equation}
\label{eq:55}
\begin{aligned}
&\tfrac12\int_{\Omega}\Ee(\ee(u(\cdot,\tilde \tau)){-}z(\cdot,\tilde
\tau)){:}(\ee(u(\cdot,\tilde \tau)){-}z(\cdot,\tilde \tau))\dd x
+\tfrac{\nu}2\norm[\L^2(\Omega)]{\nabla z(\cdot,\tilde \tau)}^2
+\int_{\calQ_{\tilde \tau}}\Me\dot z{:}\dot z\dd x\dd t\\&+
\int_{\calQ_{\tilde \tau}}\Le\ee(\dot u){:}\ee(\dot u)\dd x\dd t+
\int_{\Omega} H_1(z(\cdot,\tilde \tau))\dd x + 
\int_{\calQ_{\tilde \tau}}\theta(\alpha 
\tr(\ee(\dot u)){+}\D_zH_2(z){:}\dot z)\dd
x\dd t\\&+
\int_{\calQ_{\tilde \tau}}\Psi(\dot z)\dd x\dd t=C_0^{u,z}
+\int_{\calQ_{\tilde \tau}}\ell{\cdot}\dot u\dd x\dd t.
\end{aligned}
\end{equation}
where \(C_0^{u,z}\eqldef
\tfrac12\int_{\Omega}\Ee(\ee(u^0){-}z^0){:}
(\ee(u^0){-}z^0)\dd x+
\tfrac{\nu}2\norm[\L^2(\Omega)]{\nabla z^0}^2 + \int_{\Omega} 
H_1(z^0) \dd x\).
We integrate \eqref{eq:ent_eq1_13} over \(\calQ_{\tilde \tau}\),  
by taking into account the boundary conditions
\eqref{eq:boun_cond1}, we find 
\begin{equation*}
\begin{aligned}
&\int_{\Omega}\vartheta(\cdot,\tilde \tau)\dd x=
\int_{\Omega}\vartheta^0\dd x
+\int_{\calQ_{\tilde \tau}}\Le\ee(\dot u){:}\ee(\dot u)\dd x\dd t
+\int_{\calQ_{\tilde \tau}}\Me\dot z{:}\dot z\dd x\dd t\\&
+\int_{\calQ_{\tilde \tau}}\theta(\alpha \tr(\ee(\dot u)){+} \D_z
H_2(z){:}
\dot z)\dd x\dd t+
\int_{\calQ_{\tilde \tau}}\Psi(\dot z)\dd x\dd t.
\end{aligned}
\end{equation*}
We add this last equality to \eqref{eq:55},  
and thanks to \eqref{eq:H1}, we obtain
\begin{equation*}
\begin{aligned}
& \tfrac12\int_{\Omega}\Ee(\ee(u(\cdot,\tilde \tau)){-}z(\cdot,\tilde
\tau)){:}(\ee(u(\cdot,\tilde \tau)){-}z(\cdot,\tilde \tau))\dd x
+\tfrac{\nu}2\norm[\L^2(\Omega)]{\nabla z(\cdot,\tilde \tau)}^2 
+c^{H_1} \norm[\L^2(\Omega)]{z(\cdot,\tilde\tau)}^2 \\&+
\int_{\Omega} \vartheta (\cdot, \tilde \tau) \dd x
\leq 
 C_0^{u,z}+ \int_{\Omega} \vartheta^0 \dd x + \tilde c^{H_1} \abs{\Omega} 
+ \int_{\calQ_{\tilde \tau}}\ell{\cdot}\dot u \dd x \dd t .
\end{aligned}
\end{equation*}
Clearly there exists $C_1>0$, depending only on $c^{\Ee}$, 
$\nu$ and $c^{H_1}$ such that
\begin{equation}
\label{eq:78}
C_1 \norm[\H^1(\Omega)]{u (\cdot, \tilde \tau)}^2 
+C_1 \norm[\H^1(\Omega)]{z(\cdot, \tilde \tau)}^2 
+ \int_{\Omega} \vartheta (\cdot, \tilde \tau) \dd x \le  C_0^{u,z} 
+ \norm[\L^1(\Omega)]{\vartheta^0} 
+ \tilde c^{H_1} \abs{\Omega} + \int_{\calQ_{\tilde \tau}} 
\ell{\cdot}\dot u \dd x \dd t .
\end{equation}
Since $\ell\in\H^1(0,T;\L^2(\Omega))$, we may integrate by parts 
the last term of \eqref{eq:78}, we get
\begin{equation*}
\begin{aligned} 
& \tfrac{C_1}{2} \norm[\H^1(\Omega)]{u (\cdot, \tilde \tau)}^2 
+C_1 \norm[\H^1(\Omega)]{z(\cdot, \tilde \tau)}^2 
+ \int_{\Omega} \vartheta (\cdot, \tilde \tau) \dd x \le 
C_0^{u,z} + \norm[\L^1(\Omega)]{\vartheta^0} 
+c^{H_1}\abs{\Omega} \\
&+\norm[{\C^0([0,T]; \L^2(\Omega))}]{\ell} \norm[\L^2(\Omega)]{u^0}
+\tfrac{1}{2 C_1} \norm[{\C^0([0,T]; \L^2(\Omega))}]{\ell}^2 
+ \tfrac12 \norm[\L^2(0,T; \L^2(\Omega))]{\dot \ell}^2  
+\tfrac12 \int_0^{\tilde \tau} \norm[\L^2(\Omega)]{u}^2 \dd t,
\end{aligned}
\end{equation*}
which allows us to conclude with  Gr\"onwall's lemma 
since $\vartheta (x,t) \geq
0$ almost everywhere on $\calQ_{\tau}$.
\eproof

Note that \eqref{eqlp:1} enable us also to obtain a global estimate for the
temperature. More precisely, under the assumptions of Proposition 
\ref{thm:existence-result-global}, we have 
\begin{equation}
\label{eq:29}
\norm[\L^{\infty}(0, \tau; \L^{\beta_1} (\Omega))]{\theta} 
\le \bigl(\tfrac{\beta_1}{c^c} \tilde C \bigr)^{\frac1{\beta_1}},
\end{equation}
for any solution $(u,z, \theta)$ of problem 
\eqref{eq:ent_eq}--\eqref{eq:boun_cond} defined on $[0, \tau]$ 
with $\tau \in (0,T]$.

We assume that $\beta_1 \ge 4$. We define
\begin{equation*}
\bar R^{\theta}\eqldef T^{\frac1q} \abs{\Omega}^{\frac{\beta_1{-}4}{4\beta_1}} 
\bigl(\tfrac{\beta_1}{c^c} \tilde C\bigr)^{\frac1{\beta_1}},
\end{equation*}
and using the notations of Section \ref{sec:existence-result}, we define 
\begin{equation*}
\bar R^f\eqldef R^f(R^0,\bar R^{\theta},\norm[{\C^0([0,T];
\L^2(\Omega))}]{\ell}), 
\quad 
\bar R^{\vartheta}_{\infty}\eqldef C\bigl(T^{\frac{q{-}8}{2q}} \bar R^f{+} 
\norm[\L^2(\Omega)]{\vartheta^0}\bigr), 
\quad
\bar R^{\vartheta}\eqldef
T^{\frac{1}{\bar q}} \bar R^{\vartheta}_{\infty} +1.
\end{equation*}
Then, the results of Section \ref{sec:existence-result} 
allow us to infer that there exists $\tau \in (0, T]$ 
such that $\phi_{\tau}$ admits 
a fixed point in $\bar B_{\L^{\bar q}(0, \tau; \L^2(\Omega))} 
(0,\bar R^{\vartheta})$. 
Let us define 
\begin{equation*}
\bar\tau\eqldef\sup\bigl{\{}\tau\in(0, T]:\ \phi_{\tau} 
\textrm{ admits a fixed point in }
\bar B_{\L^{\bar q }(0,\tau;\L^2(\Omega))}(0,\bar R^{\vartheta})\bigr{\}}
\in (0,T].
\end{equation*}
It is clear that problem \eqref{eq:ent_eq1}--\eqref{eq:init_cond1} admits 
a global solution  if and only if $\bar \tau = T$. This identity is
established below by a contradiction argument. To do so, we assume that 
$\bar\tau\in(0,T)$ and we choose $\epsilon>0$ such that 
$ \bar \tau - 
\epsilon \in (0, \bar \tau)$. By definition 
of $\bar \tau$,  there exists $\tau \in (\bar\tau{-}\epsilon, \bar \tau]$ 
such that $\phi_{\tau}$ admits a fixed point 
$\vartheta = \phi_{\tau} (\vartheta)$ in  
$\bar B_{\L^{\bar q}(0, \tau; \L^2(\Omega))} (0, \bar R^{\vartheta})$, i.e., 
the problem \eqref{eq:ent_eq1}--\eqref{eq:init_cond1} admits a solution 
$(u,z, \vartheta)$ defined on $[0,\tau]$. 
We infer from the results of Section \ref{sec:global-existence-result} 
that $\norm[\L^{\infty}(0,\tau; \L^1(\Omega))]{\vartheta} \le 
\tilde C$ and $\norm[\L^q(0, \tau; \L^4(\Omega))]{\theta=\zeta(\vartheta)} 
\le \bar R^{\theta}$. Then, the definition of $\bar R^{\vartheta}_{\infty}$ 
and the results of Section \ref{sec:enthalpy-equation} imply that 
$\vartheta= \phi_{\tau} (\vartheta) \in \L^{\infty} (0, \tau;\L^2(\Omega))$ 
with 
$\norm[\L^{\infty}(0, \tau; \L^2(\Omega))]{\vartheta} \le \bar
R^{\vartheta}_{\infty}$. 

Now let $\tilde\tau\in(0,T{-}\bar \tau]$ and $\tilde R^{\vartheta}\eqldef 
( (\bar R^{\vartheta})^{\bar q}{-} \bar \tau 
(\bar R^{\vartheta}_{\infty})^{\bar q})^{\frac1{\bar q}} >0$. 
For any $\tilde \vartheta \in \bar B_{\L^{\bar q}(\tau, \tau{+}\tilde \tau;
\L^2(\Omega))} (0, \tilde R^{\vartheta})$, we define $\tilde 
\vartheta_{\text{ext}}$ as
follows $\tilde \vartheta_{\text{ext}}\eqldef\vartheta$ on $[0,\tau]$, 
$\tilde \vartheta_{\text{ext}}\eqldef\tilde \vartheta$ on 
$(\tau,\tau{+}\tilde \tau]$ 
and $\tilde \vartheta_{\text{ext}}\eqldef 0$ on $(\tau{+}\tilde \tau, T]$. 
Clearly, we have
\begin{equation*}
\norm[\L^{\bar q}(0,T; \L^2(\Omega))]{\tilde \vartheta_{\text{ext}}}^{\bar q} 
= \norm[\L^{\bar q}(0,\tau; \L^2(\Omega))]{\vartheta}^{\bar q} + 
\norm[\L^{\bar q}(\tau,\tau{+}\tilde \tau;\L^2(\Omega))]{\tilde
  \vartheta}^{\bar q} 
\le \tau (\bar R^{\vartheta}_{\infty})^{\bar q} + 
(\tilde R^{\vartheta})^{\bar q} 
\le  (\bar R^{\vartheta})^{\bar q},
\end{equation*}
and the mapping $\tilde \vartheta \mapsto \tilde 
\vartheta_{\text{ext}}$ is a contraction 
on $\L^{\bar q}(\tau,\tau{+}\tilde \tau; \L^2(\Omega))$. Let $\tilde \theta 
= \zeta (\tilde \vartheta_{\text{ext}})$. By definition of $\zeta$, we have 
$\tilde \theta = \zeta( \vartheta) = \theta$ on $[0, \tau]$, 
$\tilde \theta = \zeta (\tilde \vartheta)$ on $(\tau, \tau + \tilde \tau]$ 
and $\tilde \theta = \zeta (0) =0 $ on $(\tau + \tilde \tau, T]$. 
Hence $\tilde \theta \in \L^q(0,T; \L^4 (\Omega))$ and 
\begin{equation*}
\norm[\L^q(0,T; \L^4 (\Omega))]{\tilde \theta}^q \le (\bar R^{\theta})^q + 
\int_{\tau}^{\tau + \tilde \tau} \norm[\L^4 (\Omega)]{\zeta( \tilde
  \vartheta)}^q  \dd t \le (\bar R^{\theta})^q + 
\bigl(\tfrac{\beta_1}{c^c}\bigr)^{\frac{q}{\beta_1}}
\abs{\Omega}^{\frac{\beta_1 -2}{4\beta_1} q} \norm[\L^{\bar q}(\tau,\tau +
\tilde \tau ; \L^2 (\Omega))]{\tilde \vartheta}^{\bar q}  
\leq  (\tilde R^{\theta})^q,
\end{equation*}
with  $\tilde R^{\theta}\eqldef\bigl((\bar R^{\theta})^q{+} 
\bigl( \tfrac{\beta_1}{c^c}\bigr)^{\frac{q}{\beta_1}} 
\abs{\Omega}^{\frac{\beta_1 {-}2}{4 \beta_1} q} 
(\tilde R^{\vartheta})^{\bar q}\bigr)^{\frac1{q}}$. 
By definition of
$\phi$, 
we get immediately that the restriction of 
$\phi(\tilde \vartheta_{\text{ext}})$ on $[0,
\tau]$ 
coincide with $\phi_{\tau}(\vartheta) = \vartheta$ and we define $\tilde
\phi_{\tilde \tau} (\tilde \vartheta)$ as the restriction of $\phi(\tilde
\vartheta_{\text{ext}})$ to 
$[\tau,\tau{+}\tilde \tau]$. Furthermore, with the estimates of Section 
\ref{sec:existence-result}, we have 
$\phi(\tilde \vartheta_{\text{ext}}) \in \L^{\infty}(0,T; \L^2(\Omega))$ and
\begin{equation*}
\norm[\L^{\infty} (0,T; \L^2(\Omega))]{\phi(\tilde \vartheta_{\text{ext}})} \le
 C\bigl( T^{\frac{q{-}8}{2q}}  R^f(R^0, \tilde R^{\theta}, 
\norm[{\C^0([0,T]; \L^2(\Omega)}]{\ell}){+}\norm[\L^2(\Omega)]{\vartheta^0}
\bigr).
\end{equation*}
It follows that there exists $\tilde \tau \in (0, T{-}\bar \tau]$, independent 
of $\tau$, such that $\tilde \phi_{\tilde \tau}$ admits a 
fixed point $\tilde \vartheta$ in 
$\bar B_{\L^{\bar q}(\tau, \tau + \tilde \tau; \L^2(\Omega))} (0, \tilde
R^{\vartheta})$. 
By construction of $\tilde \phi_{\tilde \tau}$, the restriction of 
$\phi(\tilde \vartheta_{\text{ext}})$ to $[0,\tau{+}\tilde \tau]$ 
is also a fixed point
of 
$\phi_{\tau{+}\tilde \tau}$ 
in $\bar B_{\L^{\bar q}(0,\tau{+}\tilde \tau; \L^2(\Omega))}(0,\bar R^{\vartheta})$. 
But we may choose $\epsilon\in (0,\bar\tau)$ such that 
$\tau+\tilde \tau > \bar \tau - \epsilon + \tilde \tau > 
\bar \tau$, which gives a contradiction with the
definition of $\bar \tau$.

Hence we can conclude that $\bar \tau=T$. Consequently, we deduce the
following theorem:

\begin{theorem}[Global existence result] 
\label{global}
Assume that \eqref{eq:Psi}, \eqref{eq:H}, \eqref{eq:47}, \eqref{eq:LM},
\eqref{eq:ell} and \eqref{eq:36} hold. Assume moreover that $\beta_1 \ge 4$, 
$u^0 \in\V^p_0 (\Omega)$, $z^0 \in \X_{q,p} (\Omega)$, $\vartheta^0 \in
\L^2(\Omega)$ 
such that \eqref{init_enthalpy} holds.  Then the problem 
\eqref{eq:ent_eq1}--\eqref{eq:init_cond1}
admits a global solution 
$(u,z, \vartheta)$ such that \(u\in\W^{1,q}(0,T;\V^p_0 (\Omega))\),
\(z\in\L^{\infty}(0,T;\H^1(\Omega) \cap \X_{q,p}(\Omega))\cap
\H^1(0,T;\L^2(\Omega))  \), 
\(\dot z, \Delta z \in \L^{q/2} (0,T;\L^p(\Omega))
\cap\L^{q}(0,T; \L^{p/2} (\Omega))\) and
\(\vartheta \in {\mathcal W} \). Moreover $\vartheta$ remains 
strictly positive and $(u,z, \theta = \zeta(\vartheta))$ is a 
solution of problem \eqref{eq:ent_eq}--\eqref{eq:boun_cond}  on $[0,T]$.
\end{theorem}

\begin{remark}
Let us assume furthermore that there exists $\tilde C_z^{H_2}>0$ and 
$p_1  \in (0,1)$ such that
\begin{equation*}
\forall z \in \Er^{3{\times}3}_{\emph{dev}} : 
\ \abs{\D_z H_2 (z)} \le \tilde C_z^{H_2} 
(1{+}\abs{z}^{p_1}).
\end{equation*}
Then we may obtain a global existence result for any 
$\beta_1 > \max \bigl(3,\frac{6}{5(1{-}p_1)}\bigr)$. Indeed, we can establish 
a global a priori estimate for the enthalpy gradient by using a technique 
inspired by Boccardo and Gallou\"et (\cite{BocGal89NEPM}), i.e., by 
choosing the test-function $\chi = 1 - \frac{1}{(1+ \vartheta)^{\xi}}$ 
for some $\xi >0$. By reproducing the same computations as in  
\cite[Prop.~4.2]{Roub10TRIV}, we may obtain an estimate of 
$\nabla \vartheta$ in $\L^r(0,\tau; \L^r(\Omega))$, independent of $\tau$, for
any $r \in \bigl[1, \frac{d{+}2}{d{+}1}\bigr)$ provided that 
$\beta_1 > \frac{2d}{(d{+}2)(1{-}p_1)} = \frac{6}{5(1{-}p_1)}$ since $d=3$. 
Then we consider $\alpha >1$ such that $2 \mu \alpha \le r$ and 
$\frac{1}{\alpha} \ge \mu\bigl(\frac{1}{r}{-}\frac{1}{d}\bigr)+1-\mu$ for 
some real number $\mu \in (0,1)$. Using Gagliardo-Nirenberg's 
inequality, we may deduce that there exists \(C_{\emph{GN}}>0\) such that
\begin{equation*}
\norm[\L^{2 \alpha} (0, \tau; \L^{\alpha} (\Omega))]{\vartheta}^{2 \alpha} 
\le C_{\emph{GN}} \norm[\L^{\infty}(0, \tau;
\L^1(\Omega))]{\vartheta}^{(1{-}\mu)2\alpha}  \int_0^{\tau} \bigl(
\norm[\L^1(\Omega)]{\vartheta (t, \cdot)}{+} 
\norm[\L^r(\Omega)]{\nabla \vartheta (t, \cdot)} \bigr)^{2\mu \alpha} \dd t,
\end{equation*}
and reminding the estimate of $\vartheta$ in 
$\L^{\infty}(0, \tau; \L^1(\Omega))$ obtained at 
Proposition \ref{thm:existence-result-global}, 
we infer that there exists a constant $C>0$, 
depending only on the data, such that
\begin{equation*}
\norm[\L^{2 \alpha} (0, \tau; \L^{\alpha} (\Omega))]{\vartheta} \le 
C \bigl(1{+}\norm[\L^{r} (0, \tau; \L^r(\Omega))]{\nabla \vartheta}\bigr).
\end{equation*}
The three conditions
\begin{equation*}
1 \le r  <\tfrac{d+2}{d+1} , \quad 2 \mu \alpha \le r , 
\quad \tfrac{1}{\alpha} \ge \mu \bigl(\tfrac{1}{r}{-} 
\tfrac{1}{d}\bigr) +1-\mu\text{ with } \ 0< \mu <1,
\end{equation*}
allow us to choose  $\mu= \frac{rd}{d{+}r(d{+}1)} \in (0,1)$
and thus $\alpha= \frac{r}{2 \mu}= \frac{1}{2} + 
\frac{r(d+1)}{2d} \in \bigl(\frac{7}{6}, \frac{4}{3}\bigr)$
(here $d=3$). It follows that, for any $\beta_1 > 
\max\bigl(3, \frac{6}{5(1{-}p_1)}\bigr)$, we will obtain a 
global estimate of $\vartheta$ in $\L^{2\alpha}(0,\tau; \L^{\alpha}(\Omega))$ 
for any $\alpha \in \bigl(\frac{7}{6}, \frac{4}{3}\bigr)$ and of 
$\theta =\zeta(\vartheta)$ in 
$\L^{2\beta_1 \alpha}(0,\tau; \L^{\beta_1 \alpha}(\Omega))$. 
Thus, with $\alpha$ such that $\beta_1  \alpha >4$, we may 
obtain a global existence result by the same contradiction 
argument as in Section \ref{sec:global-existence-result}.
\end{remark}

\paragraph*{Acknowledgments.} 
{\hspace{-0.4cm}}
The authors warmly thank to K.~Gr\"oger and A.~Mielke for many fruitful
discussions and their encouragements.
A.P. is grateful to J.~Rehberg for stimulating discussions.
A.P. was supported  by the Deutsche Forschungsgemeinschaft 
through the projet C18 ``Analysis and numerics 
of multidimensional models for elastic phase transformation
in a shape-memory alloys'' of the Research Center {\sc Matheon}.
Moreover, L.P. wishes to thank for the kind hospitality of WIAS.

\renewcommand{\arraystretch}{0.91}\small 


\begin{thebibliography}{AA00}\itemsep0.1em

\bibitem[AuP02]{AurPet02IACR}
{\scshape F.~Auricchio {\upshape and} L.~Petrini}.
\newblock Improvements and algorithmical considerations on a recent
  three-dimensional model describing stress-induced solid phase
  transformations.
\newblock {\em Int. J. Numer. Meth. Engng.}, 55, 1255--1284, 2002.

\bibitem[AuP04]{AurPet04STPT}
{\scshape F.~Auricchio {\upshape and} L.~Petrini}.
\newblock {A three-dimensional model describing stress-temperature induced
  solid phase transformations: thermomechanical coupling and hybrid composite
  applications.}
\newblock {\em Int. J. Numer. Methods Eng.}, 61(5), 716--737, 2004.

\bibitem[BaR08]{BarRou08TPSS}
{\scshape S.~Bartels {\upshape and} T.~Roub{\'{\i}}{\v{c}}ek}.
\newblock Thermoviscoplasticity at small strains.
\newblock {\em ZAMM Z. Angew. Math. Mech.}, 88(9), 735--754, 2008.

\bibitem[BoG89]{BocGal89NEPM}
{\scshape L.~Boccardo {\upshape and} T.~Gallou{\"e}t}.
\newblock Nonlinear elliptic and parabolic equations involving measure data.
\newblock {\em J. Funct. Anal.}, 87(1), 149--169, 1989.

\bibitem[Bre73]{Brez73OMMS}
{\scshape H.~Brezis}.
\newblock {\em Op{\'e}rateurs maximaux monotones et semi-groupes de
  contractions dans les espaces de {H}ilbert}.
\newblock North-Holland Publishing Co., Amsterdam, 1973.
\newblock North-Holland Mathematics Studies, No. 5. Notas de Matem\'atica (50).

\bibitem[Bre83]{Brez83AFTA}
{\scshape H.~Brezis}.
\newblock {\em Analyse fonctionnelle}.
\newblock Collection Math\'ematiques Appliqu\'ees pour la Ma\^{i}trise.
  [Collection of Applied Mathematics for the Master's Degree]. Masson, Paris,
  1983.
\newblock Th\'eorie et applications. [Theory and applications].

\bibitem[Car90]{Car90}
{\scshape H.~Cartan}.
\newblock {\em Cours de calcul diff\'erentiel}.
\newblock Hermann, Paris, 1990.

\bibitem[CoV90]{ColVis90CDNE}
{\scshape P.~Colli {\upshape and} A.~Visintin}.
\newblock On a class of doubly nonlinear evolution equations.
\newblock {\em Comm. Partial Differential Equations}, 15(5), 737--756, 1990.

\bibitem[Dor93]{Dore91RADE}
{\scshape G.~Dore}.
\newblock {$\L^p$} regularity for abstract differential equations.
\newblock In {\em Functional analysis and related topics, 1991 ({K}yoto)},
  volume 1540 of {\em Lecture Notes in Math.}, pages 25--38. Springer, Berlin,
  1993.

\bibitem[DuL76]{DuvLio76IMP}
{\scshape G.~Duvaut {\upshape and} J.-L.~Lions}.
\newblock {\em Inequalities in mechanics and physics}.
\newblock Springer-Verlag, Berlin, 1976.
\newblock Translated from the French by C. W. John, Grundlehren der
  Mathematischen Wissenschaften, 219.

\bibitem[EfM06]{EfeMie06RILS}
{\scshape M.~Efendiev {\upshape and} A.~Mielke}.
\newblock On the rate--independent limit of systems with dry friction and small
  viscosity.
\newblock {\em J. Convex Analysis}, 13(1), 151--167, 2006.

\bibitem[FrM06]{FraMie06ERCR}
{\scshape G.~Francfort {\upshape and} A.~Mielke}.
\newblock Existence results for a class of rate-independent material models
  with nonconvex elastic energies.
\newblock {\em J. Reine Angew. Math.}, 595, 55--91, 2006.

\bibitem[GHH07]{GoHaHe07UBFE}
{\scshape S.~Govindjee, K.~Hackl, {\upshape and} R.~Heinen}.
\newblock An upper bound to the free energy of mixing by twin-compatible
  lamination for {$n$}-variant martensitic phase transformations.
\newblock {\em Contin. Mech. Thermodyn.}, 18(7-8), 443--453, 2007.

\bibitem[GMH02]{GoMiHa02FEMV}
{\scshape S.~Govindjee, A.~Mielke, {\upshape and} G.~J.~Hall}.
\newblock The free--energy of mixing for $n$--variant martensitic phase
  transformations using quasi--convex analysis.
\newblock {\em J. Mech. Physics Solids}, 50, 1897--1922, 2002.
\newblock Erratum and Correct Reprinting: 51(4) 2003, pp.~763 \& I-XXVI.

\bibitem[HaG02]{HalGov02ARFE}
{\scshape G.~Hall {\upshape and} S.~Govindjee}.
\newblock Application of the relaxed free energy of mixing to problems in shape
  memory alloy simulation.
\newblock {\em J. Intelligent Material Systems Structures}, 13, 773--782, 2002.

\bibitem[HaN75]{HalNgu75SMSG}
{\scshape B.~Halphen {\upshape and} Q.~S.~Nguyen}.
\newblock Sur les mat\'eriaux standards g\'en\'eralis\'es.
\newblock {\em J. M\'ecanique}, 14, 39--63, 1975.

\bibitem[HiR08]{HieReh08QPMB}
{\scshape M.~Hieber {\upshape and} J.~Rehberg}.
\newblock Quasilinear parabolic systems with mixed boundary conditions on
  nonsmooth domains.
\newblock {\em SIAM J. Math. Anal.}, 40(1), 292--305, 2008.

\bibitem[KoO88]{KonOle88BVPS}
{\scshape V.~Kondrat'ev {\upshape and} O.~A.~Oleinik}.
\newblock Boundary-value problems for the system of elasticity theory in
  unbounded domains. korn's inequailities.
\newblock {\em Russian Math. Surveys}, 43(5), 65--119, 1988.

\bibitem[Mie00]{Miel00EMFM}
{\scshape A.~Mielke}.
\newblock Estimates on the mixture function for multiphase problems in
  elasticity.
\newblock In A.-M.~S\"andig, W.~Schiehlen, {\upshape and} W.~Wendland, editors,
  {\em Multifield Problems}, pages 96--103, Berlin, 2000. Springer--Verlag.

\bibitem[Mie05]{Miel05ERIS}
{\scshape A.~Mielke}.
\newblock Evolution in rate-independent systems ({C}h.~6).
\newblock In C.~Dafermos {\upshape and} E.~Feireisl, editors, {\em Handbook of
  Differential Equations, Evolutionary Equations, vol.~2}, pages 461--559.
  Elsevier B.V., Amsterdam, 2005.

\bibitem[Mie07]{Miel07MTIP}
{\scshape A.~Mielke}.
\newblock A model for temperature-induced phase transformations in
  finite-strain elasticity.
\newblock {\em IMA J. Applied Math.}, 72, 644--658, 2007.

\bibitem[MiP07]{MiePet07TDPT}
{\scshape A.~Mielke {\upshape and} A.~Petrov}.
\newblock Thermally driven phase transformation in shape-memory alloys.
\newblock {\em Gakk\=otosho (Adv. Math. Sci. Appl.)}, 17, 667--685, 2007.

\bibitem[MiR06]{MieRou06RIDP}
{\scshape A.~Mielke {\upshape and} T.~Roub{\'\i}{\v{c}}ek}.
\newblock Rate-independent damage processes in nonlinear elasticity.
\newblock {\em M$^3\!$AS Math. Models Methods Appl. Sci.}, 16, 177--209, 2006.

\bibitem[MiR07]{MieRos07EURC}
{\scshape A.~Mielke {\upshape and} R.~Rossi}.
\newblock Existence and uniqueness results for a class of rate-independent
  hysteresis problems.
\newblock {\em M$^3$AS Math. Models Methods Appl. Sci.}, 17, 81--123, 2007.

\bibitem[MiT99]{MieThe99MMRI}
{\scshape A.~Mielke {\upshape and} F.~Theil}.
\newblock A mathematical model for rate-independent phase transformations with
  hysteresis.
\newblock In H.-D.~Alber, R.~Balean, {\upshape and} R.~Farwig, editors, {\em
  Proceedings of the Workshop on ``Models of Continuum Mechanics in Analysis
  and Engineering''}, pages 117--129, Aachen, 1999. Shaker-Verlag.

\bibitem[MiT04]{MieThe04RIHM}
{\scshape A.~Mielke {\upshape and} F.~Theil}.
\newblock On rate--independent hysteresis models.
\newblock {\em Nonl. Diff. Eqns. Appl. (NoDEA)}, 11, 151--189, 2004.
\newblock (Accepted July 2001).

\bibitem[MPM08]{MiPeMa08CSKV}
{\scshape A.~Mielke, A.~Petrov, {\upshape and} J.~A.~C.~Martins}.
\newblock Convergence of solutions of kinetic variational inequalities in the
  rate-independent quasi-static limit.
\newblock {\em J. Math. Anal. Appl.}, 348, 1012--1020, 2008.

\bibitem[MRS08]{MiRoSt08GLRR}
{\scshape A.~Mielke, T.~Roub\'{\i}\v{c}ek, {\upshape and} U.~Stefanelli}.
\newblock {${\Gamma}$}-limits and relaxations for rate-independent evolutionary
  problems.
\newblock {\em Calc. Var. Part. Diff. Equ.}, 31, 387--416, 2008.

\bibitem[MTL02]{MiThLe02VFRI}
{\scshape A.~Mielke, F.~Theil, {\upshape and} V.~I.~Levitas}.
\newblock A variational formulation of rate--independent phase transformations
  using an extremum principle.
\newblock {\em Arch. Rational Mech. Anal.}, 162, 137--177, 2002.
\newblock (Essential Science Indicator: Emerging Research Front, August 2006).

\bibitem[PaP11]{PaoPet11TMPV}
{\scshape L.~Paoli {\upshape and} A.~Petrov}.
\newblock Thermodynamics of multiphase problems in viscoelasticity.
\newblock {\em To appear in GAMM-Mitteilungen}, 2011.

\bibitem[PrS01]{PruSch01SMRP}
{\scshape J.~Pr{\"u}ss {\upshape and} R.~Schnaubelt}.
\newblock Solvability and maximal regularity of parabolic evolution equations
  with coefficients continuous in time.
\newblock {\em J. Math. Anal. Appl.}, 256(2), 405--430, 2001.

\bibitem[RaT83]{RavTho83}
{\scshape P.~A.~Raviart {\upshape and} J.~M.~Thomas}.
\newblock {\em Introduction \`a l'analyse num\'erique des \'equations aux
  d\'eriv\'ees partielles}.
\newblock Collection Math\'ematiques Appliqu\'ees pour la Ma\^{i}trise.
  [Collection of Applied Mathematics for the Master's Degree]. Masson, Paris,
  1983.

\bibitem[Rou09a]{Roub09RISS}
{\scshape T.~Roub{\'{\i}}{\v{c}}ek}.
\newblock Rate-independent processes in viscous solids at small strains.
\newblock {\em Math. Methods Appl. Sci.}, 32(7), 825--862, 2009.

\bibitem[Rou09b]{Roub09TVEL}
{\scshape T.~Roub{\'{\i}}{\v{c}}ek}.
\newblock Thermo-visco-elasticity at small strains with {$\L^1$}-data.
\newblock {\em Quart. Appl. Math.}, 67(1), 47--71, 2009.

\bibitem[Rou10]{Roub10TRIV}
{\scshape T.~Roub{\'{\i}}{\v{c}}ek}.
\newblock Thermodynamics of rate-independent processes in viscous solids at
  small strains.
\newblock {\em SIAM J. Math. Anal.}, 42(1), 256--297, 2010.

\bibitem[Sim87]{Sim87CSLP}
{\scshape J.~Simon}.
\newblock Compact sets in the space {$\L^p(0,T;B)$}.
\newblock {\em Ann. Mat. Pura Applic.}, 146, 65--96, 1987.

\bibitem[SMZ98]{SoMaZo98TDMS}
{\scshape A.~Souza, E.~Mamiya, {\upshape and} N.~Zouain}.
\newblock Three-dimensional model for solids undergoing stress-induced phase
  transformations.
\newblock {\em Europ. J. Mech., A/Solids}, 17, 789--806, 1998.

\bibitem[Val88]{Vale88LTEU}
{\scshape T.~Valent}.
\newblock {\em Boundary value problems of finite elasticity}, volume~31 of {\em
  Springer Tracts in Natural Philosophy}.
\newblock Springer-Verlag, New York, 1988.
\newblock Local theorems on existence, uniqueness, and analytic dependence on
  data.

\end{thebibliography}

\end{document}